\documentclass[11pt]{amsart}
\usepackage{amsfonts}
\usepackage{amssymb, amsmath}
\usepackage[letterpaper, margin=0.85in]{geometry}
\usepackage{fancyhdr}

\begin{document}
\author[]{Ali Enayat}
\title[]{Tarskian truth theories over set theory}
\maketitle

\begin{abstract}
\noindent This work is focused on Tarskian truth theories over set theory,
i.e., extensions of $\mathsf{CT}^{-}[\mathsf{ZF}]$. The theory $\mathsf{CT}%
^{-}[\mathsf{ZF}]$ is obtained by augmenting $\mathsf{ZF}$ with finitely
many Tarski-style compositional axioms for the truth predicate $\mathsf{T}$.
We pay special attention to the theory $\mathsf{CT}_{\ast }[\mathsf{ZF}],$
which is obtained by strengthening $\mathsf{CT}^{-}[\mathsf{ZF}]$ with the
sentence asserting \textquotedblleft all \textit{theorems} of $\mathsf{ZF}$
are true\textquotedblright . Our new results include the following:\medskip

\noindent \textbf{Theorem A.}~\textit{The following theories have the same
deductive closure}:\medskip

\noindent \textbf{(a) }$\mathsf{CT}_{\ast }[\mathsf{ZF}]$.\smallskip

\noindent \textbf{(b)} $\mathsf{CT}^{-}[\mathsf{ZF}]~+$ $\mathsf{Int}$-$%
\mathsf{Repl}$ $+$ $\mathsf{DC}_{\mathsf{out}}$.\smallskip

\noindent \textbf{(c)} $\mathsf{CT}^{-}[\mathsf{ZF}]~+$ $\mathsf{Int}$-$%
\mathsf{Ref}$.\smallskip

\noindent In the above, $\mathsf{Int}$-$\mathsf{Repl}$\textsf{\ }(internal
replacement) is the sentence asserting that all instances of the replacement
scheme are true; $\mathsf{DC}_{\mathsf{out}}$ is the sentence asserting that
if a finite disjunction is true, then at least one of its disjuncts is true;
and $\mathsf{Int}$-$\mathsf{Ref}$\textsf{\ }(internal reflection) is the
sentence asserting that all instances of the Montague reflection theorem are
true.\medskip

\noindent \textbf{Theorem B.}~\textit{The following are equivalent for a
sentence }$\varphi $ \textit{in the language of set theory}:\smallskip

\noindent \textbf{(a)} $\mathsf{CT}_{\ast }[\mathsf{ZFC}]\vdash \varphi .$%
\smallskip

\noindent \textbf{(b)} $\left\{ \mathsf{REF}^{n}(\mathsf{ZFC}):n\in \mathbb{N%
}\right\} \vdash \varphi .$\medskip

\noindent In the above, $\mathsf{REF}^{n}(\mathsf{ZFC})$ is
\textquotedblleft $n$-fold reflection over $\mathsf{ZFC}$.\textquotedblright
\medskip

\noindent \textbf{Theorem C.}~\textit{Assuming that }$\mathsf{ZF}$ \textit{%
has a model of the form} $\left( \mathrm{V}_{\!\alpha },\in \right) ,$ 
\textit{we have}:\smallskip

\noindent \textbf{(a) }\textit{The existence of a well-founded model of }$%
\mathsf{ZF}$\textit{\ is provable in }$\mathsf{CT}_{\ast }[\mathsf{ZF}]~+$ $%
\Delta _{0}$-$\mathsf{Sep}(\mathsf{T})$, \textit{but not in} $\mathsf{CT}%
_{\ast }[\mathsf{ZF}]$.

\noindent \textbf{(b) }\textit{The existence of a model of }$\mathsf{ZF}$%
\textit{\ of the form }$\left( \mathrm{V}_{\!\alpha },\in \right) $ \textit{%
is provable in }$\mathsf{CT}_{\ast }[\mathsf{ZF}]~+$ $\Delta _{0}$-$\mathsf{%
Sep}(\mathsf{T})+\Delta _{0}$-$\mathsf{Coll}(\mathsf{T})$, \textit{but not in%
} $\mathsf{CT}_{\ast }[\mathsf{ZF}]~+$ $\Delta _{0}$-$\mathsf{Sep}(\mathsf{T}%
)$. \medskip

\noindent In the above, $\Delta _{0}$-$\mathsf{Sep}(\mathsf{T})$ is the
separation scheme in the extended language for $\Delta _{0}$-formulae, and $%
\Delta _{0}$-$\mathsf{Coll}(\mathsf{T})$ is the collection scheme in the
extended language for $\Delta _{0}$-formulae.\medskip

\noindent \textbf{Theorem D}.~$\mathsf{CT}^{-}[\mathsf{ZF}]+\mathsf{Coll(T)}$
\textit{is conservative over} $\mathsf{ZF}$.\medskip
\end{abstract}

\begin{center}
\bigskip \bigskip

\textbf{1.~INTRODUCTION}\bigskip
\end{center}

Let $\mathsf{CT}^{-}$ be the finite list of axioms stipulating that the
truth predicate $\mathsf{T}$ satisfies Tarski's compositional clauses for
all formulae in the language of set theory, and let $\mathsf{B}$ (base
theory) be an extension of $\mathsf{KP}$ (Kripke-Platek set theory).%
\footnote{%
There are other fragments of $\mathsf{ZF}$ that can be used for this
purpose. See the beginning paragraph of Section 3.2.} We study theories of
the form $\mathsf{CT}^{-}[\mathsf{B}]+\Gamma (\mathsf{T})$, where $\mathsf{CT%
}^{-}[\mathsf{B]:=CT}^{-}+\mathsf{B},$ and $\Gamma (\mathsf{T})$ puts
further `good behavior' demands on the truth predicate $\mathsf{T}$. Our
focus will be on the special case of $\mathsf{B}=\mathsf{ZF}$
(Zermelo-Fraenkel set theory), but along the way we will encounter some
interesting results concerning $\mathsf{B}=\mathsf{KP}$ (see Section 6). The
arithmetical counterpart of our study boasts a vast literature, but in
comparison, the list of publications probing truth theories over set theory
is meager.\footnote{%
For truth theories over $\mathsf{PA}$ (Peano Arithmetic), see the monographs
by Halbach \cite{Halbach-book} and Cie\'{s}li\'{n}ski \cite{Cieslinski-book}%
, which provide technical minutiae, as well as philosophical motivations. An
updated overview of the subject is presented in the encyclopedia entry \cite%
{Halbach-Stanford} by Halbach, Leigh, and \L e\l yk.} \medskip

The formal relationship between truth and set theory was first revealed by
Tarski's celebrated theorems on definability/undefinability of truth (as
reviewed in Section 3 of this paper). Other classic results in the subject
are Montague's Reflection Theorem (Theorem 3.2.11), and Levy's Partial
Definability of Truth Theorem (Theorem 3.2.9). In retrospect, the first
major result in \textit{axiomatic} theories of truth over set theory was
obtained by Montague and Vaught \cite{Montague-Vaught}, who proved a strong
version of the reflection theorem within $\mathsf{CT}[\mathsf{ZF}],$ where $%
\mathsf{CT}[\mathsf{ZF}]$ is the result of strengthening $\mathsf{CT}^{-}[%
\mathsf{ZF}]$ by all instances of the replacement scheme in which $\mathsf{T}
$ can be mentioned (see Definition 4.16 and Theorem 4.17). Another milestone
in this subject is Krajewski's contribution \cite{Krajewski}, which includes
the proof of conservativity of $\mathsf{CT}^{-}[\mathsf{ZF}]$ over $\mathsf{%
ZF.}$ In the same paper, Krajewski introduced the key model-theoretic notion
of a \textit{satisfaction class, }which provides a potent framework for
establishing proof-theoretic results about axiomatic truth theories with the
help of model-theoretic techniques.\textrm{\ }\medskip

Somewhat more recently, the joint work of the author with Visser in \cite%
{Ali+Albert-long} and \cite{Ali+albert-short} introduced a robust
model-theoretic method, which has since come to be referred to as the
`EV-method', for building various kinds of full satisfaction classes in
order to provide a unified approach for exploring issues related to
conservativity and interpretability in the context of arbitrary base
theories.\footnote{\cite{Ali+Albert-long} was privately circulated among the
cognoscenti as a working paper, many results of which by now have appeared
in print.} A corollary of one of the general results in \cite%
{Ali+Albert-long} is that the theory $\mathsf{CT}^{-}[\mathsf{ZF}]$ +
\textquotedblleft the \textit{axioms} of $\mathsf{ZF}$ are
true\textquotedblright\ is conservative over $\mathsf{ZF}$ (see Theorem 4.6
and Corollary 4.7). On the other hand, as noted in \cite{Ali+Albert-long},
Krajewski's aforementioned conservativity result can be refined to the
conservativity of $\mathsf{CT}^{-}[\mathsf{ZF}]+\mathsf{Sep}(\mathsf{T})$
over $\mathsf{ZF}$, where $\mathsf{Sep}(\mathsf{T})$ is the natural
extension of the separation scheme to formulae that mention the truth
predicate. This theory has the feature that it proves that $\mathsf{T}$ is
closed under first order deductions (see Corollary 4.5).\footnote{%
In sharp contrast, it is known that the consistency of $\mathsf{PA}$ is
provable in the theory obtained by adding \textquotedblleft $\mathsf{T}$ is
closed under first order deductions\textquotedblright\ to $\mathsf{CT}^{-}[%
\mathsf{PA}]$.}\medskip

A deep analysis of various theories of truth -- both of the typed and
untyped varieties -- over set theoretical base theories was systematically
carried out by Fujimoto \cite{Fujimoto-APAL}, who uncovered a vibrant link
between truth theories and certain canonical class theories. The results in 
\cite{Fujimoto-APAL} are dominantly focused on `strong' theories of truth
that extend $\mathsf{CT}[\mathsf{ZF}]$. It should be noted that the
aforementioned conservativity of $\mathsf{CT}^{-}[\mathsf{ZF}]+\mathsf{Sep}(%
\mathsf{T})$ over $\mathsf{ZF}$ was independently established in \cite%
{Fujimoto-APAL}. Another significant contribution to the subject is by
Robert Van Wesep \cite{Van Wesep}, who investigated truth theories over the
universe of sets in the framework of the G\"{o}del-Bernays theory of classes.%
\footnote{%
This topic has also been studied in \cite{Ali-Mostowski-Bridge}, as well as
Section 7 of this paper. Note that \cite{Ali-Mostowski-Bridge} includes the
arithmetical inspirations for many of the set-theoretical results here.}
Truth theories over set theories have also made an appearance in
philosophical contexts; see, e.g., Fujimoto's \cite{Fujimoto-deflation}, the
joint work of Horsten, Luo, and Roberts \cite{Horsten et al}, and Heck's 
\cite{Heck}. \medskip

The results reported here address questions that naturally arise from two
sources: (a) the existing rather limited body of knowledge concerning
subtheories of $\mathsf{CT}[\mathsf{ZF}]$, and (b) the extensive results
concerning truth theories over arithmetical theories. The paper is organized
as follows:

\begin{itemize}
\item Sections 2, 3, and 4 contain preliminary definitions, conventions, and
relevant results of both basic and advanced variety.

\item The novel portion of the paper begins in Section 5, in which the
aforementioned full reflection theorem of Montague and Vaught is shown to be
provable in $\mathsf{CT}_{0}[\mathsf{ZF}]$, an intermediate system between $%
\mathsf{CT}^{-}[\mathsf{ZF}]$ and $\mathsf{CT}[\mathsf{ZF}].$ In particular, 
$\mathsf{CT}_{0}[\mathsf{ZF}]$ can prove that $\mathsf{ZF}$ has models of
the form $\left( \mathrm{V\!}_{\alpha },\in \right) $.

\item In Section 6, we study $\mathsf{CT}_{\ast }[\mathsf{ZF}]$, which is
the result of adding the sentence \textquotedblleft all \textit{theorems} of 
$\mathsf{ZF}$ are true" to $\mathsf{CT}^{-}[\mathsf{ZF}]$. The various
results of Section 6 show that $\mathsf{CT}_{\ast }[\mathsf{ZF}]$ can be
argued to be the set-theoretical analogue of the canonical theory known as $%
\mathsf{CT}_{0}[\mathsf{PA}]$ (as shown later in Section 9, $\mathsf{CT}%
_{\ast }[\mathsf{ZF}]$ is much weaker than $\mathsf{CT}_{0}[\mathsf{ZF}])$.

\item Section 7 explains the close relationship between $\mathsf{CT}_{\ast }[%
\mathsf{ZF}]$ and certain extensions of $\mathsf{GB}$ (G\"{o}del-Bernays
theory of classes).

\item Section 8 presents a natural axiomatization of the purely
set-theoretical consequences of $\mathsf{CT}_{\ast }[\mathsf{ZFC}]$, as in
Theorem B of the abstract.

\item In Section 9 we study certain natural truth theories intermediate
between $\mathsf{CT}_{\ast }[\mathsf{ZF}]$ and $\mathsf{CT}_{0}[\mathsf{ZF}%
]. $ For example, we show that even though $\mathsf{CT}_{\ast }[\mathsf{ZF}]$
proves the existence of a model of $\mathsf{ZF}$, the existence of an $%
\omega $-\textit{model} of $\mathsf{ZF}$ is unprovable in $\mathsf{CT}_{\ast
}[\mathsf{ZF}]$. In contrast, the stronger theory $\mathsf{CT}_{\ast }[%
\mathsf{ZF}]+\Delta _{0}$-$\mathsf{Sep(T)}$ is shown to prove the existence
of \textit{well-founded} models (and \textit{a fortiori}, $\omega $-models)
of $\mathsf{ZF}$, but it is incapable of proving that $\mathsf{ZF}$ has a
model of the form $\left( \mathrm{V\!}_{\alpha },\in \right) $.

\item Section 10 is devoted to the proof of conservativity of $\mathsf{CT}%
^{-}[\mathsf{ZF}]+\mathsf{Coll}(\mathsf{T})$ over $\mathsf{ZF}$, where $%
\mathsf{Coll}(\mathsf{T})$ is the natural extension of the collection scheme
to formulae that mention the truth predicate.

\item Section 11 presents some open questions that arise from this
work.\medskip
\end{itemize}

\noindent \textbf{Acknowledgments.}~The results presented here were inspired
by ideas and techniques that I have learnt over the past decade and a half
from many colleagues, including Bartosz Wcis\l o, Albert Visser, Jim
Schmerl, Fedor Pakhomov, Zach McKenzie, Adrian Mathias, Mateusz \L e\l yk,
Graham Leigh, Roman Kossak, Kentaro Fujimoto, Riki Heck, Volker Halbach,
Cezary Cie\'{s}li\'{n}ski, and Lev Beklemishev (in reverse alphabetical
order of last names).\footnote{%
In preparing this version, I benefited from helpful feedback on previous
versions from Volker Halbach, Zach McKenzie, Mateusz \L e\l yk, and Maciej G%
\l owacki.} \bigskip 

\bigskip

\begin{center}
\textbf{2.~PRELIMINARIES}\bigskip
\end{center}

\noindent \textbf{2.1.}~\textbf{Definitions and Basic Facts }(Languages and
set theories).~In what follows $\mathcal{L}\supseteq \mathcal{L}_{\mathrm{set%
}}:=\{=,\in \}.$\medskip

\begin{enumerate}
\item[\textbf{(a)}] We treat a theory as a \textit{set of axioms}, thus we
do not equate a theory with its deductive closure. \medskip

\item[\textbf{(b)}] The $\mathcal{L}$-\textit{induction scheme}, denoted $%
\mathsf{Ind}(\mathcal{L})$, consists of sentences of the form $\forall v\,%
\mathrm{Ind}_{\varphi (v,x)}$, where $\varphi (v,x)$ is an $\mathcal{L}$%
-formula, and%
\begin{equation*}
\mathrm{Ind}_{\varphi (v,x)}:=\left[ \varphi (v,0)\wedge \forall x\in \omega
\,\left( \varphi (v,x)\rightarrow \varphi (v,x+1)\right) \right] \rightarrow
\forall x\in \omega \,\varphi (v,x),
\end{equation*}%
Note that the parameter $v$ here ranges over the entire universe of
discourse, and is not limited to $\omega .$\footnote{%
Using a pairing function, one can deduce the more general versions of the
schemes considered here, in which the single parameter $v$ is replaced by a
finite tuple $\vec{v}$ of parameters.}\medskip

\item[\textbf{(c)}] The $\mathcal{L}$-\textit{separation scheme}, denoted $%
\mathsf{Sep}(\mathcal{L})$, consists of sentences of the form $\forall v\,%
\mathrm{Sep}_{\varphi (v,x)}$, where $\varphi (v,x)$ is an $\mathcal{L}$%
-formula, and%
\begin{equation*}
\mathrm{Sep}_{\varphi (v,x)}:=\forall a\,\exists b\,\left[ \forall x\,\left(
x\in b\longleftrightarrow \left( x\in a\wedge \varphi (v,x)\right) \right) %
\right] .
\end{equation*}

\item[\textbf{(d)}] The $\mathcal{L}$-\textit{replacement scheme}, denoted $%
\mathsf{Repl}(\mathcal{L})$, consists of sentences of the form $\forall v\ 
\mathrm{Repl}_{\varphi (v,x,y)}$, where $\varphi (v,x,y)$ is an $\mathcal{L}$%
-formula, and%
\begin{equation*}
\mathrm{Repl}_{\varphi (v,x,y)}:=\forall a\,\exists b\,\left[ \left( \forall
x\in a\,\exists !y\,\varphi (v,x,y)\right) \rightarrow \forall y\left( y\in
b\leftrightarrow \exists x\in a\,\varphi (v,x,y)\right) \right] .
\end{equation*}

\item[\textbf{(e)}] The $\mathcal{L}$-\textit{collection scheme}, denoted $%
\mathsf{Coll}(\mathcal{L})$, consists of sentences of the form $\forall v\ 
\mathrm{Coll}_{\varphi (v,x,y)}$, where $\varphi (v,x,y)$ is an $\mathcal{L}$%
-formula, and%
\begin{equation*}
\mathrm{Coll}_{\varphi (v,x,y)}:=\forall a\,\exists b\,\left[ \left( \forall
x\in a\text{ }\exists y\text{\ }\varphi (v,x,y)\right) \rightarrow \forall
x\in a\text{ }\exists y\in b\,\varphi (v,x,y)\right] .
\end{equation*}

\item[\textbf{(f)}] Given a class $\Phi $ of $\mathcal{L}$-formulae, we can
define the corresponding partial schemes. In particular, $\Phi $-$\mathsf{Ind%
}=\left\{ \forall v\,\mathrm{Ind}_{\varphi (v,x)}:\varphi \in \Phi \right\} $%
, $\Phi $-$\mathsf{Sep}=\left\{ \forall v\,\mathrm{Sep}_{\varphi
(v,x)}:\varphi \in \Phi \right\} $, etc. \medskip

\item[\textbf{(g)}] When $\mathsf{X}$\ is a new predicate, we write $%
\mathcal{L}_{\mathrm{set}}(\mathsf{X})$ instead of $\mathcal{L}_{\mathrm{set}%
}\cup \{\mathsf{X}\}$. For $\mathcal{L}=\mathcal{L}_{\mathrm{set}}(\mathsf{X}%
)$, we sometimes write $\mathsf{Sep}(\mathsf{X})$, $\mathsf{Coll}(\mathsf{X}%
) $, etc. instead of $\mathsf{Sep}(\mathcal{L})$, $\mathsf{Coll}(\mathcal{L})
$, etc. (respectively). \medskip

\item[\textbf{(h)}] Our axiomatization of $\mathsf{ZF}$ is as usual, except
that instead of including the scheme of replacement among the axioms of $%
\mathsf{ZF}$, we include the schemes of separation and collection, e.g., as
in \cite[Appendix A]{Chang-Keisler}. Thus each instance of the replacement
scheme is a \textit{theorem} of our $\mathsf{ZF}$. Given $\mathcal{L}%
\supseteq \mathcal{L}_{\mathrm{set}}$, we construe $\mathsf{ZF}(\mathcal{L})$
as the natural extension of $\mathsf{ZF}$ in which the schemes of separation
and collection are extended to $\mathcal{L}$-formulae, i.e., $\mathsf{ZF}(%
\mathcal{L})=\mathsf{ZF}+\mathsf{Sep}(\mathcal{L})+\mathsf{Coll}(\mathcal{L}%
).$\medskip

\item[\textbf{(i)}] For $n\in \mathbb{N}$, we employ the common notation ($%
\Sigma _{n},$ $\Pi _{n},$ $\Delta _{n}$) for the Levy hierarchy of $\mathcal{%
L}_{\mathrm{set}}$-formulae, as in the standard references in advanced set
theory such as Jech's monograph \cite{Jechbook-2003}. In particular, $\Delta
_{0}=\Sigma _{0}=\Pi _{0}$ corresponds to the collections of $\mathcal{L}_{%
\mathrm{set}}$-formulae all of whose quantifiers are bounded. For $n\in 
\mathbb{N}$ and $\mathcal{L}\supseteq \mathcal{L}_{\mathrm{set}},$ the Levy
hierarchy can be naturally extended to $\mathcal{L}$-formulae ($\Sigma _{n}(%
\mathcal{L}),$ $\Pi _{n}(\mathcal{L}),$ $\Delta _{n}(\mathcal{L})$), where $%
\Delta _{0}(\mathcal{L})$ is the smallest family of $\mathcal{L}$-formulae
that contains all atomic $\mathcal{L}$-formulae and is closed under Boolean
connectives and bounded quantification. \medskip 

\item[\textbf{(j)}] $\mathsf{KP}$ (Kripke-Platek) is the subtheory of $%
\mathsf{ZF}$ whose axioms consist of: the axioms \textrm{Extensionality}, 
\textrm{Pair}, and \textrm{Union}\textit{\ }together with the partial schemes%
\textit{\ }$\Pi _{1}$-$\mathsf{Found}$\footnote{$\mathsf{Found}$ consists of
the collection of sentences $\forall v\,\mathrm{Found}_{\varphi (v,x)}$,
where:%
\begin{equation*}
\mathrm{Found}_{\varphi (v,x)}=\left( \exists x\ \varphi (v,x)\right)
\rightarrow \left[ \exists x\ \left( \varphi (v,x)\wedge \forall y\in x\
\lnot \varphi (v,y)\right) \right] .
\end{equation*}%
Note that since we consider $\mathsf{Z}$ (Zermelo set theory) to include the
axiom of foundation (also known as regularity), each instance of $\mathsf{%
Found}$ is a theorem of $\mathsf{Z+TC}$, where $\mathsf{TC}$ (transitive
closure) states that every set is a subset of a transitive set. In
particular, each instance of $\mathsf{Found}$ is a theorem of $\mathsf{ZF}$.}%
, $\Delta _{0}$-$\mathsf{Sep}$, and $\Delta _{0}$-$\mathsf{Coll}$. Following
recent practice (initiated by Mathias \cite{Mathias-MacLane}), the
foundation scheme of $\mathsf{KP}$ is only limited to $\Pi _{1}$-formulae.
In contrast to Barwise's $\mathsf{KP}$ in \cite{Barwise}, which includes the
full scheme of foundation, our version of $\mathsf{KP}$ is finitely
axiomatizable. It is also worth pointing out that $\mathsf{KP}$ plus the
negation of the axiom of infinity is bi-interpretable with the fragment I$%
\Sigma _{1}$ of $\mathsf{PA}$; indeed the two theories can be shown to be
definitionally equivalent using the translations employed in \cite%
{Richard-Tin Lok}.\medskip\ 

\item[\textbf{(k)}] Zermelo set theory $\mathsf{Z}$ is obtained by removing
the scheme of collection from the axioms of $\mathsf{ZF}$.\footnote{%
Notice that our version of Zermelo set theory includes the axiom of
foundation (also known as regularity), but the foundation axiom is not
included in Zermelo set theory in some other sources.} $\mathsf{M}_{0}$ is a
subtheory of $\mathsf{Z}$ that is axiomatized by the axioms \textrm{%
Extensionality, Pairing, Union, Powerset}, together with the partial scheme $%
\Delta _{0}$-$\mathsf{Sep}$. Note that $\mathsf{M}_{0}$ is finitely
axiomatizable.\footnote{%
In the presence of the other axioms of $\mathsf{M}_{0}$, $\Delta _{0}$-$%
\mathsf{Sep}$ is well-known to be equivalent to the closure of the universe
under G\"{o}del-operations, see \cite[Theorem 13.4]{Jechbook-2003}.}\medskip
\end{enumerate}

\noindent \textbf{2.2.}~\textbf{Definitions and Basic Facts.~}(Model
theoretic concepts) We follow the convention of using $M$, $M^{\ast },$ $%
M_{0}$, etc.~to denote (respectively) the universes of $\mathcal{L}_{\mathrm{%
set}}$-structures $\mathcal{M}$, $\mathcal{M}^{\ast },$ $\mathcal{M}_{0},$
etc. We denote the membership relation of $\mathcal{M}$ by $\in ^{\mathcal{M}%
}$; thus an $\mathcal{L}_{\mathrm{set}}$-structure $\mathcal{M}$ is of the
form $(M,\in ^{\mathcal{M}})$. In what follows we make the blanket
assumption that $\mathcal{M}$, $\mathcal{N}$, etc.~are $\mathcal{L}$%
-structures, where $\mathcal{L}\supseteq \mathcal{L}_{\mathrm{set}}$.\medskip

\begin{enumerate}
\item[\textbf{(a)}] For an $\mathcal{L}$-formula $\varphi (\vec{x})$ with
free variables $\vec{x}=\left( x_{0},\cdot \cdot \cdot ,x_{k-1}\right) $ and 
\textit{suppressed parameters} from $M$, 
\begin{equation*}
\varphi ^{\mathcal{M}}:=\left\{ \vec{m}\in M^{k}:\mathcal{M\models \varphi }%
\left( m_{0},\cdot \cdot \cdot ,m_{k-1}\right) \right\} .
\end{equation*}%
A subset $D$ of $M^{k}$ is $\mathcal{M}$-\textit{definable} if it is of the
form $\varphi ^{\mathcal{M}}$ for some choice of $\varphi .$\medskip

\item[\textbf{(b)}] $\mathrm{Ord}^{\mathcal{M}}$ is the class of
\textquotedblleft ordinals\textquotedblright\ of $\mathcal{M}$, i.e.,%
\begin{equation*}
\mathrm{Ord}^{\mathcal{M}}:=\left\{ m\in M:\mathcal{M}\models \mathrm{Ord}%
(m)\right\} ,
\end{equation*}%
where $\mathrm{Ord}(x)$ expresses \textquotedblleft $x$ is transitive and is
well-ordered by $\in $\textquotedblright . \medskip

\item[\textbf{(c)}] We write $\mathbb{\omega }$ when referring to the set of
finite ordinals (i.e., natural numbers) of a given theory, and $\mathbb{%
\omega }^{\mathcal{M}}$ for the set of finite ordinals of a model $\mathcal{M%
}$ of set theory. We use $\mathbb{N}$ to refer to the set of natural numbers
in the real world, whose members we refer to as\textit{\ metatheoretic}%
\footnote{%
We will often use the expression `the real world' to refer to the metatheory.%
}\textit{\ natural numbers. }A model $\mathcal{M}$ of set theory is said to
be\textit{\ }$\omega $-\textit{standard} if $\left( \mathbb{\omega },\mathbb{%
\in }\right) ^{\mathcal{M}}\cong \left( \mathbb{N},\mathbb{<}\right) .$ We
identify the initial segment of $\left( \mathbb{\omega },\mathbb{\in }%
\right) ^{\mathcal{M}}$ that is isomorphic with $\left( \mathbb{N},\mathbb{<}%
\right) $ with $\mathbb{N}$. An element $k\in \omega ^{\mathcal{M}}$ is 
\textit{nonstandard} if $k\notin \mathbb{N}.$ \medskip

\item[\textbf{(d)}] For an ordinal $\alpha $, $\mathrm{V\!}_{\alpha }$ is
defined as usual as $\{x:\rho (x)<\alpha \}$, where $\rho (x)$ is the usual
rank function in set theory defined by $\rho (x)=\sup \{\rho (y)+1:y\in x\}.$%
\medskip

\item[\textbf{(e)}] We say that $X\subseteq M$ \textit{is coded }(\textit{in}
$\mathcal{M)}$ if there is some $c\in M$ such that $X=\{x\in M:x\in ^{%
\mathcal{M}}c\}.$ $X$\ is \textit{piecewise coded} in $\mathcal{M}$ if for
each $m\in M,$ $\left\{ x\in M:x\in ^{\mathcal{M}}m\right\} \cap X$ is coded$%
.$ \medskip

\item[\textbf{(f)}] Given an $\mathcal{L}_{\mathrm{set}}$-formula\textit{\ }$%
\varphi ,$ and a variable $x$ not occurring in $\varphi $, \textit{the
relativization of} $\varphi $ \textit{to} $x$, denoted $\varphi ^{x}$, is
the $\Delta _{0}$-formula obtained by restricting all the bound variables of 
$\varphi $ to $x$. \medskip

\bigskip \bigskip
\end{enumerate}

\begin{center}
\textbf{3.~TRUTH AND SET THEORY: THE RUDIMENTS }\bigskip
\end{center}

\noindent This section presents the rudiments of Tarskian satisfaction and
truth in the context of set theory. Subsection 3.1 is concerned with
satisfaction and truth over \textit{set structures}, whereas Subsection 3.2
is devoted to truth and satisfaction over the \textit{entire universe of sets%
}.\bigskip

\begin{center}
\textbf{3.1.~Tarskian satisfaction and truth over set structures}\medskip
\end{center}

\noindent The notions of truth and satisfaction are often used
interchangeably in mathematical logic, basically because in most contexts --
including the one in this section -- the two notions are interdefinable. Let
us review the textbook definition of these notions.\footnote{%
The foundational role of Tarski's definition of truth in a structure, and
its philosophical ramifications have been extensively explored from various
perspectives; my own favorite accounts given by logicians are those of
Feferman \cite{Feferman-tarskian-truth} and Hodges \cite{Hodges-on-truth}.
In the latter article, Hodges writes:
\par
\noindent \textquotedblleft I believe that the first time Tarski explicitly
presented his mathematical definition of truth in a structure was his joint
paper \cite{Tarski+Vaught} with Robert Vaught. This seems remarkably late.
Putting Tarski's \textit{Concept of truth} paper side by side with
mathematical work of the time, both Tarski's and other people's, I think
there is no doubt that Tarski had in his hand all the ingredients for the
definition of truth in a structure by 1931, twenty-six years before he
published it. [...]\ I believe there were some genuine difficulties, not all
of them completely resolved today, and they fully justify Tarski's
caution.\textquotedblright} We have two reasons for doing so; the first is
to establish notation; the other is to revisit the important fact that the
definition of satisfaction is based on \textit{recursive} clauses, which
require an appropriate set-theoretic framework so as to yield a non-circular
definition\textit{. }\medskip

\noindent Suppose $\mathcal{M}=(M,\cdot \cdot \cdot )$ is an $\mathcal{L}$%
-structure (i.e., $\mathcal{L}$ is the language/signature of $\mathcal{M}$); 
$\varphi $ is an $\mathcal{L}$-formula of first order logic, and $\alpha $
is an $\mathcal{M}$-assignment for $\varphi $, i.e.,\medskip

\begin{center}
$\alpha :\mathrm{FV}(\varphi )\rightarrow M$, \medskip
\end{center}

\noindent where $\mathrm{FV}(\varphi )$ is the set of free variables of $%
\varphi $. For such a triple $(\mathcal{M},\varphi ,\alpha )$, the \textit{%
ternary} relation $\mathcal{M}\models \varphi \lbrack {\alpha }]$ is defined
by recursion on the complexity of $\varphi $ via the following clauses. In
what follows we use $v$ and its indexed variants to range over the variables
of first order logic, and in the interest of succinctness, we assume that
first order logic is based on the logical constants $\{\lnot ,\vee ,\exists
\}.$ \medskip

\begin{enumerate}
\item[$(1)$] $\mathcal{M}\models R(v_{0},\cdot \cdot \cdot ,v_{k-1})[\alpha ]
$ \ iff $\ R^{\mathcal{M}}(\alpha (v_{0}),\cdot \cdot \cdot ,\alpha
(v_{k-1}))$; and more generally: 
\begin{equation*}
\mathcal{M}\models R(t_{0},\cdot \cdot \cdot ,t_{k-1})[\alpha ]\ \ \mathrm{%
iff}\ \ R^{\mathcal{M}}(\hat{\alpha}(t_{0}),\cdot \cdot \cdot ,\hat{\alpha}%
(t_{k-1})).
\end{equation*}%
Here $R$ is a $k$-ary relation symbol in $\mathcal{L}$, each $t_{i}$ is an $%
\mathcal{L}$-term, and $\hat{\alpha}$ is the natural extension of $\alpha $
to $\mathcal{L}$-terms $t$ such that the free variables of $t$ are in the
domain of $\alpha $.\medskip 

\item[$(2)$] $\mathcal{M}\models \lnot \varphi \lbrack \alpha ]$ \ iff ~$%
\mathcal{M}\nvDash \varphi \lbrack \alpha ]$.\medskip

\item[$(3)$] $\mathcal{M}\models (\varphi _{1}\vee \varphi _{2})[\alpha ]$ \
iff $\mathcal{M}\models \varphi _{1}[\alpha _{1}]$ or $\mathcal{M}\models
\varphi _{2}[\alpha _{2}]$, where $\alpha _{i}=\alpha \upharpoonright 
\mathrm{FV}(\varphi _{i})$.\medskip

\item[$(4)$] $\mathcal{M}\models \exists v\,\varphi \lbrack \alpha ]$ \ iff
there is some $m\in M$ such that $\mathcal{M}\models \varphi \lbrack \alpha
_{m}^{v}]$. Here $\alpha _{m}^{v}$ is the modification of $\alpha $ that
sends $v$ to $m$ and is otherwise the same as $\alpha .$\medskip
\end{enumerate}

\noindent In this context, given $(\mathcal{M},\varphi ,\alpha )$, where $%
\mathcal{M}$ is a \textit{set} (as opposed to a proper class), it is routine
-- and admittedly tedious -- to write down a formula $\mathrm{sat}(x,y,z)$
in the language of set theory such that $\mathsf{ZF}$ proves that $\mathrm{%
sat}(\mathcal{M},\varphi ,{\alpha })$ satisfies the above recursive clauses
(1) through (4) when $\mathcal{M}\models \varphi \lbrack {\alpha }]$ is
replaced with $\mathrm{sat}(\mathcal{M},\varphi ,{\alpha )}$. The
construction also makes it clear that, provably in $\mathsf{ZF}$, we have:

\begin{itemize}
\item $\mathrm{sat}(\mathcal{M},\varphi ,{\alpha })$ is equivalent to both a 
$\Sigma _{1}$-formula as well as a $\Pi _{1}$-formula of $\mathcal{L}_{%
\mathrm{set}}$ with parameter $\mathcal{M}$. Moreover, provably in $\mathsf{%
ZF}$, $\left\{ \left( \varphi ,\alpha \right) :\mathrm{Sat}\left( \mathcal{M}%
,\varphi ,{\alpha }\right) \right\} $ forms a set, which we will denote by $%
\mathrm{Sat}(\mathcal{M}{)}$. We will refer to $\mathrm{Sat}(\mathcal{M}{)}$
as \textit{the Tarskian satisfaction predicate} for $\mathcal{M}$.
\end{itemize}

\noindent By recasting the above story in model-theoretical terms, we arrive
at the following theorem.\medskip

\noindent \textbf{Theorem.~3.1.1.}~(Tarski's Definability/Codability of
Truth) \textit{Suppose} $\mathcal{N}$ \textit{is a model of} \textit{a
sufficiently strong} \textit{fragment of }$\mathsf{ZF}$ (\textit{see Remark
3.1.2}). \textit{If} $\mathcal{M}$ \textit{is a structure coded as an
element of }$\mathcal{N}$, \textit{then there is a unique }$s\in N$ \textit{%
such that }$\mathcal{N}\models \left[ s=\mathrm{Sat}(\mathcal{M}{)}\right] ,$
\textit{i.e.,} 
\begin{equation*}
\mathcal{N}\models \left[ s=\left\{ \left( \varphi ,\alpha \right) :\mathrm{%
Sat}(\mathcal{M},\varphi ,{\alpha )}\right\} \right] .
\end{equation*}%
\textit{Moreover, within }$\mathcal{N}$, $s$ \textit{is} $\Delta _{1}$ 
\textit{in the parameter} $\mathcal{M}$.\footnote{%
When $\mathcal{M}$ is the standard model of arithmetic, $s$ is $\Delta
_{1}^{1}$. More generally, within $\mathsf{Z}_{2}$\textit{\ }(second order
arithmetic, also known as first order analysis), the satisfaction predicate
for arithmetical formulae is $\Delta _{1}^{1}$; and within $\mathsf{KM}$%
\textit{\ }(Kelley-Morse theory of classes), the satisfaction predicate for
set-theoretical formulae is $\Delta _{1}^{1}.$} \medskip

\noindent \textbf{Remark.~3.1.2.}~There are two canonical fragments of $%
\mathsf{ZF}$ that are `sufficiently strong' for the purposes of Theorem
3.1.1, namely:\medskip

\begin{enumerate}
\item[\textbf{(a)}] $\mathsf{KP}$\ (see Definition 2.1(j)), as shown by
Friedman, Li, and Wong in \cite[Lemma 4.1]{Friedman at al}.\medskip

\item[\textbf{(b)}] The fragment $\mathsf{M}_{0}+\mathsf{Infinity}$ of $%
\mathsf{Z}$ (see Definition 2.1(k)), as shown by Mathias \cite[Proposition
3.10]{Mathias-MacLane}.\medskip
\end{enumerate}

\noindent Note, however, that much weaker systems suffice if one only wishes
to have a set theory within which the Tarskian satisfaction relation of
every internal set structure is \textit{definable}, as opposed to \textit{%
coded as a set}. As shown by Mathias \cite[Proposition 10.37]%
{Mathias-fixingDevlin} the system $\mathsf{DS}$ (for `Devlin Strengthened'),
and even the weaker system $\mathsf{MW}$ (for `Middle Way') are capable of
defining the Tarskian satisfaction predicate for set structures in a $\Delta
_{1}$-manner. \medskip 

\noindent \textbf{Remark. 3.1.3.}~In model theory one often uses the notion
of the \emph{theory of an }$\mathcal{L}$\emph{-structure $\mathcal{M}$},
denoted $\mathrm{Th}(\mathcal{M}).$ Officially speaking, $\mathrm{Th}(%
\mathcal{M})$ is defined as the set of $\mathcal{L}$-sentences $\sigma $
such that $(\sigma ,\varnothing )\in \mathrm{Sat}(\mathcal{M}{)}$; here we
are using the common practice of using the term `sentence'\ to mean `formula
with no free variables'. Model theorists also commonly use the notion of 
\emph{the elementary diagram of a model $\mathcal{M}$}, here denoted $%
\mathrm{ED}(\mathcal{M})$.\footnote{%
In Chang and Keisler's text \cite{Chang-Keisler}, the elementary diagram of $%
\mathcal{M}$ is denoted $\mathrm{Th}(\mathcal{M},m)_{m\in M}.$ Hodges' text 
\cite{Hodges-text} \ uses the notation $\mathrm{Th}(M_{M})$ for the
elementary diagram of a model $M$.}

\begin{itemize}
\item Officially speaking, $\mathrm{ED}(\mathcal{M})$ is defined as the set
of sentences $\varphi (\dot{m}_{0},\cdot \cdot \cdot ,\dot{m}_{k-1})$ in the
language $\mathcal{L}_{M},$ obtained by enriching $\mathcal{L}$ with
constant symbols $\dot{m}$ for each $m\in M$, such that $\left( \varphi
(x_{0},...,x_{k-1}),\alpha \right) \in \mathrm{sat}_{\mathcal{M}}$, where $%
\alpha (x_{i})=m_{i}$ for $i<k.$
\end{itemize}

\noindent However, one can readily turn the tables around and directly
define $\mathrm{ED}(\mathcal{M})$ without a detour through $\mathrm{Sat}(%
\mathcal{M}{)}$, e.g., as in Shoenfield's classic textbook \cite[Section 2.5]%
{Shoenfield-textbook}. More specifically, the relation $\mathcal{M}\models
\sigma $ can be alternatively defined by recursion on the complexity of $%
\mathcal{L}_{M}$-\textit{sentences} $\sigma $ using the following clauses.
\medskip

\begin{enumerate}
\item[$(i)$] $\mathcal{M}\models R(\dot{m}_{0},\cdot \cdot \cdot ,\dot{m}%
_{k-1})$ iff $R^{\mathcal{M}}(m_{0},\cdot \cdot \cdot ,m_{k-1});$ and more
generally: 
\begin{equation*}
\mathcal{M}\models R(t_{0},\cdot \cdot \cdot ,t_{k-1})\ \ \mathrm{iff}\ R^{%
\mathcal{M}}(t_{0}^{\mathcal{M}},\cdot \cdot \cdot ,t_{k-1}^{\mathcal{M}}).
\end{equation*}%
Here $R$ is a $k$-ary relation symbol in $\mathcal{L}$, and each $t_{i}$ is
a closed $\mathcal{L}_{M}$-term.\medskip

\item[$(ii)$] $\mathcal{M}\models \lnot \sigma $ iff $\mathcal{M}\nvDash
\sigma $.\medskip

\item[$(iii)$] $\mathcal{M}\models \sigma _{1}\vee \sigma _{2}$ iff $%
\mathcal{M}\models \sigma _{1}$ or $\mathcal{M}\models \sigma _{2}$.\medskip

\item[$(iv)$] $\mathcal{M}\models \exists v\,\varphi (v)$ iff there is some $%
m\in M$ such that $\mathcal{M}\models \varphi (\dot{m}/v)$. \bigskip
\end{enumerate}

\begin{center}
\textbf{3.2.~Truth and satisfaction over the universe of sets}\medskip
\end{center}

In the previous subsection we reviewed the basics of satisfaction and truth
applied to `small' structures within set theory, i.e., \textit{set-structures%
} (as opposed to proper class structures). In this subsection we examine the
notions of \textit{satisfaction classes} and \textit{truth classes}, which
respectively capture the notions of satisfaction and truth over the \textit{%
entire universe of sets}. As we shall see in Proposition 3.2.6, a truth
class is essentially an \textit{extensional} satisfaction class. The theory
of sets required for this purpose can be quite modest, for definiteness we
have chosen $\mathsf{KP}$. The weakest set theory for the purposes of
Definition 3.2.1 is a \textit{sequential} theory, the canonical example of
which is $\mathsf{AS}$ (Adjunctive Set Theory), but the verification of
sequentiality of $\mathsf{AS}$ is quite laborious; see, e.g., Visser's [Vi].
To get an idea of the work involved in the bootstrapping necessary to
accommodate a full truth predicate over a sequential theory, see Fangjing
Xiong's thesis \cite{Xiong}, in which the sequential theory $\mathsf{PA}^{-}$
serves as the base theory. \medskip 

\noindent \textbf{3.2.1.~Definition.}~We will use the following
abbreviations relating to the set-theoretic coding of syntax; note that all
the formulae in the list below are $\mathcal{L}_{\mathrm{set}}$%
-formulae.\medskip

\begin{enumerate}
\item[\textbf{(a)}] $\mathrm{Form}(x)$ expresses \textquotedblleft $x$ is an 
$\mathcal{L}_{\mathrm{set}}$-formula\textquotedblright , and $\mathrm{Form}%
_{k}(x)$ is the conjunction of $\mathrm{Form}(x)$ and \textquotedblleft $x$
has $k$ free variables\textquotedblright .\footnote{%
We assume that $\mathsf{KP}\vdash \forall x\left( \mathrm{Form}%
(x)\rightarrow x\in \mathrm{V}_{\omega }\right) .$ Note that there is a
definable bijection in $\mathsf{KP}$\ between $\omega $ and $\mathrm{V}%
_{\omega },$ and coding on the hereditarily finite sets is much easier than
coding on $\omega .$}\medskip

\item[\textbf{(b)}] $\mathrm{Sent}(x)$ is the conjunction of $\mathrm{Form}%
(x)$ and \textquotedblleft $x$ has no free variables".\medskip

\item[\textbf{(c)}] $\mathrm{Var}(x)$ expresses \textquotedblleft $x$ is a
variable\textquotedblright .\medskip

\item[\textbf{(d)}] $\mathrm{Asn}(\alpha )$ expresses \textquotedblleft $%
\alpha $ is an assignment\textquotedblright , where an assignment here
simply refers to a function whose domain consists of a (finite) set of
variables.\medskip 

\item[\textbf{(e)}] $y\in \mathrm{FV}(x)$ is the conjunction of $\mathrm{Form%
}(x)$ and \textquotedblleft $y$ is a free variable of $x$\textquotedblright
.\medskip

\item[\textbf{(f)}] $y\in \mathrm{Dom}(\alpha )$ expresses \textquotedblleft
the domain of $\alpha $ includes $y$\textquotedblright .\medskip

\item[\textbf{(g)}] $\mathrm{Asn}(\alpha ,x)$ expresses \textquotedblleft $%
\alpha $ is an assignment for $x$\textquotedblright , i.e. it is the
conjunction of $\mathrm{Form}(x)$, $\mathrm{Asn}(\alpha )$, and
\textquotedblleft the domain of $\alpha $ is $\mathrm{FV}(x)$%
\textquotedblright .\medskip

\item[\textbf{(h)}] For assignments $\alpha $ and $\beta $, $\beta \supseteq
\alpha $ expresses \textquotedblleft the domain of $\beta $ extends the
domain of $\alpha $ and $\alpha (v)=\beta (v)$ for all $v\in \mathrm{Dom}%
(\alpha )$\textquotedblright .\medskip

\item[\textbf{(i)}] $x\vartriangleleft y$ expresses \textquotedblleft $x$ is
an immediate subformula of $y$\textquotedblright , i.e., $x\vartriangleleft
y $ abbreviates the conjunction of $y\in \mathrm{Form}$ and the following
disjunction:\medskip

$\left( y=\lnot x\right) \vee \exists z\left( \left( y=x\vee z\right) \vee
\left( y=z\vee x\right) \right) \vee \exists v\in \mathrm{Var}\mathsf{\ }%
\left( y=\exists v\ x\right) .$\medskip

\item[\textbf{(j)}] In a context where $\mathrm{Form}(x)$ holds, we write $%
(x=\lnot y)$ instead of the more formal 
\begin{equation*}
\exists y(\mathrm{Form}(y)\wedge x=\dot{\lnot}y),
\end{equation*}%
where $\dot{\lnot}$ is the definable function whose output is the code of $%
\lnot \varphi $, when given the code of $\varphi $ as input$.$ We follow a
similar practice for the expressions $\left( x=y_{1}\vee y_{2}\right) $ and $%
(x=\exists v\,y).$\medskip
\end{enumerate}

\noindent \textbf{3.2.2.~Definition.}~The theory $\mathsf{CS}^{-}\left( 
\mathsf{F}\right) $ defined below is formulated in an \textit{expansion} of $%
\mathcal{L}_{\mathrm{set}}$ by adding a fresh \textit{binary} predicate $%
\mathsf{S}(x,y)$ (denoting satisfaction) and a fresh unary predicate $%
\mathsf{F}$ (denoting a specified collection of formulae). The binary/unary
distinction is of course not an essential one since $\mathsf{KP}$ has access
to a definable pairing function. However, the binary/unary distinction 
\textit{at the conceptual level} marks the key difference between the
concepts of satisfaction and truth.\medskip

\begin{enumerate}
\item[\textbf{(a)}] $\mathsf{CS}^{-}\left( \mathsf{F}\right) $ is the
conjunction of the universal generalizations of the formulae $(1)$ through $%
(5)$ listed below. In what follows $v$ and $w$ range over variables, while $%
\alpha $ and $\beta $ range over assignments. It is helpful to bear in mind
that the axioms of $\mathsf{CS}^{-}(\mathsf{F})$ collectively express:
\textquotedblleft $\mathsf{F}$ is a subset of $\mathcal{L}_{\mathrm{set}}$%
-formulae that is closed under immediate subformulae; each member of $%
\mathsf{S}$ is an ordered pair of the form $(x,\alpha )$, where $x$ is in $%
\mathsf{F}$ and $\alpha $ is an assignment for $x;$ and $\mathsf{S}$
satisfies Tarski's compositional clauses for a satisfaction
predicate\textquotedblright .\medskip

\begin{enumerate}
\item[(1)] $\left[ \mathsf{F}(x)\rightarrow \mathrm{Form}(x)\right] \wedge %
\left[ y\vartriangleleft x\wedge \mathsf{F}(x)\rightarrow \mathsf{F}(y)%
\right] \wedge \left[ \mathsf{S}(x,\alpha )\rightarrow \left( \mathrm{Form}%
(x)\wedge \alpha \in \mathrm{Asn}(x)\right) \right] .$\medskip

\item[(2)] $\left( 
\begin{array}{c}
\left( \mathsf{S}\left( \left( v=w\right) ,\alpha \right) \leftrightarrow %
\left[ \mathrm{Dom}(\alpha )=\{v,w\}\wedge (\alpha (v)=\alpha (w)\right]
\right) \wedge \\ 
\left( \mathsf{S}\left( \left( v\in w\right) ,\alpha \right) \leftrightarrow %
\left[ \mathrm{Dom}(\alpha )=\{v,w\}\wedge \alpha (v)\in \alpha (w)\right]
\right)%
\end{array}%
\right) $.\medskip

\item[(3)] $\left[ \mathsf{F}(x)\wedge (x=\lnot y)\wedge \alpha \in \mathrm{%
Asn}(x)\right] \rightarrow \left[ \mathsf{S}(x,\alpha )\leftrightarrow \lnot 
\mathsf{S}(y,\alpha )\right] .$\medskip

\item[(4)] $\left[ \mathsf{F}(x)\wedge \left( x=y_{1}\vee y_{2}\right)
\wedge \alpha \in \mathrm{Asn}(x)\right] \rightarrow \left[ \mathsf{S}%
(x,\alpha )\leftrightarrow \left( \mathsf{S}\left( y_{1},\alpha
\upharpoonright \mathrm{FV}(y_{1})\right) \vee \mathsf{S}\left( y_{2},\alpha
\upharpoonright \mathrm{FV}(y_{2})\right) \right) \right] .$\medskip

\item[(5)] $\left[ \mathsf{F}(x)\wedge (x=\exists v\,y)\wedge \alpha \in 
\mathrm{Asn}(x))\right] \rightarrow \left[ \mathsf{S}(x,\alpha
)\leftrightarrow \exists \beta \supseteq \alpha \ \mathsf{S}(y,\beta )\right]
$\medskip
\end{enumerate}

\item[\textbf{(b)}] $\mathsf{CS}^{-}$ is the theory whose axioms are
obtained by substituting the predicate $\mathsf{F}(x)$ by the $\mathcal{L}_{%
\mathrm{set}}$-formula $\mathrm{Form}(x)$ in the axioms of $\mathsf{CS}^{-}(%
\mathsf{F})$. \textit{Thus the axioms in} $\mathsf{CS}^{-}$ \textit{are
formulated in the language obtained by adding} $\mathsf{S}$ \textit{to} $%
\mathcal{L}_{\mathrm{set}}$ (\textit{with no mention of} $\mathsf{F}$).$%
\medskip $

\item[\textbf{(c)}] Given any base theory $\mathsf{B}\supseteq \mathsf{KP},$
we write $\mathsf{CS}^{-}[\mathsf{B}]$ as a shorthand for $\mathsf{CS}%
^{-}\cup \mathsf{B}.$\medskip
\end{enumerate}

\noindent \textbf{3.2.3.~Definition}.~Suppose $\mathcal{M}\models \mathsf{KP}
$, and let $F\subseteq \mathrm{Form}^{\mathcal{M}}=\left\{ m\in M:\mathcal{M}%
\models \mathrm{Form}(m)\right\} $. $\medskip $

\begin{enumerate}
\item[\textbf{(a)}] A subset $S$ of $M^{2}$ is said to be an $F$-\textit{%
satisfaction class on }$\mathcal{M}$\textit{\ }if $(\mathcal{M},F,S)\models 
\mathsf{CS}^{-}(\mathsf{F})$, here the interpretation of $\mathsf{F}$ is $F$
and the interpretation of $\mathsf{S}$ is $S.$ $S$ is a \textit{satisfaction
class} on $\mathcal{M}$ if $S$ is an $F$-satisfaction class for some $F$.$%
\medskip $

\item[\textbf{(b)}] An $F$-satisfaction class $S$ is \textit{extensional }if
for all $\varphi _{0}$ and $\varphi _{1}$ in $F$, and for any $\varphi _{0}$%
-assignment $\alpha _{0}$ and any $\varphi _{1}$-assignment $\alpha _{1}$,
we have: 
\begin{equation*}
(\mathcal{M},F,S)\models \left[ \left[ (\varphi _{0},\alpha _{0})\thicksim
(\varphi _{1},\alpha _{1})\right] \longrightarrow \left[ \mathsf{S}(\varphi
_{0},\alpha _{0})\leftrightarrow \mathsf{S}(\varphi _{1},\alpha _{1})\right] %
\right] ,
\end{equation*}%
where $(\varphi _{0},\alpha _{0})\thicksim (\varphi _{1},\alpha _{1})$ means
that $\varphi _{0}$ and $\varphi _{1}$ are the same except for their free
variables, and for all variables $x$ and $y$, if $x$ occurs freely in the
same position in $\varphi _{0}$ as $y$ does in $\varphi _{1}$, then $\alpha
_{0}(x)=\alpha _{1}(y).$ $\medskip $

\item[\textbf{(c)}] A subset $S$ of $M$ is said to be a \textit{full} 
\textit{satisfaction class on }$\mathcal{M}$\textit{\ }if $(\mathcal{M}%
,S)\models \mathsf{CS}^{-}$ for $F=\mathrm{Form}^{\mathcal{M}}$.\medskip 
\end{enumerate}

\noindent \textbf{3.2.4.~Definition.}~The theory $\mathsf{CT}^{-}\left( 
\mathsf{F}\right) $ defined below is formulated in the language obtained by
augmenting $\mathcal{L}_{\mathrm{set}}$ with a fresh \textit{unary}
predicate $\mathsf{T}(x)$ (denoting truth) and a fresh unary predicate $%
\mathsf{F}$ (denoting a specified collection of formulae). \medskip

\begin{enumerate}
\item[\textbf{(a)}] Reasoning within the theory $\mathsf{KP}$ (and not
within the metatheory) we fix a function $a\mapsto \dot{a}$, that designates
constant symbols $\dot{a}$ for each object $a$ in the universe of sets
(e.g., $\dot{a}$ is defined as the ordered pair $\left\langle
a,3\right\rangle $ in Devlin's monograph \cite{Devlin-book-on-L} ). \medskip

\item[\textbf{(b)}] $\mathrm{Sent}_{\mathsf{F}}^{+}(x)$ expresses
\textquotedblleft $x$ is an $\mathcal{L}_{\mathrm{set}}$-sentence obtained
by substituting constants from $\left\{ \dot{a}:a\in \mathrm{V}\right\} $
for each free variable of some formula in $\mathsf{F}$\textquotedblright .
We write $\mathrm{Sent}^{+}(x)$ if $\mathsf{F}=\mathrm{Form}$.\medskip

\item[\textbf{(c)}] $y\vartriangleleft x$ expresses \textquotedblleft $y$ is
an immediate subformula of $x$\textquotedblright\ as in Definition
3.2.1(i).\medskip 

\item[\textbf{(d)}] ${\mathnormal{\mathsf{F}}}_{\leq 1}(\varphi (v))$
expresses \textquotedblleft $\mathsf{F}(\varphi )$ and $\varphi $ has at
most one free variable $v$\textquotedblright .\medskip

\item[\textbf{(e)}] $\varphi \lbrack \dot{x}/v]$ is (the code of) the
formula obtained by substituting all occurrences of the variable $v$ in $%
\varphi $ with the constant symbol $\dot{x}$ representing $x.$\medskip

\item[\textbf{(f)}] $\mathsf{CT}^{-}(\mathsf{F})$ satisfies the universal
generalizations of the conjunction of $(1)$ through $(5)$ below, where $%
\varphi (v)$ ranges over (codes of) $\mathcal{L}_{\mathrm{set}}$%
-formulae:\medskip

\begin{itemize}
\item[$(1)$] $\left[ \mathsf{T}(\varphi )\rightarrow \mathrm{Sent}%
^{+}(\varphi )\right] \wedge \left[ \psi \vartriangleleft \varphi \wedge 
\mathsf{F}(\varphi )\rightarrow \mathsf{F}(\psi )\right] .$\medskip

\item[$(2)$] $\left[ \left( \mathsf{T}\left( \dot{x}=\dot{y}\right)
\leftrightarrow x=y\right) \wedge \left( \mathsf{T}\left( \dot{x}\in \dot{y}%
\right) \leftrightarrow x\in y\right) \right] .$\medskip

\item[$(3)$] $\left[ \mathrm{Sent}_{\mathsf{F}}^{+}(\varphi )\wedge \mathrm{%
Sent}_{\mathsf{F}}^{+}(\psi )\right] \rightarrow \left[ \left( \varphi
=\lnot \psi \right) \rightarrow \left[ \mathsf{T}(\varphi )\leftrightarrow
\lnot \mathsf{T}(\psi )\right] \right] \mathsf{.}$\medskip 

\item[$(4)$] $\left[ \mathrm{Sent}_{\mathsf{F}}^{+}(\varphi )\wedge \mathrm{%
Sent}_{\mathsf{F}}^{+}(\psi _{1})\wedge \mathrm{Sent}_{\mathsf{F}}^{+}(\psi
_{2})\right] \rightarrow \left[ \left( \varphi =\psi _{1}\vee \psi
_{2}\right) \rightarrow \left( \mathsf{T}(\varphi )\leftrightarrow \left( 
\mathsf{T}(\psi _{1})\vee \left( \mathsf{T}(\psi _{2})\right) \right)
\right) \right] \mathsf{.}$\medskip

\item[$(5)$] $\left[ \mathrm{Sent}_{\mathsf{F}}^{+}(\varphi )\wedge \mathsf{F%
}_{\leq 1}(\psi (v))\right] \rightarrow \left[ \left( \varphi =\exists
v\,\psi (v)\right) \rightarrow \left[ \mathsf{T}(\varphi )\leftrightarrow
\exists x\,\mathsf{T}(\mathsf{\psi (}\dot{x}/v\mathsf{))}\right] \right] .$%
\medskip 
\end{itemize}

\item[\textbf{(g)}] $\mathsf{CT}^{-}$ is the theory whose axioms are
obtained by substituting the predicate $\mathsf{F}(x)$ by the $\mathcal{L}_{%
\mathrm{set}}$-formula $\mathrm{Form}(x)$ in the axioms of $\mathsf{CT}^{-}(%
\mathsf{F})$. \textit{Thus the axioms of} $\mathsf{CT}^{-}$ \textit{are
formulated in the language} $\mathcal{L}_{\mathrm{set}}(\mathsf{T})$ \textit{%
obtained by adding} $\mathsf{T}$ \textit{to} $\mathcal{L}_{\mathrm{set}}$ (%
\textit{with no mention of} $\mathsf{F}$). \medskip

\item[\textbf{(h)}] $\mathsf{CT}^{-}[\mathsf{B}]$ is a shorthand for $%
\mathsf{CT}^{-}+\mathsf{B}$, where $\mathsf{CT}^{-}$ is as in Definition
3.2.4 and $\mathsf{B}$ is an $\mathcal{L}_{\mathrm{set}}$-theory (referred
to as a \textit{base theory}) extending $\mathsf{KP}$.\medskip 
\end{enumerate}

\noindent \textbf{3.2.5.~Definition}.~Let $\mathcal{M}\models \mathsf{KP}$,
and suppose $F\subseteq \mathrm{Form}^{\mathcal{M}}$, and $F$ is closed
under direct subformulae of $\mathcal{M}$. Recall that $\left( \mathrm{Sent}%
_{\mathsf{F}}^{+}\right) ^{(\mathcal{M},F)}$ consists of $x\in M$ such that $%
(\mathcal{M},F)$ satisfies \textquotedblleft $x$ is an $\mathcal{L}_{\mathrm{%
set}}$-sentence obtained by substituting constants from $\left\{ \dot{m}%
:m\in \mathrm{V}\right\} $ for the free variables of a formula in $\mathsf{F}
$\textquotedblright .$\medskip $

\begin{enumerate}
\item[\textbf{(a)}] A subset $T$ of $M$ is an $F$-\textit{truth class }on%
\textit{\ }$\mathcal{M}$\textit{\ }if $(\mathcal{M},F,T)\models \mathsf{CT}%
^{-}(\mathsf{F})$, here the interpretation of $\mathsf{F}$ is $F$ and the
interpretation of $\mathsf{T}$ is $T.$ $T$ is a \textit{truth class} on $%
\mathcal{M}$ if $T$ is an $F$-truth class for some $F$.\medskip

\item[\textbf{(b)}] For $k\in \mathbb{\omega }^{\mathcal{M}}$, $T$ is a $%
\Sigma _{k}$-\textit{truth class }on\textit{\ }$\mathcal{M}$ if $T$ is an $F$%
-truth class on $\mathcal{M}$, where $F$ is the collection of all $\Sigma
_{k}$-formulae in the sense of $\mathcal{M}.$ The notion of $\mathrm{Depth}%
_{k}$-truth class is defined similarly, where $\mathrm{Depth}_{k}$ is the
collection of formulae whose logical depth is at most $k$ (see Definition
7.2).\medskip\ 

\item[\textbf{(c)}] $T$ is a \textit{full} \textit{truth class }on\textit{\ }%
$\mathcal{M}$\textit{\ }if $(\mathcal{M},T)\models \mathsf{CT}^{-}$;
equivalently: if $(\mathcal{M},F,T)\models \mathsf{CT}^{-}(\mathsf{F})$ for $%
F=\mathrm{Form}^{\mathcal{M}}$.\medskip
\end{enumerate}

The following proposition codifies the inter-definability of truth classes
and extensional satisfaction classes. The relationship between extensional
satisfaction classes and truth classes (in the context of arithmetic) was
first made explicit in \cite{Ali+albert-short}; for further elaborations see 
\cite{Cieslinski-book} and \cite{Wcislo-Collection}. Let $\mathrm{dot}(x)$
be the $\mathsf{KP}$-definable function\textit{\ }$m\mapsto \ \dot{m}$ where 
$\dot{m}$ is the constant associated with $m\in \mathrm{V}$, and $\varphi
\ast \alpha $ be the sentence $\varphi (\mathrm{dot}\circ \alpha )$, i.e.,
the sentence obtained by replacing each occurrence of a free variable $x$ of 
$\varphi $ with the constant symbol $\dot{m}$, where $\alpha (x)=m.$\medskip

\noindent \textbf{3.2.6.~Proposition.}~\textit{Suppose }$\mathcal{M}\models 
\mathsf{KP}$, $T$\textit{\ is an }$F$\textit{-truth class on }$\mathcal{M}$%
\textit{, and }$S$\textit{\ is an extensional }$F$\textit{-satisfaction
class on }$\mathcal{M}.$\textit{\medskip }

\begin{enumerate}
\item[$(a)$] $\mathcal{S}(T)$\textit{\ is an extensional }$F$\textit{%
-satisfaction class on }$\mathcal{M}$\textit{, where }$\mathcal{S}(T)$ 
\textit{is defined as the collection of ordered pairs }$(\varphi ,\alpha )$ 
\textit{such that} $\varphi \ast \alpha \in T.$\medskip

\item[$(b)$] $\mathcal{T}(S)$\textit{\ is an }$F$\textit{-truth class on }$%
\mathcal{M}$\textit{, where }$\mathcal{T}(S)$\textit{\ is defined as the
collection of }$\varphi $ \textit{such that} $\left( \varphi ,\varnothing
\right) \in S$ (\textit{where} $\varnothing $ \textit{is the empty assignment%
}).\medskip

\item[$(c)$] $\mathcal{S}(\mathcal{T}(S))=S$, \textit{and} $\mathcal{T}(%
\mathcal{S(}T))=T.$\textit{\medskip }
\end{enumerate}

\noindent \textbf{3.2.7}.~\textbf{Theorem}~(Model-theoretic formulation of
Tarski's Undefinability of Satisfaction). \textit{Suppose} $\mathcal{M}$ 
\textit{is an }$\mathcal{L}$\textit{-structure. Fix some }$m\in M$, \textit{%
and let} $c_{m}$ \textit{be a constant added to} $\mathcal{L}$ \textit{for
denoting} $m$. \textit{Finally, let} 
\begin{equation*}
\varphi (x)\mapsto \#(\varphi (x))\in M
\end{equation*}%
\textit{be a mapping that assigns an element of} $M$ \textit{to each unary} $%
\mathcal{L}(c_{m})$-\textit{formula}. \textit{There is no binary }$\mathcal{L%
}(c_{m})$-\textit{formula} $S(x,y)$ \textit{such that for all unary }$%
\mathcal{L}(c_{m})$-\textit{formulae} $\varphi (x)$, \textit{we have}:

\begin{equation*}
\left( \mathcal{M},m\right) \models \forall x\,\left[ S(\#(\varphi
),x)\leftrightarrow \varphi (x)\right] .
\end{equation*}

\noindent \textbf{Proof\footnote{%
This proof is reminiscent of Russell's Paradox (1901), and of the proof of
Cantor's theorem (1891) on nonexistence of a surjection of a set $X$ onto $%
\mathcal{P}(X).$ I learned about it from Kossak's \cite[Theorem 2.3]%
{Kossak-TUT}, where it is attributed to Schmerl. A similar proof can be
found in Kripke's lecture notes \cite[page 66]{Kripke-book}. The usual proof
of undefinability of truth is based on the \textit{parametric form} of the
fixed point theorem (see, e.g., \cite[Ch.III, Theorem 2.1]{Hajek and Pudlak}%
), which is provable in theories that interpret Robinson's \textsf{Q. }Also
note that in certain models $\mathcal{M}$ of set theory, $\mathrm{Th}(%
\mathcal{M}$) is $\mathcal{M}$-definable (using a parameter), such models
include recursively saturated ones, and models of the form $\left( \mathrm{%
V\!}_{\alpha },\in \right) $ where $\alpha >\omega .$}.}~Suppose not, and
let $S(x,y)$ be such a formula, and let $R(x):=\lnot S(x,x)$, and let $r\in
M $ such that $r:=\#\left( R(x)\right) .$ Then we have:\textit{\medskip }

\noindent (1) $\ \ \left( \mathcal{M},m\right) \models \forall x\,\left[
S\left( \#(R),x\right) \leftrightarrow R(x)\right] .$\textit{\medskip }

\noindent By (1) and the definition of $R$ we obtain the following
contradiction:\textit{\medskip }

\noindent (2) $\ \ \left( \mathcal{M},m\right) \models S(r,r)\leftrightarrow
R(r)\leftrightarrow \lnot S(r,r).$

\hfill $\square $\textit{\medskip }

\noindent \textbf{3.2.8}.~\textbf{Corollary.}~\textit{If} $S$ \textit{is an }%
$F$\textit{-satisfaction class on a model }$\mathcal{M}$ \textit{of} $%
\mathsf{KP}$ \textit{that} $F$ \textit{includes all standard} $\mathcal{M}$-%
\textit{formulae}, \textit{then} $S$ \textit{is not} $\mathcal{M}$-\textit{%
definable. In particular, no full satisfaction/truth class on }$\mathcal{M}$%
\textit{\ is }$\mathcal{M}$\textit{-definable.}\medskip

\noindent It is a well-known result of Levy \cite{Levy-book} that if $%
\mathcal{M}\models \mathsf{ZF}\mathrm{,}$ then there is a $\Delta _{0}$%
-satisfaction class for\textit{\ }$\mathcal{M}$\textit{\ }that is\textit{\ }%
definable in $\mathcal{M}$ both by a $\Sigma _{1}$-formula and a $\Pi _{1}$%
-formula (see \cite[p.~186]{Jechbook-2003} for a proof). This makes it clear
that for each $n\geq 1,$ there is a $\Sigma _{n}$-satisfaction class for%
\textit{\ }$\mathcal{M}$\textit{\ }that is\textit{\ }definable in $\mathcal{M%
}$ by a $\Sigma _{n}$-formula. This leads to the theorem below.\medskip

\noindent \textbf{3.2.9}.\textbf{~Theorem }(Levy's Partial Definability of
Truth).~ \textit{For each} $n\in \omega $ \textit{there is an} $\mathcal{L}_{%
\mathrm{set}}$-\textit{formula} $\mathrm{True}_{\Sigma _{n}}$ \textit{such
that for all models} $\mathcal{M}$ \textit{of} \textit{a sufficiently strong 
}(\textit{see Remark 3.2.10 below}) $\mathsf{ZF}$, $\mathrm{True}_{\Sigma
_{n}}^{\mathcal{M}}$ \textit{is a} $\Sigma _{n}$-\textit{truth class for} $%
\mathcal{M}$. \textit{Furthermore, for} $n\geq 1$, $\mathrm{True}_{\Sigma
_{n}}$\textit{\ is }$\Sigma _{n}\mathrm{.}$\medskip

\noindent \textbf{3.2.10}.\textbf{~Remark.~}There are two canonical
fragments of $\mathsf{ZF}$ that are `sufficiently strong' for the purposes
of Theorem 3.2.9. One is $\mathsf{KP}$, and the other one is $\mathsf{M}_{0}+%
\mathsf{Infinity}+\mathsf{TC}$, where $\mathsf{M}_{0}$ is as in part (k) of
Definition 2.1, and $\mathsf{TC}$ is the sentence asserting that every set
has a transitive closure. See, e.g., Definitions 2.9 and 2.10 of McKenzie's 
\cite{Zach collection} for the case of $\mathsf{KP}$; a similar construction
works for $\mathsf{M}_{0}+\mathsf{Infinity+TC}$. What is needed in both
cases is $\mathsf{TC}$, plus the ability of the theory to define the
Tarskian satisfaction predicate for set structures. Also note that, as shown
by Pudl\'{a}k and Visser, all sequential theories support partial
satisfaction predicates for formulae of a given quantifier complexity (for
more detail and references, see \cite[Fact F]{Enayat-Visser-PAMS}). The weak
set theory $\mathsf{AS}$ (Adjunctive Set Theory) is a sequential theory; see
e.g., \cite{Visser-Pairs-sets}.\medskip 

\noindent The following result is fundamental. For an exposition, see \cite[%
Theorem 12.14]{Jechbook-2003}.\footnote{%
Levy \cite{Levy-book} refined Theorem 3.2.11 by showing that for each $n\in
\omega $, there is a formula $\theta _{n}(x)$ such $\mathsf{ZF}$ proves that
the collection of ordinals satisfying $\theta _{n}(x)$ is c.u.b.~in the
class of ordinals, and if $\theta _{n}(\alpha )$ then \textrm{V}$_{\alpha
}\prec _{\Sigma _{n}}\mathrm{V}.$} \medskip

\noindent \textbf{3.2.11.~Montague Reflection Theorem \cite{Montague1961}.}~%
\textit{For each} \textit{formula\ }$\varphi (\vec{x}),$ \textit{there is a
formula} $\theta _{\varphi }(y)$ \textit{such that} $\mathsf{ZF}\vdash 
\mathrm{Ref}_{\varphi ,\theta }$, \textit{where} $\mathrm{Ref}_{\varphi
,\theta }$ \textit{is the sentence that expresses}:\medskip

\begin{center}
\textquotedblleft $\left\{ \theta _{\varphi }(\alpha ):\alpha \in \mathrm{Ord%
}\right\} $ is c.u.b.\footnote{%
Here \textquotedblleft c.u.b.\textquotedblright\ stands for
\textquotedblleft closed and unbounded\textquotedblright . Recall for $%
X\subseteq \mathrm{Ord},$ $X$ is said to be \textit{closed} if for each
limit ordinal $\alpha $, if $X\cap \alpha \in X$, then $\alpha \in X.$}~in $%
\mathrm{Ord}$\textquotedblright\ and $\forall \alpha \left( \theta _{\varphi
}(\alpha )\rightarrow \left( \mathrm{V\!}_{\alpha }\prec _{\varphi }\mathrm{V%
}\right) \right) ,$\medskip
\end{center}

\noindent \textit{where} $\left( \mathrm{V\!}_{\alpha }\prec _{\varphi }%
\mathrm{V}\right) $ \textit{is shorthand for the following} $\mathcal{L}_{%
\mathrm{set}}$-\textit{formula}:\medskip

\begin{center}
$\forall x_{0}\in \mathrm{V\!}_{\alpha }\cdot \cdot \cdot \forall x_{k-1}\in 
\mathrm{V\!}_{\alpha }\,\left[ \varphi \left( x_{0},\cdot \cdot \cdot
,x_{k-1}\right) \longleftrightarrow \overset{(\mathrm{V\!}_{\alpha },\in
)~\models ~\varphi \left( \dot{x}_{0},\cdot \cdot \cdot ,\dot{x}%
_{k-1}\right) }{\overbrace{\varphi \left( \dot{x}_{0},\cdot \cdot \cdot ,%
\dot{x}_{k-1}\right) \in \mathrm{ED}(\mathrm{V\!}_{\alpha },\in )\ }}\right]
.$\footnote{%
As noted by Montague, the proof of Theorem 3.2.11 shows that the \textrm{%
V\negthinspace }$_{\alpha }$-hierarchy can be replaced by any definable,
monotone, and continuous hierarchy \textrm{W\negthinspace }$_{\alpha }$
whose union is $\mathrm{V}$. For example, in the presence of the axiom of
choice, we can let \textrm{W\negthinspace }$_{\alpha }=\mathrm{H}_{\aleph
_{\alpha }}$, where $\mathrm{H}_{\kappa }$ is the collection of sets that
are hereditarily of cardinality less than $\kappa $.}\medskip
\end{center}

\noindent \textbf{3.2.12~Remark.}~In the displayed biconditional in
Reflection Theorem 3.2.11, the right-hand-side can be replaced with $\varphi
^{\mathrm{V\!}_{\alpha }}\left( x_{0},\cdot \cdot \cdot ,x_{k-1}\right) $.
This is because of the fact that$,$ $\mathsf{ZF}\vdash \pi _{\varphi }$,
where: $\mathsf{\ }$%
\begin{equation*}
\pi _{\varphi }:=\forall x_{0}\in m\cdot \cdot \cdot \forall x_{k-1}\in m\, 
\left[ \left( \varphi ^{m}\left( x_{0},\cdot \cdot \cdot ,x_{k-1}\right)
\leftrightarrow \left( \varphi \left( \dot{x}_{0},\cdot \cdot \cdot ,\dot{x}%
_{k-1}\right) \in \mathrm{ED}(m,\in )\right) \right) \right] .
\end{equation*}%
\noindent \textbf{3.2.13.~Remark.~}The proof of the Montague Reflection
Theorem shows a slightly stronger result than the statement of Theorem
3.2.11, since it shows that for each formula\textit{\ }$\varphi (\vec{x}),$
there is a formula $\theta _{\varphi }(y)$ such that for \textit{all
subformulae} $\psi $ \textit{of} $\varphi $ we have $\mathsf{ZF}\vdash 
\mathrm{Ref}_{\psi ,\theta }.$\bigskip \bigskip

\begin{center}
\textbf{4.~}$\mathsf{CT}^{-}\mathbf{[}\mathsf{ZF}\mathbf{]}$ \textbf{AND }$%
\mathsf{CT}\mathbf{[}\mathsf{ZF}\mathbf{]}$\bigskip
\end{center}

\noindent In this section we review known results about the two most
`famous' Tarskian theories of truth over $\mathsf{ZF}$, namely $\mathsf{CT}%
^{-}[\mathsf{ZF}],$ and its strengthening $\mathsf{CT}[\mathsf{ZF}]$.
\medskip

\noindent \textbf{4.1.~Definition}.~Recall that $\mathcal{L}_{\mathrm{set}}(%
\mathsf{T})$ is the extension of $\mathcal{L}_{\mathrm{set}}$ with a fresh
unary predicate $\mathsf{T}(x)$. For unexplained notation in the items
below, see Definition 2.1$\medskip $

\begin{enumerate}
\item[\textbf{(a)}] $\mathsf{Int}$\textsf{-}$\mathsf{Sep}$ (internal
separation) is the $\mathcal{L}_{\mathrm{set}}(\mathsf{T})$-sentence that
asserts that every instance of the separation scheme is true, i.e.,%
\begin{equation*}
\forall \varphi (v,x)\in \mathrm{Form}\ \mathsf{T}(\forall v\,\mathrm{Sep}%
_{\varphi }).
\end{equation*}

\item[\textbf{(b)}] $\mathsf{Int}$\textsf{-}$\mathsf{Coll}$ (internal
collection) is the $\mathcal{L}_{\mathrm{set}}(\mathsf{T})$-sentence that
asserts that every instance of the collection scheme is true, i.e., 
\begin{equation*}
\forall \varphi (v,x,y)\in \mathrm{Form}\ \mathsf{T}(\forall v\,\mathrm{Coll}%
_{\varphi }).
\end{equation*}

\item[\textbf{(c)}] $\mathsf{Int}$\textsf{-}$\mathsf{Repl}$ (internal
replacement) is the $\mathcal{L}_{\mathrm{set}}(\mathsf{T})$-sentence that
asserts that every instance of the replacement scheme is true, i.e.,%
\begin{equation*}
\forall \varphi (v,x,y)\in \mathrm{Form}\ \mathsf{T}(\forall v\,\mathrm{Repl}%
_{\varphi }).
\end{equation*}

\item[\textbf{(d)}] $\mathsf{Int}$\textsf{-}$\mathsf{Ind}$\textsf{\ }%
(internal induction) is the $\mathcal{L}_{\mathrm{set}}(\mathsf{T})$%
-sentence that asserts that every instance of the induction scheme is true,
i.e.,%
\begin{equation*}
\forall \varphi (v,x)\in \mathrm{Form}\ \mathsf{T}(\forall v\,\mathrm{Ind}%
_{\varphi }).
\end{equation*}

\item[\textbf{(e)}] $\mathsf{Ind(T)}$ is the full scheme of induction over $%
\omega $ in the language $\mathcal{L}_{\mathrm{set}}(\mathsf{T}).\medskip $

\item[\textbf{(f)}] $\mathsf{Sep(T)}$ is the full scheme of separation in
the language $\mathcal{L}_{\mathrm{set}}(\mathsf{T})$.$\medskip $

\item[\textbf{(g)}] $\mathsf{Coll(T)}$ is the full scheme of collection in
the language $\mathcal{L}_{\mathrm{set}}(\mathsf{T})$.$\medskip $

\item[\textbf{(h)}] $\mathsf{CT}[\mathsf{ZF}]:=\mathsf{CT}^{-}[\mathsf{ZF}]+%
\mathsf{Sep(T)+Coll(T)}.$\footnote{%
In Fujimoto's paper \cite{Fujimoto-APAL}, $\mathsf{CT}[\mathsf{ZF}]$ is
named $\mathsf{TC.}$ Also note that $\mathsf{CT}[\mathsf{ZF}]$ is
definitionally equivalent to $\mathsf{CS}[\mathsf{ZF}]:=\mathsf{CS}^{-}[%
\mathsf{ZF}]+\mathsf{Sep(S)+Coll(S).}$}$\medskip $
\end{enumerate}

\noindent \textbf{4.2.}~\textbf{Remark.}~In the presence of $\mathsf{CT}^{-}[%
\mathsf{ZF}]$, $\mathsf{Int}$\textsf{-}$\mathsf{Repl}$ is equivalent to the
conjunction of $\mathsf{Int}$\textsf{-}$\mathsf{Sep}$\textsf{\ }and $\mathsf{%
Int}$\textsf{-}$\mathsf{Coll}$\textsf{. }This can be readily verified by a
minor variant of the usual proof of equivalence of the replacement scheme
with the union of the schemes of separation and collection (in the presence
of the finitely many remaining axioms of $\mathsf{ZF}$). Thus, in the
presence of $\mathsf{CT}^{-}[\mathsf{ZF}],$ each of the sentences $\mathsf{%
Int}$\textsf{-}$\mathsf{Repl}$ and $\mathsf{Int}$-$\mathsf{Sep\wedge Int}$-$%
\mathsf{Coll}$\textsf{\ }is equivalent to the sentence that expresses
\textquotedblleft all the \textit{axioms} of the usual axiomatization of $%
\mathsf{ZF}$ are true\textquotedblright . Indeed $\mathsf{CT}^{-}[\mathsf{ZF}%
_{0}]$ suffices for this purpose, where $\mathsf{ZF}_{0}$ is the result of
deleting the separation and collection schemes from our formulation of the
axioms of $\mathsf{ZF}$ (in Definition 3.2.4). $\medskip $

\noindent Krajewski \cite{Krajewski} used the Montague-Vaught Reflection
Theorem to show that $\mathsf{CT}^{-}[\mathsf{ZF}]$ is conservative over $%
\mathsf{ZF}$; a close examination of his proof reveals the following
stronger result, as independently noted in \cite[Theorem 20]{Fujimoto-APAL}
and \cite{Ali+Albert-long}. We will see in Section 10 that $\mathsf{CT}^{-}[%
\mathsf{ZF}]+\mathsf{Coll(T)}$ is also conservative over $\mathsf{ZF}$%
.\medskip

\noindent \textbf{4.3}.~\textbf{Theorem.}~$\mathsf{CT}^{-}[\mathsf{ZF}]+%
\mathsf{Sep(T)}$ \textit{is conservative over} $\mathsf{ZF}$.\medskip

\noindent \textbf{4.4}.~\textbf{Definition.}~$\mathsf{GRef}_{\mathsf{T}}$
(the \textit{global reflection principle over }$\mathsf{T}$) is the $%
\mathcal{L}_{\mathrm{set}}(\mathsf{T})$-sentence that asserts that $\mathsf{T%
}$ is closed under first order proofs. More formally: 
\begin{equation*}
\mathsf{GRef}_{\mathsf{T}}:=\forall \varphi \in \mathrm{Sent}^{+}\,\left( 
\mathrm{Prov}_{\mathsf{T}}(\varphi )\rightarrow \mathsf{T}(\varphi )\right) ,
\end{equation*}%
where $\mathrm{Prov}_{\mathsf{T}}(\varphi )$ is the $\mathcal{L}_{\mathrm{set%
}}(\mathsf{T})$-sentence that expresses \textquotedblleft there is a proof
of $\varphi $ from premises in $\mathsf{T}$\textquotedblright , and $\varphi
\in \mathrm{Sent}^{+}$ expresses \textquotedblleft $\varphi $ is a sentence
in the language obtained by adding the proper class of constants $\left\{ 
\dot{a}:a\in \mathrm{V}\right\} $ to $\mathcal{L}_{\mathrm{set}}$%
\textquotedblright\ (as in Definition 3.2.4(b))\medskip

\noindent The following corollary is in sharp contrast with the fact that $%
\mathrm{Con}\mathsf{(PA)}$ is provable in $\mathsf{CT}^{-}\mathsf{[PA]}$ +
\textquotedblleft $\mathsf{T}$ is closed under first order
proofs\textquotedblright ; see \cite{Cieslinski-book}$\mathsf{.}$\medskip

\noindent \textbf{4.5}.~\textbf{Corollary.}~\textit{The following theories
are conservative over} $\mathsf{ZF}$:\medskip

\begin{enumerate}
\item[$(a)$] $\mathsf{CT}^{-}\mathsf{[ZF]}$ $+\ \mathsf{Ind}(\mathsf{T).}$%
\medskip

\item[$(b)$] $\mathsf{CT}^{-}\mathsf{[ZF]}$ + $\mathsf{GRef}_{\mathsf{T}}$%
.\medskip
\end{enumerate}

\noindent \textbf{Proof.}~$(a)$ follows from Theorem 4.3, since $\mathsf{%
Ind(T)}$ is provable in $\mathsf{CT}^{-}[\mathsf{ZF}]+\mathsf{Sep(T).}$ $(b)$
readily follows from $(a)$ by an induction on lengths of proofs.

\hfill $\square $\medskip

\noindent The following conservativity result is a special case (for set
theory) of a general result established in \cite{Ali+Albert-long} ; it
generalizes the conservativity of $\mathsf{CT}^{-}[\mathsf{PA}]+\mathsf{Int}$%
-$\mathsf{Ind}$ over $\mathsf{PA}$, first established in \cite{Kotlarski et
al}.\medskip 

\noindent \textbf{4.6}.~\textbf{Theorem.}~\textit{If }$\mathsf{B\supseteq KP}
$\textit{, and} $\mathsf{S}$ \textit{is a scheme all of whose instances are
provable in} $\mathsf{B}$, \textit{then} $\mathsf{CT}^{-}[\mathsf{B}]$ + $%
\mathsf{Int}$-$\mathsf{S}$ \textit{is conservative over} $\mathsf{B}$. 
\textit{Here }$\mathsf{Int}$-$\mathsf{S}$ \textit{is the} $\mathcal{L}_{%
\mathrm{set}}(\mathsf{T})$-\textit{sentence that asserts that every instance
of }$\mathsf{S}$\textit{\ is true.}\footnote{%
Given a language $\mathcal{L}$, an $\mathcal{L}$-\textit{template} for a
scheme $\mathsf{S}$ is given by a sentence $\tau (P)$ formulated in the
language obtained by augmenting $\mathcal{L}$ with an $n$-ary predicate $%
P(x_{1},...,x_{n})$. A sentence $\psi $ is then said to be \emph{an instance
of }$\mathsf{S}$ if $\psi $ is of the form $\forall v~\tau \lbrack \varphi
(v,x_{1},\cdot \cdot \cdot ,x_{n})/P]$, where $\tau \lbrack \varphi
(v,x_{1},\cdot \cdot \cdot ,x_{n})/P]$ is the result of substituting all
subformulae of the form $P(t_{1},\cdot \cdot \cdot ,t_{n})$, where each $%
t_{i}$ is a term, with $\varphi (v,t_{1},t_{2},\cdot \cdot \cdot ,t_{n})$
(and re-naming bound variables of $\varphi $ to avoid unintended clashes).
For more detail, and related results, see \cite{Enayat+Lelyk-JSL}.}\medskip 

\noindent \textbf{4.7}.~\textbf{Corollary.}~$\mathsf{CT}^{-}[\mathsf{ZF}]~+$ 
$\mathsf{Int}$\textsf{-}$\mathsf{Repl}$ \textit{is conservative over} $%
\mathsf{ZF}$.\medskip

\noindent \textbf{4.8}.~\textbf{Remark.}~Even though each of the theories $%
\mathsf{CT}^{-}[\mathsf{ZF}]+\mathsf{Sep(T)}$ and $\mathsf{CT}^{-}[\mathsf{ZF%
}]+\mathsf{Int}$\textsf{-}$\mathsf{Repl}$ is conservative over $\mathsf{ZF}$
(by Corollary 4.3 and Corollary 4.7), in light of Corollary 4.5 their union
implies $\mathrm{Con}\mathsf{(ZF)}$, and is thus not conservative over $%
\mathsf{ZF}$. \medskip

\noindent The next results shows that the two conservative theories $\mathsf{%
CT}^{-}[\mathsf{ZF}]+\mathsf{Sep(T)}$ and $\mathsf{CT}^{-}[\mathsf{ZF}]+%
\mathsf{Int}$\textsf{-}$\mathsf{Repl}$ behave differently with respect to
interpretability in $\mathsf{ZF}$.\medskip

\noindent \textbf{4.9}.~\textbf{Theorem.}~$\mathsf{CT}^{-}[\mathsf{ZF}]+%
\mathsf{Sep(T)}$ \textit{is interpretable in} $\mathsf{ZF}$, \textit{but }$%
\mathsf{CT}^{-}[\mathsf{ZF}]+\mathsf{Int}$\textsf{-}$\mathsf{Repl}$ \textit{%
is not interpretable in} $\mathsf{ZF}$.\medskip

\noindent \textbf{Proof.}~An inspection of the proof of Theorem 4.3 shows
that $\mathsf{CT}^{-}[\mathsf{ZF}]+\mathsf{Sep(T)}$ is locally interpretable
in $\mathsf{ZF.}$\footnote{%
This observation implies that $\mathsf{CT}^{-}[\mathsf{ZF}]$ is not finitely
axiomatizable.} Since $\mathsf{ZF}$ is a reflective theory (i.e., proves the
formal consistency of each of its finite subtheories), the global
interpretability of $\mathsf{CT}^{-}[\mathsf{ZF}]+\mathsf{Sep(T)}$ in $%
\mathsf{ZF}$ then follows by Orey's compactness theorem. In light of the
well-known fact that $\mathsf{GB}$ is not interpretable in $\mathsf{ZF}$,
the failure of interpretability of $\mathsf{CT}^{-}[\mathsf{ZF}]~+$ $\mathsf{%
Int}$\textsf{-}$\mathsf{Repl}$ in $\mathsf{ZF}$ follows from Lemma 7.9\
(which implies that $\mathsf{GB}$ is interpretable in $\mathsf{CT}^{-}[%
\mathsf{ZF}]~+$ $\mathsf{Int}$\textsf{-}$\mathsf{Repl}$).

\hfill $\square $\medskip

\noindent \textbf{4.10}.~\textbf{Remark.}~The interpretation of $\mathsf{CT}%
^{-}[\mathsf{ZF}]+\mathsf{Sep(T)}$ in $\mathsf{ZF}$ can be readily shown to
be a \textit{polynomial interpretation} (in the sense of \cite[Definition
2.4.4(2)]{Trio on feasible red.}). This implies that (1) $\mathsf{CT}^{-}[%
\mathsf{ZF}]+\mathsf{Sep(T)}$ has at most polynomial speed-up over $\mathsf{%
ZF}$,\textsf{\ }and (2) the conservativity of $\mathsf{CT}^{-}[\mathsf{ZF}]+%
\mathsf{Sep(T)}$ over $\mathsf{ZF}$ can be verified in $\mathrm{I}\Delta
_{0}+\mathsf{Exp.}$ In contrast, using the observation that $\mathsf{CT}^{-}[%
\mathsf{ZF}]~+$ $\mathsf{Int}$\textsf{-}$\mathsf{Repl}$ is deductively
equivalent to the finitely axiomatized theory $\mathsf{CT}^{-}[\mathsf{KP}]~+
$ $\mathsf{Int}$\textsf{-}$\mathsf{Repl}$\textsf{, }usual methods (as in 
\cite[Corollary 8]{Fischer-speed-up}) show that $\mathsf{CT}^{-}[\mathsf{ZF}%
]~+$ $\mathsf{Int}$\textsf{-}$\mathsf{Repl}$ has superexponential speed-up
over $\mathsf{ZF}$, and therefore the conservativity of $\mathsf{CT}^{-}[%
\mathsf{ZF}]~+$ $\mathsf{Int}$\textsf{-}$\mathsf{Repl}$ over $\mathsf{ZF}$
cannot be established in $\mathrm{I}\Delta _{0}+\mathsf{Exp.}$\footnote{%
Using Leigh's methodology in \cite{Graham cons.}, one should be able to show
that the conservativity of $\mathsf{CT}^{-}[\mathsf{ZF}]~+$ $\mathsf{Int}$%
\textsf{-}$\mathsf{Repl}$ over $\mathsf{ZF}$ to be verifiable in $\mathrm{I}%
\Delta _{0}+\mathsf{Supexp.}$} \medskip 

\noindent \textbf{4.11.~Definition.}~For $n\in \mathbb{N},$ $\mathsf{\Sigma }%
_{n}^{1}$-$\mathsf{Sep}$ is the scheme of separation for $\mathsf{\Sigma }%
_{n}^{1}$-formulae; and $\mathsf{\Sigma }_{n}^{1}$-$\mathsf{Coll}$ is the
scheme of collection for $\mathsf{\Sigma }_{n}^{1}$-formula. The full
separation scheme in this context will be denoted by $\mathsf{\Sigma }%
_{\infty }^{1}$-$\mathsf{Sep}$, and the full collection scheme will be
denoted by $\mathsf{\Sigma }_{\infty }^{1}$-$\mathsf{Coll}.$\medskip 

\noindent The following result, due to Fujimoto \cite[Theorem 20]%
{Fujimoto-APAL}, is the set-theoretic analogue of the well-known mutual $%
\omega $-interpretability of $\mathsf{ACA}$ and $\mathsf{CT}[\mathsf{PA}]$%
.\medskip

\noindent \textbf{4.12.~Theorem.}~$\mathsf{CT}[\mathsf{ZF}]$ \textit{and} $%
\mathsf{GB+\Sigma }_{\infty }^{1}$-$\mathsf{Sep}+\mathsf{\Sigma }_{\infty
}^{1}$-$\mathsf{Coll}$ \textit{are mutually} $\mathrm{V}$-\textit{%
interpretable, i.e., they can be interpreted in each other via
interpretations that are the identity on the class of }$\mathrm{V}$\textit{\
of all sets. Consequently they have the same }$\mathcal{L}_{set}$-\textit{%
consequences}.\medskip

\noindent \textbf{4.13.~Remark. }A proof\footnote{%
The proof is similar to the well-known argument in the context of arithmetic
that shows that in the presence of $\mathrm{I}\Delta _{0}$ and the
collection scheme, each instance of the induction scheme is derivable.}
based on the complexity of formulae shows that:\textbf{\ }\medskip

\begin{center}
$\mathsf{GB}+\mathsf{\Delta }_{0}^{0}$-$\mathsf{Sep}+\mathsf{\Sigma }%
_{\infty }^{1}$-$\mathsf{Coll}\vdash \mathsf{\Sigma }_{\infty }^{1}$-$%
\mathsf{Sep.}$ \medskip
\end{center}

\noindent A similar proof shows that $\mathsf{CT}^{-}[\mathsf{ZF}]+\mathsf{%
\Delta }_{0}$-$\mathsf{Sep(T)}+\mathsf{\Delta }_{0}$-$\mathsf{Coll(T)\vdash
Sep(T).}$\medskip

\noindent \textbf{4.14.~Definition.}~For $n\in \omega $, $\mathsf{CT}_{n}[%
\mathsf{ZF}]:=\mathsf{CT}^{-}[\mathsf{ZF}]+\mathsf{\Sigma }_{n}$-$\mathsf{%
Sep(T)+\mathsf{\Sigma }_{n}}$-$\mathsf{Coll(T).}$\medskip

\noindent \textbf{4.15.~Remark.}~$\mathsf{CT}_{n+2}$ $[\mathsf{ZF}]$ proves
the consistency of $\mathsf{CT}_{n}[\mathsf{ZF}]+\mathsf{Sep(T)}$ for\textit{%
\ }$n\geq 1$. This can be established by a straightforward adaptation of the
proof of Theorem 4.6 of McKenzie's \cite{Zach collection}.\medskip

\noindent \textbf{4.16.~Definition.}~$\mathsf{FRef}$ (Full Reflection) is
the $\mathcal{L}_{\mathrm{set}}(\mathsf{T})$-sentence:%
\begin{equation*}
\forall \alpha _{0}\in \mathrm{Ord}\ \exists \alpha \in \mathrm{Ord\,}\left[ 
\mathrm{(}\alpha _{0}<\alpha )\wedge \forall \varphi \in \mathrm{Form\ }%
\mathsf{T}\left( \mathrm{V\!}_{\alpha }\prec _{\varphi }\mathrm{V}\right) %
\right] ,
\end{equation*}

\noindent where $\mathrm{V\!}_{\alpha }\prec _{\varphi }\mathrm{V}$ is
shorthand for $\forall \vec{x}\in \mathrm{V\!}_{\alpha }\ \left[ \varphi
\left( \vec{x}\right) \longleftrightarrow \varphi ^{\mathrm{V\!}_{\alpha }}(%
\vec{x})\right] .$\footnote{%
It is implicit in this notation that $\vec{x}$ lists the free variables of $%
\varphi $.}\medskip 

\noindent The following result was first established by Montague and Vaught,
who formulated their result in terms of the definitionally equivalent theory 
$\mathsf{CS}[\mathsf{ZF}]$ instead of $\mathsf{CT}[\mathsf{ZF}].$\footnote{%
Note that in \cite[Theorem 7.1]{Montague-Vaught}, our $\mathsf{ZF}$ is
denoted $\mathsf{ZFS}$\ ($\mathsf{S}$ for `set theory'), and our $\mathsf{%
CS[ZF]}$ is denoted $\mathsf{ZFS}^{\prime }$.} This result was revisited by
Fujimoto \cite[Theorem 23]{Fujimoto-APAL}.\medskip 

\noindent \textbf{4.17.~Theorem.}~(Montague and Vaught) $\mathsf{CT}[\mathsf{%
ZF}]\vdash \mathsf{FRef.}$\bigskip \bigskip

\begin{center}
\textbf{5.~FULL REFLECTION IN }$\mathsf{CT}_{0}[\mathsf{ZF}]$\bigskip
\end{center}

\noindent Recall from Definition 4.14 that $\mathsf{CT}_{0}[\mathsf{ZF}]:=%
\mathsf{CT}^{-}[\mathsf{ZF}]+\mathsf{\Delta }_{0}$-$\mathsf{Sep(T)+\mathsf{%
\Delta }_{0}}$-$\mathsf{Coll(T).}$ In this section we fine-tune Theorem 4.17
by showing that $\mathsf{FRef}$, and its iterations, are provable in $%
\mathsf{CT}_{0}[\mathsf{ZF}].$ For this purpose we have the occasion to
introduce the weaker theory $\mathsf{CT}_{\ast }[\mathsf{ZF}]$. We will
further explore $\mathsf{CT}_{\ast }[\mathsf{ZF}]$ in Sections 6, 7, and 8%
\textsf{.}\medskip

\noindent \textbf{5.1.~Definition.}~Let $\mathrm{Prov}_{\mathrm{ZF}}(x)$ be
the $\mathcal{L}_{\mathrm{set}}$-formula that expresses \textquotedblleft $x$
is provable in $\mathsf{ZF}$\textquotedblright , and let 
\begin{equation*}
\mathsf{CT}_{\ast }[\mathsf{ZF}]:=\mathsf{CT}^{-}[\mathsf{ZF}]+\mathrm{GRef}%
_{\mathrm{ZF}},
\end{equation*}%
where:\medskip

\begin{center}
$\mathrm{GRef}_{\mathrm{ZF}}:=\forall x\in \mathrm{Sent}^{+}(\mathrm{Prov}_{%
\mathrm{ZF}}(x)\rightarrow \mathsf{T}(x)).$\footnote{%
Recall from Definition 3.2.4(b) that $\mathrm{Sent}^{+}$ is the proper class
of sentences of the language obtained by enriching $\mathcal{L}_{\mathrm{set}%
}$ with a constant symbol for each element of the universe. Notice that
since the language $\mathcal{L}_{\mathrm{set}}$ in which $\mathsf{ZF}$ is
formulated has no constant symbols, $\varphi \left( c_{1},\cdot \cdot \cdot
,c_{k}\right) $ is provable in $\mathsf{ZF}$, where each $c_{i}$ is a
constant symbol, iff $\forall x_{1}\cdot \cdot \cdot \forall x_{k}\,\varphi
\left( x_{1},\cdot \cdot \cdot ,x_{k}\right) $ is provable in $\mathsf{ZF}.$}%
\medskip\ 
\end{center}

\noindent \textbf{5.2.~Proposition.}~$\mathsf{CT}_{\ast }[\mathsf{ZF}]\vdash
\sigma $, where:\medskip

\begin{center}
$\sigma :=\left[ 
\begin{array}{c}
\forall m\ \forall k\in \omega \ \forall \psi (\vec{v})\in \mathrm{Form}%
_{k}\ \forall x_{0}\in m\cdot \cdot \cdot \forall x_{k-1}\in m\smallskip  \\ 
\left( \psi ^{m}\left( x_{0},\cdot \cdot \cdot ,x_{k-1}\right) \in \mathsf{T}%
\leftrightarrow \left( \psi \left( \dot{x}_{0},\cdot \cdot \cdot ,\dot{x}%
_{k-1}\right) \in \mathrm{ED}(m,\in )\right) \right) 
\end{array}%
\right] .$\footnote{%
\noindent Recall that within $\mathsf{ZF}$, given a set $m$, $\mathrm{ED}%
(m,\in )$ is the elementary diagram of $(m,\in ).$}\medskip 
\end{center}

\noindent \textbf{Proof.}~This follows from $\mathrm{GRef}_{\mathsf{ZF}}$
once we observe that following holds, which is a formalization of the
assertion in Remark 3.2.12. 
\begin{equation*}
\mathsf{ZF}\vdash \forall m\,\forall k\in \omega \,\forall \psi (\vec{v})\in 
\mathrm{Form}_{k}\,\mathrm{Prov}_{\mathrm{ZF}}\left( \pi _{\psi }\right) ,
\end{equation*}%
where:$\mathsf{\ }$%
\begin{equation*}
\pi _{\psi }:=\forall x_{0}\in m\cdot \cdot \cdot \forall x_{k-1}\in m\ 
\left[ \left( \psi ^{m}\left( x_{0},\cdot \cdot \cdot ,x_{k-1}\right)
\leftrightarrow \left( \psi \left( \dot{x}_{0},\cdot \cdot \cdot ,\dot{x}%
_{k-1}\right) \in \mathrm{ED}(m,\in )\right) \right) \right] .
\end{equation*}

\hfill $\square $\medskip

\noindent \textbf{5.3.} \textbf{Definition.}~$\mathsf{Int}$-$\mathsf{Ref}$
(internal reflection) is the following $\mathcal{L}_{\mathrm{set}}(\mathsf{T}%
)$-sentence: 
\begin{equation*}
\forall k\in \omega \ \forall \varphi \in \mathrm{Form}_{k}\,\exists \theta
\in \mathrm{Form}_{1}\,\mathsf{T}\left( \mathrm{Ref}_{\varphi ,\theta
}\right) ,
\end{equation*}%
where $\mathrm{Ref}_{\varphi ,\theta }$ is as in Reflection Theorem
3.2.11.\medskip

\noindent \textbf{5.4.} \textbf{Theorem.}~$\mathsf{CT}_{\ast }[\mathsf{ZF}%
]\vdash \mathsf{Int}$-$\mathsf{Ref}$.\medskip

\noindent \textbf{Proof.}~Observe that Theorem 3.2.11 is a syntactic fact
that is already provable in modest arithmetical theories, and \textit{a
fortiori} we have: 
\begin{equation*}
\mathsf{ZF}\vdash \forall k\in \omega \,\forall \varphi \in \mathrm{Form}%
_{k}\,\exists \theta \in \mathrm{Form}_{1}\,\mathsf{T}\left( \mathrm{Ref}%
_{\varphi ,\theta }\right) .
\end{equation*}%
Hence the result follows from $\mathrm{GRef}_{\mathsf{ZF}}$.

\hfill $\square $\medskip

\noindent \textbf{5.5.~Remark.}~We will see in Theorem 6.12 that a weaker
form of internal reflection, implies the strong form $\mathsf{Int}$-$\mathsf{%
Ref}$ within $\mathsf{CT}^{-}[\mathsf{ZF}]$. \medskip

\noindent The rest of this section is focused on $\mathsf{CT}_{0}[\mathsf{ZF}%
].$ We begin with the following observation.\medskip

\noindent \textbf{5.6.~Remark.}~Recall that $\mathsf{\Delta }_{0}$-$\mathsf{%
Coll(T)}$ is equivalent to $\mathsf{\Sigma }_{1}$-$\mathsf{Coll(T)}$, with
the usual trick of collapsing two consecutive existential quantifiers into a
single one with the help of Kuratowski pairing function. This fact enables
the theory $\mathsf{CT}_{0}[\mathsf{ZF}]$ to prove the existence of
functions that are defined by $\Sigma _{1}(\mathsf{T})$-recursions. In
particular, if $(\mathcal{M},T)\models \mathsf{CT}_{0}[\mathsf{ZF}],$ and 
\begin{equation*}
G:\mathrm{Ord}^{\mathcal{M}}\times \omega ^{\mathcal{M}}\rightarrow M
\end{equation*}%
is a function whose graph is $\Sigma _{1}(\mathsf{T})$-definable in $(%
\mathcal{M},T)$ (parameters allowed), then for any given $\alpha _{0}\in 
\mathrm{Ord}^{\mathcal{M}}$ there is some $F\in M$ such that $\mathcal{M}%
\models \left[ F:\omega \rightarrow \mathrm{V}\right] $, and:%
\begin{equation*}
\mathcal{M}\models \left[ F(0)=\alpha _{0}\wedge \forall k\in \omega \
F(k+1)=G(F(k),k)\right] .
\end{equation*}%
Recall from Theorem 4.17 that $\mathsf{CT}[\mathsf{ZF}]\vdash \mathsf{FRef}$%
. We refine this result in the next theorem.\medskip

\noindent \textbf{5.7.~Theorem.}~$\mathsf{CT}_{0}[\mathsf{ZF}]\vdash \mathsf{%
FRef}.$\medskip

\noindent \textbf{Proof.}~Suppose $(\mathcal{M},T)\models \mathsf{CT}_{0}[%
\mathsf{ZF}]$. It is easy to see that there is an enumeration $\left\langle
\varphi _{i}:i\in \omega ^{\mathcal{M}}\right\rangle $ of all $\mathcal{L}_{%
\mathrm{set}}$-formulae within $\mathcal{M}$ such that if $\varphi _{i}$ is
a subformula of $\varphi _{j}$ then $i<j$, and furthermore, the map $%
i\mapsto \varphi _{i}$ is $\Delta _{1}$ in $\mathcal{M}$. Fix some $\alpha
_{0}\in \mathrm{Ord}^{\mathcal{M}}$. Since by Theorem 5.4, $\mathrm{GRef}_{%
\mathsf{ZF}}$ is provable in $\mathsf{CT}_{0}[\mathsf{ZF}]$, there is a
function $G:\mathrm{Ord}^{\mathcal{M}}\times \omega ^{\mathcal{M}%
}\rightarrow \mathrm{Ord}^{\mathcal{M}}$ such that:%
\begin{equation*}
G(\delta ,k)=\beta _{0}\Leftrightarrow (\mathcal{M},T)\models \left[ \beta
_{0}=\min \left\{ \beta \in \mathrm{Ord}:\left( \beta >\delta \right) \wedge
\forall i\leq k\,\mathsf{T}\left( \mathrm{V\!}_{\beta }\prec _{\varphi _{i}}%
\mathrm{V}\right) \right\} \right] .
\end{equation*}

\noindent Note that the graph of $G$ is $\mathsf{\Sigma }_{1}\mathsf{(T)}$%
-definable in $(\mathcal{M},T)\mathsf{.}$ Given $\alpha _{0}\in \mathrm{Ord}%
^{\mathcal{M}}$, we wish to show that there is some $F\in M$ such that $%
\mathcal{M}\models \left[ F:\omega \rightarrow \mathrm{Ord}\right] $ such
that: 
\begin{equation*}
(\mathcal{M},T)\models \left[ F(0)=\alpha _{0}\wedge \forall k\in \omega \
F(k+1)=G(F(k),k)\right] .
\end{equation*}

\noindent Therefore by Remark 5.6 the desired $F$ exists as a set in $%
\mathcal{M}$. Next let $\alpha \in \mathrm{Ord}^{\mathcal{M}}$ such that:%
\begin{equation*}
(\mathcal{M},T)\models \left[ \alpha =\bigcup\limits_{k\in \omega }F(k)%
\right] \mathsf{.}
\end{equation*}%
With the help of $\mathrm{GRef}_{\mathsf{ZF}}$, we can then verify that $%
\left( \mathcal{M},T\right) \models \left[ \forall \varphi \in \mathrm{Form\ 
}\mathsf{T}\left( \mathrm{V\!}_{\alpha }\prec _{\varphi }\mathrm{V}\right) %
\right] .$

\hfill $\square \medskip $

\noindent \textbf{5.8.~Definition.}~For each $n\in \mathbb{N}$, $\mathsf{FRef%
}^{n}$ is the $\mathcal{L}_{\mathrm{set}}(\mathsf{T})$-sentence recursively
defined by: $\medskip $

\begin{itemize}
\item $\mathsf{FRef}^{1}:=\mathsf{FRef,}$ and$\medskip $

\item $\mathsf{FRef}^{n+1}:=\forall \alpha _{0}\in \mathrm{Ord}\,\exists
\alpha \in \mathrm{Ord\,}\left[ \mathrm{(}\alpha _{0}<\alpha )\wedge \mathrm{%
V\!}_{\alpha }\prec \mathrm{V}\wedge \left( \mathsf{FRef}^{^{n}}\right)
^{\left( \mathrm{V\!}_{\alpha },\in ,\mathsf{T}\cap \mathrm{V\!}_{\alpha
}\right) }\right] .$
\end{itemize}

\noindent In the above, $\left( \mathsf{FRef}^{^{n}}\right) ^{\left( \mathrm{%
V\!}_{\alpha },\in ,\mathsf{T}\cap \mathrm{V\!}_{\alpha }\right) }$ is the
relativization of $\mathsf{FRef}^{n}$ to $\left( \mathrm{V\!}_{\alpha },\in ,%
\mathsf{T}\cap \mathrm{V\!}_{\alpha }\right) $; i.e., it is the result of
relativing all quantifiers in $\mathsf{FRef}^{n}$ to $\mathrm{V\!}_{\alpha }$%
, and every occurrence of $\mathsf{T}$ by $\mathsf{T}\cap \mathrm{V\!}%
_{\alpha }.$ Note that in light of Proposition 5.2, provably in $\mathsf{CT}%
_{\ast }[\mathsf{ZF}],$ if $\forall \varphi \in \mathrm{Form\ }\mathsf{T}%
\left( \mathrm{V\!}_{\alpha }\prec _{\varphi }\mathrm{V}\right) $, then $%
\mathsf{T}\cap \mathrm{V\!}_{\alpha }=\mathrm{ED}($\textrm{V}$_{\alpha },\in
)$\textrm{.}$\medskip $

\noindent \textbf{5.9.~Theorem.}~\textit{For each} $n\in \mathbb{N}$, $%
\mathsf{CT}_{0}[\mathsf{ZF}]\vdash \mathsf{FRef}^{n}.\medskip $

\noindent \textbf{Proof.}~We will use Theorem 5.7 to derive $\mathsf{FRef}%
^{2}$ in $\mathsf{CT}_{0}[\mathsf{ZF}].$ A similar inductive reasoning shows
that $\mathsf{CT}_{0}[\mathsf{ZF}]\vdash \mathsf{FRef}^{n}$ for all $n\in 
\mathbb{N}$. Suppose $(\mathcal{M},T)\models \mathsf{CT}_{0}[\mathsf{ZF}].$
Let $G_{1}:\mathrm{Ord}^{\mathcal{M}}\rightarrow \mathrm{Ord}^{\mathcal{M}}$
be defined by:%
\begin{equation*}
(\mathcal{M},T)\models \left[ G_{1}(\delta )=\min \left\{ \beta \in \mathrm{%
Ord}:\left( \beta >\delta \right) \wedge \forall \varphi \in \mathrm{Form}\ 
\mathsf{T}\left( \mathrm{V\!}_{\alpha }\prec _{\varphi }\mathrm{V}\right)
\right\} \right] .
\end{equation*}%
By Theorem 5.7, $F_{1}$ is well-defined. Note that the graph of $G_{1}$ is $%
\mathsf{\Sigma }_{1}\mathsf{(T)}$-definable in $(\mathcal{M},T).$ Therefore,
given $\alpha _{0}\in \mathrm{Ord}^{\mathcal{M}}$, there is some $F_{1}\in M$
such that $\mathcal{M}\models \left[ F_{1}:\omega \rightarrow \mathrm{Ord}%
\right] $ and $F_{1}$ satisfies the following: 
\begin{equation*}
(\mathcal{M},T)\models \left[ F_{1}(0)=\alpha _{0}\wedge \forall k\in \omega
\,F_{1}(k+1)=G_{1}(F_{1}(k))\right] .
\end{equation*}

\noindent Next let $\alpha \in \mathrm{Ord}^{\mathcal{M}}$ such that:%
\begin{equation*}
(\mathcal{M},T)\models \left[ \alpha =\bigcup\limits_{k\in \omega }F_{1}(k)%
\right] \mathsf{.}
\end{equation*}%
Usual arguments then show that $(\mathcal{M},T)\models \left[ \left( \alpha
_{0}<\alpha \right) \wedge \left( \mathrm{V\!}_{\alpha }\prec \mathrm{V}%
\right) \wedge \mathsf{FRef}^{\left( \mathrm{V\!}_{\alpha },\in ,\mathsf{T}%
\cap \mathrm{V\!}_{\alpha }\right) }\right] .$

\hfill $\square \medskip $

\noindent \textbf{5.10.~Remark.}~One can formulate sentences $\mathsf{FRef}%
^{\alpha }$ for appropriate transfinite $\alpha $, and then using a
reasoning similar to the proof of Theorem 5.9 shows that $\mathsf{CT}_{0}[%
\mathsf{ZF}]\vdash \mathsf{FRef}^{\alpha }.\bigskip $\bigskip 

\begin{center}
\textbf{6.~THE MANY FACES OF }$\mathsf{CT}_{\ast }[\mathsf{ZF}]$\bigskip
\end{center}

\noindent Recall from the previous section that $\mathsf{CT}_{\ast }[\mathsf{%
ZF}]=\mathsf{CT}^{-}[\mathsf{ZF}]+\mathsf{GRef}_{\mathrm{ZF}}.$ In this
section we establish various results concerning $\mathsf{CT}_{\ast }[\mathsf{%
ZF}]$ that culminate in Theorem 6.12 that reveals the `many faces' of this
theory. Along the way, we will also meet the well-behaved and much weaker
theory $\mathsf{CT}^{-}[\mathsf{ZF}]+\Delta _{0}^{\mathrm{fin}}$-$\mathsf{%
Ind(T)}$, which also exhibits a `many faces' feature, as in Theorem
6.10.\medskip 

The beginning material in this section pertains to extensions of $\mathsf{KP}
$ formulated in the language $\mathcal{L}_{\mathrm{set}}(\mathsf{P})$, where 
$\mathcal{L}_{\mathrm{set}}(\mathsf{P})$ is obtained by adding a unary
predicate $\mathsf{P}$ to $\mathcal{L}_{\mathrm{set}}.$ Of course our main
interest is when $\mathsf{P}$ is further assumed to behave like a truth
predicate, which we will turn to after dealing with some general results
concerning theories formulated in $\mathcal{L}_{\mathrm{set}}(\mathsf{P})$%
.\medskip

\noindent \textbf{6.1.~Definition.}~Suppose $\mathsf{P}$ is a unary
predicate. \medskip

\begin{enumerate}
\item[\textbf{(a)}] $\mathrm{Fin}(v)$ is shorthand for the $\mathcal{L}_{%
\mathrm{set}}$-formula that expresses \textquotedblleft $v$ is
finite\textquotedblright , i.e., $\mathrm{Fin}(v)=\left[ \exists k\in \omega
\ \left\vert v\right\vert =k\right] .$ Given an $\mathcal{L}_{\mathrm{set}}(%
\mathsf{P})$-formula $\varphi (v)$, $\forall ^{\mathrm{fin}}v\mathrm{\,}%
\varphi (v)$ is shorthand for $\forall v\ \left( \mathrm{Fin}(v)\rightarrow
\varphi (v)\right) .$ \medskip

\item[\textbf{(b)}] $\forall ^{\mathrm{fin}}s\mathrm{\,}\left( s\cap \mathsf{%
P}\in \mathrm{V}\right) $ is shorthand for $\forall ^{\mathrm{fin}}s\mathrm{%
\,}\left[ \exists x\,\overset{x\ =\ s\ \cap \ \mathsf{P}}{\overbrace{\left(
\forall y(y\in x\leftrightarrow \left( \left( y\in s\right) \wedge \mathsf{P}%
(y)\right) \right) }}\right] .$

\item[\textbf{(c)}] An $\mathcal{L}_{\mathrm{set}}(\mathsf{P})$-formula is $%
\Delta _{0}^{\mathrm{fin}},$ if it is of the form $\delta ^{\ast }$, where $%
\delta $ is a $\Delta _{0}$-formula in the language $\mathcal{L}_{\mathrm{set%
}}(\mathsf{P})$, and the operation $\ast $ is defined on $\Delta _{0}$%
-formulae in the language $\mathcal{L}_{\mathrm{set}}(\mathsf{P})$ by the
following inductive clauses. In what follows $x$ and $y$ range over
variables of first order logic.$\smallskip $

\begin{enumerate}
\item[(1)] $\left( x\in y\right) ^{\ast }=x\in y.\smallskip $

\item[(2)] $\left( x=y\right) ^{\ast }=x=y.\smallskip $

\item[(3)] $\left( \mathsf{P}\left( x\right) \right) ^{\ast }=\mathsf{P}%
\left( x\right) .\smallskip $

\item[(4)] $\left( \lnot \delta \right) ^{\ast }=\lnot \delta ^{\ast
}.\smallskip $

\item[(5)] $\left( \delta _{1}\vee \delta _{2}\right) ^{\ast }=\delta
_{1}^{\ast }\vee \delta _{2}^{\ast }.\smallskip $

\item[(6)] $\left( \exists x\in y\,\delta \right) ^{\ast }=\exists x\in
y\,\left( \mathrm{Fin}(y)\wedge \delta ^{\ast }\right) .$\footnote{%
As usual, we construe $\left( \exists x\in y\,\delta (x)\right) $ as an
abbreviation of $\exists x(x\in y\wedge \delta (x))$.}\medskip 
\end{enumerate}
\end{enumerate}

\noindent \textbf{6.2.~Remark. }The following lists some basic facts about $%
\Delta _{0}^{\mathrm{fin}}$-formulae that can be readily verified from first
principles.\medskip

\begin{enumerate}
\item[\textbf{(a)}] Recall that in our formulation of first order logic, $%
\left\{ \lnot ,\vee ,\exists \right\} $ are the only logical constants, thus 
$\forall x\,\delta (x)$ is understood to be an abbreviation of $\lnot
\exists x\,\lnot \delta (x).$ With this convention in mind, clause (7) of
part (c) of Definition 6.1 implies that $\left( \forall x\in y\,\delta
(x)\right) ^{\ast }$ can be written as $\forall x\in y\,\left( \mathrm{Fin}%
(y)\rightarrow \delta ^{\ast }(x)\right) $.

\item[\textbf{(b)}] If $\varphi $ and $\psi $ are equivalent to $\Delta
_{0}^{\mathrm{fin}}$-formulae, then so are $\varphi \vee \psi $, $\lnot
\varphi $, and $\exists x\in y\,\left( \mathrm{Fin}(y)\wedge \varphi \right)
.$\medskip

\item[\textbf{(c)}] Each of the following formulae is equivalent (over $%
\mathsf{KP}$) to a $\Delta _{0}^{\mathrm{fin}}$-formula:\medskip

\begin{enumerate}
\item[(1)] $\mathrm{Fin}(x)\wedge (y\in \cup x)$ $\smallskip $

\item[(2)] $\left\langle x,y\right\rangle =z$, where $\left\langle
x,y\right\rangle $ is the Kuratowski ordered pair of $x$ and $y.\smallskip $

\item[(3)] $\left\langle x,y\right\rangle \in z\wedge \mathrm{Fin}(z).$%
\medskip
\end{enumerate}
\end{enumerate}

\noindent \textbf{6.3.~Definition.}~Let $\mathcal{L}=\mathcal{L}_{\mathrm{set%
}}(\mathsf{P}),$ where $\mathsf{P}$ is a unary predicate. \medskip

\begin{enumerate}
\item[\textbf{(a)}] $\Delta _{0}^{\mathrm{fin}}$-$\mathsf{Ind(P)}$ consists
of $\mathcal{L}$-sentences of the following form, where $\varphi (v,x)$ is $%
\Delta _{0}^{\mathrm{fin}}(\mathsf{P)}$.\footnote{%
Using a pairing function, one can deduce the more general versions of the
schemes considered here, in which the single parameter $v$ is replaced by a
finite tuple $\vec{v}$ of parameters.}%
\begin{equation*}
\forall v\left[ \left( \varphi (v,0)\wedge \forall x\in \omega \,\left(
\varphi (v,x)\rightarrow \varphi (v,x+1)\right) \right) \rightarrow \forall
x\in \omega \,\varphi (v,x)\right] .
\end{equation*}

\item[\textbf{(b)}] $\Delta _{0}^{\mathrm{fin}}$-$\mathsf{Min(P)}$ consists
of $\mathcal{L}$-sentences of the following form, where $\varphi (v,x)$ is $%
\Delta _{0}^{\mathrm{fin}}(\mathsf{P)}.$%
\begin{equation*}
\forall v\left[ \left( \exists x\in \omega \ \varphi (v,x)\right)
\rightarrow \left( \exists x\in \omega \,\varphi (v,x)\wedge \forall y\in
\omega \,\left( y\in x\rightarrow \lnot \varphi (v,x)\right) \right) \right]
.
\end{equation*}

\item[\textbf{(c)}] $\Delta _{0}^{\mathrm{fin}}$-$\mathsf{Sep}(\mathsf{P})$
consists of $\mathcal{L}$-sentences of the following form, where $\varphi
(v,x)$ is $\Delta _{0}^{\mathrm{fin}}(\mathsf{P)}$. Note that in the formula
below, $w$ does not need not be a subset of $\omega .$%
\begin{equation*}
\forall v\left[ \forall ^{\mathrm{fin}}w\,\exists y\,\left[ \forall
x\,\left( x\in y\leftrightarrow \left( x\in w\wedge \varphi (v,x)\right)
\right) \right] \right] .
\end{equation*}
\end{enumerate}

\noindent \textbf{6.4.~Remark.~}As pointed out in Corollary 4.5, $\mathsf{CT}%
^{-}[\mathsf{ZF}]+\mathsf{Ind(T)}$ is conservative over $\mathsf{ZF}$.
However, there is an instance of $\Delta _{0}^{\mathrm{fin}}$\textbf{-}$%
\mathsf{Ind(T)}$ that is unprovable in $\mathsf{CT}^{-}[\mathsf{ZF}].$ This
can be readily demonstrated using the existence of `pathological' models $(%
\mathcal{M},T)$ of $\mathsf{CT}^{-}[\mathsf{ZF}]$ in which some sentence $%
\varphi $ is deemed true by $T$, and yet $\mathrm{D}(k,\varphi )$ is false
for some nonstandard $k\in \omega ^{\mathcal{M}}$ (or even for all
nonstandard $k\in \omega ^{\mathcal{M}}$)$,$ where $\mathrm{D}(k,\varphi )$
is defined inside $\mathcal{M}$ by the recursion via: 
\begin{equation*}
\mathrm{D}(1,\varphi ):=\left[ \varphi \vee \varphi \right] \text{;}\ \ \ 
\mathrm{D}(i+1,\varphi ):=\left[ \mathrm{D}(i,\varphi )\vee \ \mathrm{D}%
(i,\varphi )\right] .
\end{equation*}%
Such pathological models can be readily constructed using the EV-method.
Notice that models of $\mathsf{CT}^{-}[\mathsf{KP}]+\Delta _{0}^{\mathrm{fin}%
}$\textbf{-}$\mathsf{Ind(T)}$ do not exhibit such a pathology since we will
see in Theorem 6.10 that the closure of $\mathsf{T}$ under first order
deductions is provable in $\mathsf{CT}^{-}[\mathsf{KP}]+\Delta _{0}^{\mathrm{%
fin}}$\textbf{-}$\mathsf{Ind(T)}$. \medskip 

\noindent The next result (and its proof) sheds light on the expressive
power of $\Delta _{0}^{\mathrm{fin}}$-formulae.\medskip

\noindent \textbf{6.5.~Lemma.~}\textit{The following are equivalent in }$%
\mathsf{KP}$:\medskip

\begin{enumerate}
\item[$(a)$] $\Delta _{0}^{\mathrm{fin}}$-$\mathsf{Ind(P)}$.\medskip

\item[$(b)$] $\Delta _{0}^{\mathrm{fin}}$-$\mathsf{Min(P).}$\medskip

\item[$(c)$] $\Delta _{0}^{\mathrm{fin}}$-$\mathsf{Sep}(\mathsf{P})$.
\medskip

\item[$(d)$] $\forall ^{\mathrm{fin}}s\mathrm{\,}\left( s\cap \mathsf{P}\in 
\mathrm{V}\right) .$\medskip
\end{enumerate}

\noindent \textbf{Proof.}~The verification of $(a)\Leftrightarrow (b)$ is
straightforward, and is similar to the well-known equivalence of the
induction principle and the minimum principle in arithmetic for $\mathsf{%
\mathsf{\Delta }}\mathsf{_{0}}$-formulas of arithmetic (see, e.g., \cite[%
Lemma I.2.4]{Hajek and Pudlak}). Also note that $(c)\Rightarrow (d)$ is
trivial. \bigskip

\noindent The proof will be complete once we verify $(a)\Rightarrow (c)$, $%
(c)\Rightarrow (b)$, and $(d)\Rightarrow (c).$\medskip

\noindent $(a)\Rightarrow (c):$ Assume $\left( \mathcal{M},P\right) \models
\Delta _{0}^{\mathrm{fin}}$-$\mathsf{Ind(P)}$, where $\mathcal{M}$ is a
model of $\mathsf{KP}$. To show that $\left( \mathcal{M},P\right) \models
\Delta _{0}^{\mathrm{fin}}$-$\mathsf{Sep}(\mathsf{P})$, suppose $s\in M$,
and for some $k\in \omega ^{\mathcal{M}},$ $\mathcal{M}$ satisfies $%
\left\vert s\right\vert =k$. Thus we can fix some $f\in M$ such that $%
\mathcal{M}$ satisfies \textquotedblleft $f:k\rightarrow s$, and $f$ is a
bijection\textquotedblright . Let $p\in M$ such that $\mathcal{M}$ satisfies
\textquotedblleft $p$ is the powerset of $k$\textquotedblright . Note that $%
\mathcal{M}$ satisfies $\left\vert f\right\vert =k$ and $\left\vert
p\right\vert =2^{k}$, in particular $f$ and $p$ are finite in the sense of $%
\mathcal{M}$. Given some $\mathsf{\mathsf{\Delta }}_{0}^{\mathrm{fin}}(%
\mathsf{P})$ formula $\varphi (x)$ with suppressed parameters we need to
show: \medskip 

\begin{center}
\textquotedblleft $\left\{ x\in s:\varphi (x)\right\} $ exists%
\textquotedblright . \medskip
\end{center}

\noindent For this purpose, it suffices to show that $\{i\in k:\varphi
(f(i))\}$ exists, since if for some $a\in M,$

\begin{equation*}
\left( \mathcal{M},P\right) \models \left[ a=\{i\in k:\varphi (f(i))\}\right]
,
\end{equation*}%
then $\left( \mathcal{M},P\right) \models \left[ \{f(x):x\in a\}=\left\{
x\in s:\varphi (x)\right\} \right] .$ Consider formula $\theta (i)$ below,
whose free variable is $i$, and which has parameters $k$, $f$ and $p$.

\begin{equation*}
\theta (i):=\left[ \exists y\in p\overset{y\,=\,\left\{ j\in i\,:\,\varphi
(f(j))\right\} }{\overbrace{\,\forall j\ \left[ j\in y\leftrightarrow \left(
j\in i\wedge \varphi (f(j))\right) \right] }}\right] .
\end{equation*}%
\noindent Since $s$, $k$, $f$, and $p$ are finite in $\mathcal{M}$, $\theta
(i)$ is readily seen to be equivalent to a $\Delta _{0}^{\mathrm{fin}}$%
-formula in $\left( \mathcal{M},P\right) $. More explicitly, since $p$ is
finite, it suffices to check that the formula $\psi (y,i)$ expressing $%
y\,=\,\left\{ j\in i\,:\,\varphi (f(j))\right\} $ is equivalent to a $\Delta
_{0}^{\mathrm{fin}}$-formula in $\left( \mathcal{M},P\right) .$ $\psi (y,i)$
is the conjunction of $\psi _{1}(y,i)$ and $\psi _{2}(y,i),$ where $\psi
_{1}(y,i):=\forall j\ \left[ j\in y\rightarrow \left( j\in i\wedge \varphi
(f(j))\right) \right] $ and $\psi _{2}(y,i):=\forall j\ \left[ \left( j\in
i\wedge \varphi (f(j))\right) \rightarrow j\in y\right] .$\medskip

\noindent We first analyze $\psi _{1}(y,i).$ Since $y\in p=\mathcal{P}(k)$, $%
\psi _{1}(y,i)$ is equivalent to $\forall j\in y\ \left[ \mathrm{Fin}%
(y)\rightarrow \left( j\in i\wedge \varphi (f(j))\right) \right] ,$ so it
suffices to check that $\varphi (f(j))$ is equivalent to a $\Delta _{0}^{%
\mathrm{fin}}$-formula. Looking further in, $\varphi (f(j))$ is equivalent
to $\exists z\in s(\left\langle j,z\right\rangle \in f\wedge \varphi (z)),$
and $s$ is finite, so it remains to verify that $\left\langle
j,z\right\rangle \in f$ is equivalent to a $\Delta _{0}^{\mathrm{fin}}$%
-formula, which follows from the last item of Remark 6.2.\medskip

We now turn to $\psi _{2}(y,i)$. $\psi _{2}(y,i)$ can be written as $\forall
j\in i\ \left[ \varphi (f(j))\rightarrow j\in y\right] .$ Note that $i$
ranges over elements of the finite ordinal $k$, so $\psi _{2}(y,i)$ is
equivalent $\forall j\in i\ \left[ \mathrm{Fin}(i)\rightarrow \left( \varphi
(f(j))\rightarrow j\in y\right) \right] .$ Since we already verified that $%
\varphi (f(j))$ is equivalent to $\Delta _{0}^{\mathrm{fin}}$-formula, this
concludes our verification for $\psi _{2}(y,i).$\medskip

Clearly $\theta (0)$ holds in $\left( \mathcal{M},P\right) $. $\left( 
\mathcal{M},P\right) $ also satisfies.$\forall i<k-1\ \left[ \theta
(i)\rightarrow \theta (i+1)\right] $, since given sets $a$ and $b$, $\mathsf{%
KP}$ proves that $a\cup \{b\}$ exists. Thus by $\Delta _{0}^{\mathrm{fin}}$-$%
\mathsf{Ind(P)}$, $\theta (k)$ holds, as desired.\medskip

\noindent $(c)\Rightarrow (b)$. Assume $\left( \mathcal{M},P\right) \models
\Delta _{0}^{\mathrm{fin}}$-$\mathsf{Sep}(\mathsf{P})$. Suppose that $\delta
(x,\mathsf{P})$ is a $\Delta _{0}^{\mathrm{fin}}\mathsf{(P)}$-formula and $%
\delta (m,\mathsf{P})$ holds in $\left( \mathcal{M},P\right) $ for some $%
m\in \omega ^{\mathcal{M}}.$ Then, as viewed by $\left( \mathcal{M},P\right) 
$, $\{i\in m+1:\delta (i,\mathsf{P})\}$ is coded by some $y$ by the veracity
of $\Delta _{0}^{\mathrm{fin}}$-$\mathsf{Sep}(\mathsf{P})$ in $\left( 
\mathcal{M},P\right) $. On the other hand, $\mathsf{KP}$ proves that a
non-empty subset of an ordinal has a $\in $-minimum element, so $y$ has a
minimum member, as desired.\medskip 

\noindent $(d)\Rightarrow (c).$\footnote{%
This is the set-theoretical analogue of \cite[Lemma 4.2]{Ali+Fedya}.} Assume
(d). We wish to verify (c) using induction on the depth of $\Delta _{0}^{%
\mathrm{fin}}(\mathsf{P})$-formulae. The atomic case is guaranteed by (d),
and the Boolean cases go through since, provably in $\mathsf{KP}$, the
universe is closed under relative complements and intersections. The
existential case, in turn, goes through since provably in $\mathsf{KP}$, the
class of finite sets is closed under Cartesian products, and under finite
unions. Below, we provide further detail for the existential case of the
induction argument. \medskip

\noindent Suppose $\left( \mathcal{M},P\right) \models \forall ^{\mathrm{fin}%
}s\mathrm{\,}\left( s\cap \mathsf{P}\in \mathrm{V}\right) $, and we wish to
show that $\left\{ x\in a:\theta (x)\right\} $ exists in $\mathcal{M}$,
where $a$ is a finite set in the sense of $\mathcal{M}$. We need to consider
two cases:\medskip

\noindent Case 1: $\theta (x)$ is of the form $\exists y\in b\,\left( 
\mathrm{Fin}(b)\wedge \varphi (x,y)\right) $, for some parameter $b\in M.$%
\medskip

\noindent Case 2: $\theta (x)$ is of the form $\exists y\in x\,\left( 
\mathrm{Fin}(x)\wedge \varphi (x,y)\right) .$\medskip

\noindent If Case 1 holds, then if $\mathrm{Fin}(b)$ is false in $\mathcal{M}
$, then $\theta (x)$ does not hold for any $x$ and therefore $\left\{ x\in
a:\theta (x)\right\} $ is the empty set in the sense of $\mathcal{M}$. On
the other hand, if $b$ is finite in the sense of $\mathcal{M}$, then so is $%
a\times b$, and therefore by the inductive assumption there is for some $%
s\in M$ such that:%
\begin{equation*}
\left( \mathcal{M},P\right) \models \left[ s=\left\{ \left\langle
x,y\right\rangle \in a\times b:\varphi (x,y)\right\} \right] .
\end{equation*}

\noindent This makes it clear that if we define $t$ in $\mathcal{M}$ as $%
\left\{ x\in a:\exists y\in b\,\left\langle x,y\right\rangle \in s\right\} $%
, then we have:%
\begin{equation*}
\left( \mathcal{M},P\right) \models \left[ t=\left\{ x\in a:\exists y\in
b\,\left( \mathrm{Fin}(b)\wedge \varphi (x,y)\right) \right\} \right] .
\end{equation*}

\noindent If Case 2 holds, we can find $c$ and $d$ in $M$ such that $%
\mathcal{M}$ satisfies $c=\{x\in a:\mathrm{Fin}(x)\},$ and $d=\cup c.$ Since
\textquotedblleft a finite union of finite sets is finite\textquotedblright\
holds in $\mathcal{M}$, $d$ is finite in the sense of $\mathcal{M}$.
Therefore, by the inductive assumption there is for some $s\in M$ we have:%
\begin{equation*}
\left( \mathcal{M},P\right) \models \left[ s=\left\{ \left\langle
x,y\right\rangle \in a\times d:\varphi (x,y)\right\} \right] .
\end{equation*}

\noindent Hence if we define $t$ in $\mathcal{M}$ as $\left\{ x\in a:\exists
y\in d\,\left\langle x,y\right\rangle \in s\right\} $, then we have:%
\begin{equation*}
\left( \mathcal{M},P\right) \models \left[ t=\left\{ x\in a:\exists y\in
x\,\left( \mathrm{Fin}(x)\wedge \varphi (x,y)\right) \right\} \right] .
\end{equation*}

\hfill $\square $\medskip

\noindent \textbf{6.6.~Definition.}~\medskip

\begin{enumerate}
\item[\textbf{(a)}] $\mathsf{GRef}_{\varnothing }$ (Global Reflection over
First Order Logic) is the following $\mathcal{L}_{\mathrm{set}}(\mathsf{T})$
that expresses \textquotedblleft $\mathsf{T}$ contains all sentences of $%
\mathrm{Sent}^{+}$ that are instances of theorems of first order
logic\textquotedblright .%
\begin{equation*}
\forall \varphi \in \mathrm{Sent}^{+}\,\left( \mathrm{Prov}_{\mathsf{%
\varnothing }}(\varphi )\rightarrow \mathsf{T}(\varphi )\right) .
\end{equation*}

\item[\textbf{(b)}] $\mathsf{GRef}_{\mathsf{T}}^{\text{\textrm{Prop}}}$
(Global Propositional Reflection over $\mathsf{T}$), is the following $%
\mathcal{L}_{\mathrm{set}}(\mathsf{T})$ that expresses \textquotedblleft $%
\mathsf{T}$ is closed under proofs in propositional logic\textquotedblright .%
\begin{equation*}
\forall \varphi \in \mathrm{Sent}^{+}\,\left( \mathrm{Prov}_{\mathsf{T}}^{%
\text{\textrm{Prop}}}(\varphi )\rightarrow \mathsf{T}(\varphi )\right) .
\end{equation*}

\item[\textbf{(c)}] $\mathsf{PI}$ (\textit{Propositional Induction}\footnote{%
This principle is dubbed \textit{Sequential Induction} in \cite{DC paper of
trio}, and is denoted $\mathsf{SI}$.})\textit{\ }is the $\mathcal{L}_{%
\mathrm{set}}(\mathsf{T})$-sentence that asserts that\textit{\ }for all
finite sequences $\left\langle \varphi _{i}:i\leq k\right\rangle ,$ where
each $\varphi _{i}\in \mathrm{Sent}^{+}$, the following holds:%
\begin{equation*}
\left[ \mathsf{T}(\varphi _{0})\wedge \forall i<k\ \left[ \mathsf{T}\left( 
\mathsf{\varphi }_{i}\right) \rightarrow \mathsf{T}\left( \mathsf{\varphi }%
_{i+1}\right) \right] \right] \rightarrow \mathsf{T}(\varphi _{k}).
\end{equation*}

\item[\textbf{(d)}] $\mathsf{SPI}$ (\textit{Strong Propositional Induction}%
\footnote{%
This principle is dubbed \textit{Sequential Order Induction} in \cite{DC
paper of trio}, and is denoted $\mathsf{SOI}$.}) is the $\mathcal{L}_{%
\mathrm{set}}(\mathsf{T})$-sentence that asserts that\textit{\ }for all
finite sequences $\left\langle \varphi _{i}:i\leq k\right\rangle ,$ where
each $\varphi _{i}\in \mathrm{Sent}^{+}$, the following holds:%
\begin{equation*}
\left[ \mathsf{T}(\varphi _{0})\wedge \forall j\leq k\ \left[ \forall i<j\ 
\mathsf{T}\left( \mathsf{\varphi }_{i}\right) \right] \rightarrow \mathsf{T}%
\left( \mathsf{\varphi }_{j}\right) \right] \rightarrow \mathsf{T}(\varphi
_{k}).
\end{equation*}

\item[\textbf{(e)}] Given a finite sequence $\left\langle \varphi
_{i}:i<k\right\rangle $ of sentences, and $1\leq j<k$, $\bigvee\limits_{i<j}%
\varphi _{i}$ is defined by induction on $j$ via:%
\begin{equation*}
\bigvee\limits_{i<1}\varphi _{i}:=\varphi _{0};\ \ \
\bigvee\limits_{i<j+1}\varphi _{i}:=\left( \bigvee\limits_{i<j}\varphi
_{i}\right) \vee \varphi _{j}.
\end{equation*}%
\smallskip

\item[\textbf{(f)}] $\mathsf{DC}$ (\textit{Disjunctive Correctness}) is the
conjunction of $\mathsf{DC}_{\mathrm{out}}$ and $\mathsf{DC}_{\mathrm{in}}$,
where $\mathsf{DC}_{\mathrm{out}}$\textit{\ }is the sentence asserting that
for all finite sequences $\left\langle \varphi _{i}:i<k\right\rangle ,$
where each $\varphi _{i}\in \mathrm{Sent}^{+}$, the following holds:$%
\smallskip $%
\begin{equation*}
\mathsf{T}(\bigvee\limits_{i<k}\varphi _{i})\rightarrow \left( \exists i<k\,%
\mathsf{T}(\varphi _{i})\right) ,\smallskip
\end{equation*}%
and $\mathsf{DC}_{\mathrm{in}}$\textit{\ }is the sentence asserting that for
all finite sequences $\left\langle \varphi _{i}:i<k\right\rangle $, where
each $\varphi _{i}\in \mathrm{Sent}^{+}$, the following holds:%
\begin{equation*}
\left( \exists i<k\ \mathsf{T}(\varphi _{i})\right) \rightarrow \mathsf{T}%
(\bigvee\limits_{i<k}\varphi _{i}).
\end{equation*}%
\medskip
\end{enumerate}

\noindent \textbf{6.7.~Lemma.}\footnote{%
This is the analogue of \cite[Proposition 7]{DC paper of trio}, but the
proof presented here for $(b)\Rightarrow (a)$ uses a different strategy.}%
\textbf{~}\textit{The following are equivalent in }$\mathsf{CT}^{-}[\mathsf{%
KP}]$:\medskip

\begin{enumerate}
\item[$(a)$] $\Delta _{0}^{\mathrm{fin}}$-$\mathsf{Ind(T)}$.\medskip

\item[$(b)$] $\mathsf{SPI.}$\medskip
\end{enumerate}

\noindent \textbf{Proof.}~We first verify $(a)$ $\Rightarrow (b)$. Suppose
to the contrary that $(\mathcal{M},T)\models \Delta _{0}^{\mathrm{fin}}$-$%
\mathsf{Ind(T)}$, but $\mathsf{SPI}$ fails in $(\mathcal{M},T).$ Thus%
\begin{equation*}
(\mathcal{M},T)\models \left[ \exists x\in \omega \,\exists s\,\lnot \mathsf{%
SPI(}s\mathsf{,}x\mathsf{)}\right] \mathsf{,}
\end{equation*}%
where $\mathsf{SPI(}x,s\mathsf{)}$ asserts that $s$ is a finite sequence of
sentences $\left\langle \varphi _{i}:i\leq y\right\rangle $, where $x\leq y$%
, and $\mathsf{SPI}$ holds for $\left\langle \varphi _{i}:i\leq
x\right\rangle .$ It is easy to see that $\mathsf{SPI(}x,s\mathsf{)}$ is $%
\Delta _{0}^{\mathrm{fin}}$ (with parameter $s$), which makes it clear that $%
\lnot \mathsf{SPI(}x,s\mathsf{)}$ is also $\Delta _{0}^{\mathrm{fin}}$.
Since by Lemma 6.5 $(\mathcal{M},T)$ satisfies $\Delta _{0}^{\mathrm{fin}}$-$%
\mathsf{Min(T),}$ there is some $k\in \omega ^{\mathcal{M}}$ such that%
\begin{equation*}
(\mathcal{M},T)\models \left[ \lnot \mathsf{SPI(}k,s\mathsf{)}\wedge \,%
\mathsf{SPI(}k-1,s\mathsf{)}\right] .
\end{equation*}%
It is now easy to derive a contradiction from above using the compositional
properties of the truth predicate.\medskip 

\noindent The proof of $(b)$ $\Rightarrow (a)$ takes more effort. We begin
with the observation that $\mathsf{CT}^{-}[\mathsf{KP}]+\mathsf{SPI}$ can
readily prove $\mathsf{GRef}_{\mathsf{T}}^{\text{\textrm{Prop}}}$ (for more
detail, see the proof of \cite[Proposition 7]{DC paper of trio}). Next we
observe that both $\mathsf{DC}$\ and $\mathsf{Int}$-$\mathsf{Ind}$ (internal
induction; see Definition 4.1) are provable in $\mathsf{CT}^{-}[\mathsf{KP}]+%
\mathsf{GRef}_{\mathsf{T}}^{\text{\textrm{Prop}}}$. The proof of $\mathsf{DC}
$\ in $\mathsf{CT}^{-}[\mathsf{KP}]+\mathsf{GRef}_{\mathsf{T}}^{\text{%
\textrm{Prop}}}$ is elementary. The proof of $\mathsf{Int}$-$\mathsf{Ind}$
in $\mathsf{CT}^{-}[\mathsf{KP}]+\mathsf{GRef}_{\mathsf{T}}^{\text{\textrm{%
Prop}}}$ is based on the observation that $\varphi (\dot{k})$ is derivable
in propositional logic from the assumptions 
\begin{equation*}
\left\{ \varphi (\dot{0})\}\cup \{\varphi (\dot{i})\rightarrow \varphi (\dot{%
j}):i<k-1,\ j=i+1\right\} .
\end{equation*}%
We will show, using $\mathsf{DC}$\ and $\mathsf{Int}$-$\mathsf{Ind}$, that $%
\forall ^{\mathrm{fin}}s\ \left( s\cap \mathsf{T}\in \mathrm{V}\right) $
holds. By Lemma 6.5, this will establish $\Delta _{0}^{\mathrm{fin}}$-$%
\mathsf{Ind(T).}$ For this purpose, we resort to the strategy employed in
the proof of \cite[Lemma 4.4]{Ali+Fedya}. Reasoning in $\mathsf{CT}^{-}[%
\mathsf{KP}]+\mathsf{DC}$\ and $\mathsf{Int}$-$\mathsf{Ind}$, suppose $s$ is
a finite set, and fix a bijection $f:k\rightarrow s$, where $k\in \omega .$
Define a sequence $\left\langle \varphi _{i}:i<k\right\rangle $ of sentences
by: $\varphi _{i}=f(i)$ if $f(i)\in \mathrm{Sent}^{+}$, and otherwise let $%
\varphi _{i}$ is the sentence expressing $0\neq 0.$ Note that there is a set 
$a$ such that $a=s\cap \mathsf{T}\ $iff there is a set $b$ such that $%
b=\left\{ i<k:\mathsf{T}(\varphi _{i})\right\} .$ We will now show that such
a set $b$ exists. Towards this goal, consider the unary formula $\theta (x)$
given by: 
\begin{equation*}
\theta (x):=\bigvee_{i<k}\left( (x=\dot{i})\wedge \varphi _{i}\right) .
\end{equation*}%
\textbf{Claim }$\mathbf{(\ast )}$ $\forall i\in k\,\left[ \mathsf{T}(\varphi
_{i})\leftrightarrow \mathsf{T}(\theta (\dot{i}))\right] .$\medskip

\noindent $(\rightarrow )$ Suppose $\mathsf{T}(\varphi _{i})$ for some $i\in
k.$ Then $\mathsf{T}\left( (\dot{i}=\dot{i})\wedge \varphi _{i}\right) ,$
and hence by $\mathsf{DC}$ we have $\mathsf{T}(\theta (\dot{i})).\medskip $

\noindent $(\leftarrow )$ Suppose $\mathsf{T}(\theta (\dot{j}))$ for some $%
i<\in k.$ Then by $\mathsf{DC}$, there is some $j\in k$ such that $\mathsf{T}%
\left( (\dot{i}=\dot{j})\wedge \varphi _{j}\right) .$ So $\mathsf{T}(\varphi
_{i})$ holds since $\mathsf{T}$ commutes with conjunction and $\mathsf{T}(%
\dot{i}=\dot{j})$ holds iff $i=j.$\medskip

\noindent Using $\mathsf{Int}$-$\mathsf{Ind,}$ we can readily show that for
each $k\in \omega $, there is some $b$ such that $b=\left\{ i\in k:\mathsf{T}%
(\theta (\dot{i}))\right\} $. The inductive step goes through since $\mathsf{%
KP}$ can prove that given $a$ and $b$, $a\cup \{b\}$ exists. Thus, by Claim $%
(\ast )$, there is some $b$ such that $b=\left\{ i\in k:\mathsf{T}\left(
\theta (\dot{i})\right) \right\} $, as desired.

\hfill $\square $\bigskip

\noindent \textbf{6.8.~Lemma.} \textit{The following are provable in} $%
\mathsf{CT}^{-}[\mathsf{KP}]+\mathsf{DC}_{\mathrm{out}}$.\medskip

\begin{enumerate}
\item[$(a)$] $\mathsf{PI.}$\medskip

\item[$(b)$] $\mathsf{DC}_{\mathrm{in}}.$\medskip

\item[$(c)$] $\mathsf{SPI.}$\medskip
\end{enumerate}

\noindent \textbf{Proof.}~The proof of $(a)$ has the following two steps.%
\footnote{%
The proof of (a) is the same as the remarkable proof of \cite[Theorem 8]%
{Cieslinski-book}, it is therefore presented in abridged form.}\medskip

\noindent Step 1. In the first step, given $\left\langle \varphi _{i}:i\leq
k\right\rangle $ such that: 
\begin{equation*}
\left[ \mathsf{T}(\varphi _{0})\wedge \forall i<k\ \left( \mathsf{T}\left( 
\mathsf{\varphi }_{i}\right) \rightarrow \mathsf{T}\left( \mathsf{\varphi }%
_{i+1}\right) \right) \right] \rightarrow \mathsf{T}(\varphi _{k}),
\end{equation*}%
we construct a new sequence of sentences $\left\langle \psi _{i}:i\leq
k\right\rangle $ given by:\medskip 

\begin{center}
$\psi _{0}:=\lnot \varphi _{0},$ and $\psi _{i}:=\left( \lnot \varphi
_{i}\rightarrow \bigvee\limits_{j<i}\lnot \varphi _{j}\right) $ for $1\leq
i\leq k$. \smallskip
\end{center}

\noindent We show in $\mathsf{CT}^{-}[\mathsf{KP}]$ alone that $\mathsf{T}%
(\psi _{i})$ for all $i\leq k.$\medskip

\noindent Step 2. We use $\mathsf{DC}_{\mathsf{out}}$ to show that $\mathsf{T%
}(\varphi _{i})$ for all $i\leq k.$ This concludes the proof of\ $(a)$%
.\medskip

\noindent To prove $(b)$, by $(a)$ it suffices to show that $\mathsf{CT}^{-}[%
\mathsf{KP}]\,+$ $\mathsf{SPI\vdash \mathsf{DC}_{\mathsf{in}},}$ which is
straightforward.\medskip

\noindent (c) By part $(b)$, it suffices to show that $\mathsf{CT}^{-}[%
\mathsf{KP}]\,+$ $\mathsf{DC\vdash \mathsf{SPI.}}$

\hfill $\square $\medskip

\noindent \textbf{6.9.~Corollary.~}$\mathsf{CT}^{-}[\mathsf{ZF}]+\mathsf{DC}%
_{\mathsf{out}}\vdash \Delta _{0}^{\mathrm{fin}}$-$\mathsf{Ind(T).}$\medskip

\noindent \textbf{Proof.}~This follows immediately from Lemmas 6.7 and 6.8.

\hfill $\square $\medskip

\noindent \textbf{6.10.~Theorem.~}\textit{The following are equivalent over} 
$\mathsf{CT}^{-}[\mathsf{KP}]$:\medskip

\begin{enumerate}
\item[$(a)$] $\forall ^{\mathrm{fin}}s\mathrm{\,}\left( s\cap \mathsf{T}\in 
\mathrm{V}\right) .$\medskip

\item[$(b)$] $\Delta _{0}^{\mathrm{fin}}$-$\mathsf{Ind(T).}$ \medskip

\item[$(c)$] $\mathsf{GRef}_{\mathsf{T}}.$\medskip

\item[$(d)$] $\mathsf{GRef}_{\varnothing }.$\medskip

\item[$(e)$] $\mathsf{DC}$\textit{.}\medskip

\item[$(f)$] $\mathsf{GRef}_{\mathsf{T}}^{\text{\textrm{Prop}}}\mathsf{.}$%
\medskip

\item[$(g)$] $\mathsf{DC}_{\mathrm{out}}.$\medskip
\end{enumerate}

\noindent \textbf{Proof.}~$(a)\Rightarrow (b)$: This follows from
Proposition 6.5. \medskip

\noindent $(b)\Rightarrow (c):$ This was demonstrated in the arithmetical
context by \L e\l yk \cite{Mateusz-prolongable} using substantial
bootstrapping. The proof turns out to be somewhat simpler in our context,
since (1) we are working with a truth predicate rather than a satisfaction
predicate, (2) we are already in a set-theoretical context and do not need
to resort to coding in order to simulate set-theoretical notions, and (3)
the language of set theory has no terms other than variables. However, as in
the arithmetical case, the main stumbling block in the `obvious' proof of
closure of $\mathsf{T}$ under first order proofs (using induction on the
lengths of proofs) is that it is not overtly clear that the commutation of $%
\mathsf{T}$ with arbitrary blocks of existential quantifiers is provable
within $\mathsf{CT}^{-}[\mathsf{KP}]+\Delta _{0}^{\mathrm{fin}}$-$\mathsf{%
Ind(T)}$.\footnote{%
Note that since the universal quantifier $\forall x\,\varphi $ in our
context is defined as $\lnot \exists x\lnot \varphi $, using $\Delta _{0}^{%
\mathrm{fin}}$-$\mathsf{Ind(T)}$ we can easily show that $\mathsf{EC}$
implies the dual property $\mathsf{UC}$ (Universal Correctness), which
stipulates that $\mathsf{T}$ commutes with blocks of universal quantifiers.}
More explicitly, within $\mathsf{KP}$, consider the map $\mathrm{ecl}:%
\mathrm{Form}\rightarrow \mathrm{Sent}$, where $\mathrm{ecl}$ stands for
`existential closure', defined as follows: given $\varphi \in \mathrm{Form}%
_{k},$ 
\begin{equation*}
\mathrm{ecl}(\varphi (v_{0},\cdot \cdot \cdot ,v_{k-1})):=\exists v_{0}\cdot
\cdot \cdot \exists v_{k-1}\ \varphi (v_{0},\cdot \cdot \cdot ,v_{k-1}).
\end{equation*}%
Recall (from the paragraph preceding Proposition 3.2.6) that within $\mathsf{%
KP}$ we can define the map $\left\langle \varphi ,\alpha \right\rangle
\mapsto \varphi \ast \alpha $, where $\varphi \in \mathrm{Form}$, and $%
\alpha $ is an assignment for $\varphi $ (i.e., $\alpha $ is a function
whose domain is the set \textrm{FV}($\varphi )$ of free variables of $%
\varphi )$, where $\varphi \ast \alpha $ is the sentence obtained by
replacing each occurrence of a free variable $x$ of $\varphi $ with $\dot{m}$%
, where $\alpha (x)=m.$ The following statement $\mathsf{EC}$ (Existential
Correctness\footnote{%
This terminology was suggested in \cite{Ali+Albert-long}, where it was shown
that it $\mathsf{EC}$ can be conservatively added to $\mathsf{CT}^{-}$ over
an arbitrary base theory that is formulated in a relational language. The
same concept was dubbed `quantifier correctness' by Wcis\l o \cite%
{Wcislo-NDJFL}, who showed that it can be conservatively added to $\mathsf{PA%
}$ (in its usual formulation using function symbols).}) expresses the
commutation of $\mathsf{T}$ with arbitrary blocks of existential
quantifiers. In the statement below $\mathrm{Asn}(\alpha ,x)$ expresses
\textquotedblleft $\alpha $ is an assignment for $\varphi $%
\textquotedblright , as in Definition 3.2.1.\medskip 

\noindent $(\mathsf{EC})\ \ \ \forall \varphi \in \mathrm{Form}\,\left[ 
\mathsf{T}(\mathrm{ecl}(\varphi ,x))\right] \leftrightarrow \left[ \exists
\alpha \,\mathrm{Asn}(\alpha ,x)\wedge \mathsf{T}(\varphi \ast \alpha )%
\right] .$\medskip\ 

\noindent We will employ \L e\l yk's strategy \cite{Mateusz-prolongable} for
showing that $\mathsf{EC}$ is provable in $\mathsf{CT}^{-}[\mathsf{KP}%
]+\Delta _{0}^{\mathrm{fin}}$-$\mathsf{Ind(T).}$ For this purpose, consider
the map $\mathrm{becl}:\mathrm{Form\times Var}\rightarrow \mathrm{Sent},$
where $\mathrm{becl}$ stands for `bounded existential closure', and is
defined by: 
\begin{equation*}
\mathrm{becl(}\varphi (v_{0},\cdot \cdot \cdot ,v_{k-1}),x)=:\exists
v_{0}\in x\cdot \cdot \cdot \exists v_{k-1}\in x\ \varphi (v_{0},\cdot \cdot
\cdot ,v_{k-1}).
\end{equation*}%
The following principle $\mathsf{BEC}$ (Bounded Existential Correctness) is
readily provable in $\mathsf{CT}^{-}[\mathsf{KP}]+\Delta _{0}^{\mathrm{fin}}$%
-$\mathsf{Ind(T)}$:\medskip 

\noindent $(\mathsf{BEC})\ \ \ \forall \varphi \in \mathrm{Form}\,\exists
x\,\left( \mathsf{T}(\mathrm{becl}(\varphi ,x))\leftrightarrow \left[
\exists \alpha \left[ \mathrm{Asn}(\alpha ,x)\,\wedge \mathsf{T}(\varphi
\ast \alpha )\wedge \forall i\in k\,\alpha (i)\in x\,\right] \right] \right)
.$\medskip 

\noindent With $\mathsf{BEC}$ at hand, the $\Sigma _{1}^{\mathrm{fin}}$%
-statement $\left[ \exists \alpha \,\left( \mathrm{Asn}(\alpha ,\varphi
)\,\wedge \mathsf{T}(\varphi \ast s)\right) \right] $ is equivalent to the
following $\Delta _{0}^{\mathrm{fin}}$-statement 
\begin{equation*}
\mathsf{T}(\exists x\,\mathrm{becl}(\varphi ,x)),
\end{equation*}%
where $x$ is not a free variable of $\varphi .$ This shows that $\mathsf{EC}$
is provable in $\mathsf{CT}^{-}[\mathsf{KP}]+\Delta _{0}^{\mathrm{fin}}$-$%
\mathsf{Ind(T)}$, which in turn allows one to prove that $\mathsf{T}$ is
closed under first order proofs, similar to the proof in \cite%
{Mateusz-prolongable} for the arithmetical case. \medskip

\noindent $(c)\Rightarrow (d):$ Trivial.\medskip

\noindent $(d)\Rightarrow (e):$ We will use a variant of the construction
used earlier in the proof of Lemma 6.7. Let $\mathrm{Sqsent}(s,k)$ be $%
\mathrm{Seq}(s,k)\wedge \forall i\in k\,\left( s(i)\in \mathrm{Sent}%
^{+}\right) ,$ where $\mathrm{Seq}(s,k)$ expresses \textquotedblleft $s$ is
a sequence of length $k$\textquotedblright . Within $\mathsf{KP}$, we can
define a map that associates a unary formula $\theta _{s}(x)$ to any $%
s=\left\langle s_{i}:i<k\right\rangle $, where $\mathrm{Sqsent}(s,k)$ holds,
via: 
\begin{equation*}
\theta _{s}(x)=\bigvee\limits_{i<k}\left[ \left( x=\dot{i}\right) \wedge
s_{i}\right] .
\end{equation*}

\noindent Let $\alpha _{\mathsf{KP}}$ be the conjunction of the finitely
many axioms of $\mathsf{KP}$. The following can be readily verified:\medskip

\begin{enumerate}
\item[$(i)$] $\mathsf{KP}\vdash \left[ \forall k\in \omega \ \forall
s\,\forall i<k\,\left[ \mathrm{Sqsent}(s,k)\rightarrow \mathrm{Prov}%
_{\varnothing }\left( \alpha _{\mathsf{KP}}\rightarrow \left( \theta _{s}(%
\dot{i})\leftrightarrow s_{i}\right) \right) \right] \right] .\medskip $

\item[$(ii)$] $\mathsf{KP}\vdash \left[ \forall k\in \omega \ \forall
s\,\forall i<k\,\left[ \mathrm{Sqsent}(s,k)\rightarrow \mathrm{Prov}%
_{\varnothing }\left( \alpha _{\mathsf{KP}}\rightarrow \left( \left( \exists
x\in i\,\theta _{s}(x)\right) \leftrightarrow
\bigvee\limits_{i<k}s_{i}\right) \right) \right] \right] $.\medskip 
\end{enumerate}

\noindent Thanks to the above, the following are provable in $\mathsf{CT}%
^{-}[\mathsf{KP}]+\mathsf{GRef}_{\mathsf{\varnothing }}$ :$\medskip $

\begin{enumerate}
\item[$(\ast )$] $\forall k\in \omega \ \forall s\,\forall i<k\,\left[ 
\mathrm{Sqsent}(s,k)\rightarrow \left( \mathsf{T}\left( \theta _{s}(\dot{i}%
)\leftrightarrow \mathsf{T}(s_{i})\right) \right) \right] .\medskip $

\item[$(\ast \ast )$] $\forall k\in \omega \ \forall s\,\forall i<k\,\left[ 
\mathrm{Sqsent}(s,k)\rightarrow \left( \mathsf{T}\left( \exists x\in
i\,\theta _{s}(x)\right) \leftrightarrow \mathsf{T}\left(
\bigvee\limits_{i<k}s_{i}\right) \right) \right] .\medskip $
\end{enumerate}

\noindent Note that $\mathsf{DC}$ follows immediately from $(\ast )$\ and $%
(\ast \ast )$.\medskip

\noindent $(e)\Rightarrow (f):$ This follows from Corollary 6.9, since $%
\mathsf{GRef}_{\mathsf{T}}^{\text{\textrm{Prop}}}$ is easily provable in $%
\Delta _{0}^{\mathrm{fin}}$-$\mathsf{Ind(T)}$. \medskip

\noindent $(f)\Rightarrow (g):$ This can be verified by an elementary 
\textit{reductio ad absurdum} argument. \medskip

\noindent $(g)\Rightarrow (a):$ This follows from Lemma 6.5 and Corollary
6.9.

\hfill $\square $\medskip

\noindent \textbf{6.11.~Remark.~}In the proof of $(b)\Rightarrow (c)$ of
Theorem 6.10, we noted that $\mathsf{EC}$ is provable in $\mathsf{CT}^{-}[%
\mathsf{KP}]+\Delta _{0}^{\mathrm{fin}}$-$\mathsf{Ind(T)}$. Together with
the provability of $\mathsf{DC}$ in $\mathsf{CT}^{-}[\mathsf{KP}]+\Delta
_{0}^{\mathrm{fin}}$-$\mathsf{Ind(T)}$, this shows that within $\mathsf{CT}%
^{-}[\mathsf{KP}]+\Delta _{0}^{\mathrm{fin}}$-$\mathsf{Ind(T),}$ $\mathsf{T}$
is uniquely determined by the proper subset $\mathsf{T}_{1}$ of itself that
consists of sentences in $\mathrm{Sent}^{+}$ that contain precisely one
constant symbol. This is because within $\mathsf{CT}^{-}[\mathsf{KP}]+%
\mathsf{EC}+\mathsf{DC}$, given $\varphi (x_{0},\cdot \cdot \cdot
,x_{k-1})\in \mathrm{Form}_{k}$ and $\alpha \in \mathrm{Asn}(\varphi )$, we
have: 
\begin{equation*}
\mathsf{T}(\varphi \ast \alpha )\longleftrightarrow \mathsf{T}(\widehat{%
\varphi }\ast \widehat{\alpha }),
\end{equation*}%
where $\widehat{\varphi }(y):=\left[ \text{\textrm{\textquotedblleft }}%
y:k\rightarrow \mathrm{V}\text{\textrm{\textquotedblright }}\wedge \exists
x_{0}\cdot \cdot \cdot \exists x_{k-1}\left( \varphi (x_{0},\cdot \cdot
\cdot ,x_{k-1})\wedge \bigwedge\limits_{i<k}y(i)=x_{i}\right) \right] $, and 
$\widehat{\alpha }(y)=f,$ where $f$ is the function with domain $k$ such
that $f(i)=\alpha (x_{i})$ for each $i\in k.$\medskip

\noindent \textbf{6.12.~Many Faces Theorem for }$\mathsf{CT}_{\ast }[\mathsf{%
ZF}].$~\textit{The following axiomatize the same theory} \textit{over} $%
\mathsf{CT}^{-}[\mathsf{KP}]$:\medskip

\begin{enumerate}
\item[$(a)$] $\mathsf{Int}$-$\mathsf{Repl+GRef}_{\mathsf{T}}$.\medskip

\item[$(b)$] $\mathsf{Int}$-$\mathsf{Repl+GRef}_{\varnothing }$.\medskip

\item[$(c)$] $\mathsf{Int}$-$\mathsf{Repl+GRef}_{\mathsf{T}}^{\text{\textrm{%
Prop}}}.$\medskip

\item[$(d)$] $\mathsf{GRef}_{\mathrm{ZF}}\mathrm{.}$\medskip

\item[$(e)$] $\mathsf{Int}$-$\mathsf{Repl}$\textsf{\ +}$\mathrm{\,}\forall ^{%
\mathrm{fin}}s$ $\mathsf{T}\cap s\in \mathrm{V.}$\medskip

\item[$(f)$] $\mathsf{Int}$-$\mathsf{Repl}$\textsf{\ + }$\mathsf{DC}.$%
\medskip

\item[$(g)$] $\mathsf{Int}$-$\mathsf{Repl}$\textsf{\ + }$\mathsf{DC}_{%
\mathsf{out}}.$\medskip

\item[$(h)$] $\mathsf{Int}$-$\mathsf{Ref}$.\footnote{%
Recall that $\mathsf{Int}$-$\mathsf{Ref}$ (internal reflection) was defined
in Definition 5.3.}\medskip

\item[$(i)$] $\mathsf{Int}$-$\mathsf{Ref}^{\mathrm{weak}}:=\forall \varphi
\in \mathrm{Form\,}\forall \alpha _{0}\in \mathrm{Ord\,}\exists \alpha \in 
\mathrm{Ord\,}\left[ \left( \alpha _{0}<\alpha \right) \wedge \mathsf{T}(%
\mathrm{V\!}_{\alpha }\prec _{\varphi }\mathrm{V})\right] .$\medskip
\end{enumerate}

\noindent \textbf{Proof.}~The equivalence of $(a)$ through $(g)$ follows
from Theorem 6.9. Also $(i)$ easily follows from $(h)$, and in Theorem 5.3
we saw that $\mathsf{Int}$-$\mathsf{Ref}$ is provable in $\mathsf{CT}_{\ast
}[\mathsf{ZF}]$. \medskip

\noindent So the proof of the theorem is complete once we verify that $%
(i)\Rightarrow (b)$. It is not hard to see that internal collection and
internal separation hold within $\mathsf{CT}^{-}[\mathsf{ZF}]+\mathsf{Int}$-$%
\mathsf{Ref}^{\mathrm{weak}}$, which by Remark 4.2, shows that internal
replacement holds in $\mathsf{CT}^{-}[\mathsf{ZF}]+\mathsf{Int}$-$\mathsf{Ref%
}^{\mathrm{weak}}.$ It remains to show that $\mathsf{GRef}_{\varnothing }$
is provable in $\mathsf{CT}^{-}[\mathsf{ZF}]+\mathsf{Int}$-$\mathsf{Ref}^{%
\mathrm{weak}}$\textsf{.} For this purpose, suppose $(\mathcal{M},T)\models 
\mathsf{CT}^{-}[\mathsf{ZF}]+\mathsf{Int}$-$\mathsf{Ref}^{\mathrm{weak}}$,
and $\mathcal{M}\models \mathrm{Prov}_{\varnothing }(\varphi )$ for some $%
\varphi \in \left( \mathrm{Sent}^{+}\right) ^{\mathcal{M}}.$ Let $\alpha
_{0} $ be the first ordinal in $\mathrm{Ord}^{\mathcal{M}}$ such that $%
\mathcal{M} $ satisfies \textquotedblleft $\mathrm{V\!}_{\alpha }$ contains
all the constants in $\varphi $\textquotedblright . By the soundness theorem
of first order logic within $\mathcal{M}$, for all $\alpha $ in $\mathrm{Ord}%
^{\mathcal{M}}$ greater than $\alpha _{0}$ we have: 
\begin{equation*}
\mathcal{M}\models \left[ \varphi \in \mathrm{ED}(\mathrm{V\!}_{\alpha },\in
)\right] .
\end{equation*}%
By the veracity of $\mathsf{Int}$-$\mathsf{Ref}^{\mathrm{weak}}$ in $(%
\mathcal{M},T),$ there is some $\alpha \in \mathrm{Ord}^{\mathcal{M}}$ that
exceeds $\alpha _{0}$ such that $\mathsf{T}(\mathrm{V\!}_{\alpha }\prec
_{\varphi }\mathrm{V}).$ This makes it clear that $(\mathcal{M},T)\models 
\mathsf{T}(\varphi )$.

\hfill $\square $\medskip

\noindent \textbf{6.13.~Theorem.} $\mathsf{CT}^{-}[\mathsf{ZF}]+\mathsf{Int}$%
-$\mathsf{Repl}+\mathsf{DC}_{\mathrm{in}}$ \textit{is conservative over} $%
\mathsf{ZF}$.\medskip

\noindent \textbf{Proof.}~The proof strategy of \cite[Theorem 20]{DC paper
of trio} is readily adaptable to our context.\footnote{%
As in \cite[Remark 27]{DC paper of trio}, further `good properties' can be
added to this conservativity result; e.g., existential correctness (defined
in the proof of Theorem 6.9), and the agreement of $\mathsf{T}$ with the all
internal truth predicates $\left\{ \mathrm{True}_{n}:n\in \mathbb{N}\right\}
.$}

\hfill $\square $\bigskip \bigskip

\begin{center}
\textbf{7.~}$\mathsf{CT}_{\ast }[\mathsf{ZF}]$ \textbf{AND} $\mathsf{GB}%
^{\ast }$\bigskip
\end{center}

\noindent This section uses some results established in \cite[Lemma 4.3]%
{Ali-Mostowski-Bridge}, which are reviewed below. We use $\mathsf{GB}$ to
denote the G\"{o}del-Bernays theory of classes. We view $\mathsf{GB}$ as a
two-sorted theory, whose models can be represented in the form $(\mathcal{M},%
\mathfrak{X})$, where $\mathcal{M}$ is an $\mathcal{L}_{\mathsf{set}}$%
-structure, $\mathfrak{X}\subseteq \mathcal{P}(M)\mathfrak{.}$ It is
well-known that $(\mathcal{M},\mathfrak{X})\models \mathsf{GB}$ iff the
following two conditions hold:\medskip

\begin{enumerate}
\item[\textbf{(a)}] If $X_{1},\cdot \cdot \cdot ,X_{n}\in \mathfrak{X}$,
then $(\mathcal{M},X_{1},\cdot \cdot \cdot ,X_{n})\models \mathsf{ZF}(%
\mathsf{X}_{1},\cdot \cdot \cdot ,\mathsf{X}_{n})$.\medskip

\item[\textbf{(b)}] $X_{1},\cdot \cdot \cdot ,X_{n}\in \mathfrak{X},$ and $Y$
is parametrically definable in $(\mathcal{M},X_{1},\cdot \cdot \cdot ,X_{n})$%
, then $Y\in \mathfrak{X}$.\medskip
\end{enumerate}

\noindent The following result is well-known.\medskip

\noindent \textbf{7.1}.~\textbf{Theorem.}~$\mathsf{GB}$ \textit{is
conservative over} $\mathsf{ZF}$; \textit{but }$\mathsf{GB}$ \textit{is not
interpretable in }$\mathsf{ZF}.$\medskip

\noindent Mostowski \cite{Mostowski-impredicative} showed that there is a
formula $\mathrm{T}_{\mathrm{Most}}(x)$ -- dubbed the \textit{Mostowski
truth predicate} here -- such that for all sentences $\varphi $ in the
language of $\mathcal{L}_{\mathrm{set}}$ of $\mathsf{ZF}$, we have:\medskip

\begin{center}
$\mathsf{GB}\vdash \varphi \leftrightarrow \mathrm{T}_{\mathrm{Most}%
}(\ulcorner \varphi \urcorner ).$\medskip
\end{center}

\noindent The conservativity of $\mathsf{GB}$ over $\mathsf{ZF}$ combined
with G\"{o}del's second incompleteness theorem makes it clear that $\mathrm{%
Con}(\mathsf{ZF})$\ is unprovable in $\mathsf{GB}$. Together with the above
result about $\mathrm{T}_{\mathrm{Most}}$, we witness a striking phenomenon: 
$\mathsf{GB}$ possesses a `truth-predicate' for $\mathsf{ZF}$, and yet the
formal consistency of $\mathsf{ZF}$ is unprovable in $\mathsf{GB}$. \medskip

\noindent \textbf{7.2.~Definition.}~The following definitions should be
understood to be carried out within $\mathsf{GB}$.\medskip

\begin{enumerate}
\item[\textbf{(a)}] $\mathrm{Depth}_{k}$ is the collection of $\mathcal{L}_{%
\mathrm{set}}$-formulae $\varphi $ with $\mathrm{depth}(\varphi )\leq k$;
here $\mathrm{depth}(\varphi )$ is the length of the longest path in the
parsing tree of $\varphi $ (also known as the formation/syntactic tree).
Thus $\mathrm{Depth}_{0}=\varnothing $, and $\mathrm{Depth}_{1}$ consists of
atomic formulae.

\item[\textbf{(b)}] The \textit{Mostowski cut}\footnote{%
In the terminology of models of arithmetic, a \textit{cut} of a nonstandard
model $\mathcal{M}$ of arithmetic is \textit{an initial segment of }$%
\mathcal{M}$ \textit{that is closed under immediate successors}. This
terminology can be readily applied to $\omega $-nonstandard models of set
theory.}, denoted $\mathrm{C}_{\mathrm{Most}}$, consists of $k\in \omega $
such that there is a class $T$ with the property that $T$ is a $\mathrm{Depth%
}_{k}$-truth class for the structure $($\textrm{V}$,\in )$. $\mathrm{Depth}_{%
\mathrm{Most}}$ consists of $\mathcal{L}_{\mathrm{set}}$-formulae $\varphi $
such that $\mathrm{depth}(\varphi )\in \mathrm{C}_{\mathrm{Most}}.$

\item[\textbf{(c)}] The \textit{Mostowski truth} \textit{predicate, }denoted 
$\mathrm{T}_{\mathrm{Most}}(x)$ expresses:\medskip
\end{enumerate}

\begin{center}
$x$ is (the code of) an $\mathcal{L}_{\mathrm{set}}^{+}$-formula $\varphi (%
\dot{m}_{0},...,\dot{m}_{k-1})$ and \medskip

$\exists p\geq \mathrm{depth}(\varphi )$ $\exists T$ $[\varphi (\dot{m}%
_{0},\cdot \cdot \cdot ,\dot{m}_{k-1})\in T$ and $T$ is a $\mathrm{Depth}%
_{p} $-truth class over $($\textrm{V}$,\in )]$. \medskip
\end{center}

\noindent \textbf{7.3.~Lemma.}~\cite[Lemma 3.2]{Ali-Mostowski-Bridge} 
\textit{Provably in} $\mathsf{GB}$, $\mathsf{C}_{\mathsf{Most}}$\ \textit{is
a cut of }$\omega $.\medskip

\noindent \textbf{7.4.~Theorem.}~\textit{Provably in} $\mathsf{GB}$, $%
\mathrm{T}_{\mathrm{Most}}(x)$ \textit{is an }$\mathsf{F}$-\textit{truth
class for }$\mathsf{F}=\mathrm{Depth}_{\mathrm{Most}}$. \textit{In
particular,} \textit{if} $(\mathcal{M},\mathfrak{X})\models \mathsf{GB}$, 
\textit{then for every standard }$\mathcal{L}_{\mathrm{set}}$-\textit{%
formula }$\varphi (x_{1},...,x_{n})$, \textit{and any sequence} $%
\left\langle m_{0},...m_{k-1}\right\rangle $ \textit{of elements of} $%
\mathcal{M}$, \textit{the following equivalence holds}:\medskip 

\begin{center}
$\mathcal{M}\models \varphi (m_{0},\cdot \cdot \cdot ,m_{k-1})$ iff $(%
\mathcal{M},\mathfrak{X})\models \left[ \varphi (\dot{m}_{0},\cdot \cdot
\cdot ,\dot{m}_{k-1})\in \mathsf{T}_{\mathsf{Most}}\right] .$ \medskip
\end{center}

\noindent \textbf{7.5.~Remark.}~As shown in \cite[Theorem 3.7]%
{Ali-Mostowski-Bridge}, \textit{given }$(\mathcal{M},\mathfrak{X})\models 
\mathsf{GB}$, if $\mathfrak{X}=\mathfrak{X}_{\mathrm{Def}\mathsf{(}\mathcal{M%
}\mathsf{)}}$ or $\mathcal{M}$ is not recursively saturated, then $\mathrm{C}%
_{\mathrm{Most}}^{(\mathcal{M},\mathfrak{X)}}=$\textit{\ }$\omega $\textit{.}%
\medskip

\noindent \textbf{7.6.~Definition.}~The following are the set-theoretical
analogue of the extensions $\mathsf{ACA}_{0}^{\ast }$ and $\mathsf{ACA}%
_{0}^{\prime }$ of $\mathsf{ACA}_{0}$ (in the notation of \cite[Theorem 3.7]%
{Ali-Mostowski-Bridge}).\medskip

\noindent \textbf{(a)} $\mathsf{GB}^{\ast }=\mathsf{GB\,}+\,\forall k\in
\omega \ \exists T\ \mathrm{Tr}(T,k)$, where $\mathrm{Tr}(T,k)$ expresses
\textquotedblleft $T$ is a $\mathrm{Depth}_{k}$-truth class over $($\textrm{V%
}$,\in )$\textquotedblright .\footnote{%
Thus $\mathsf{GB}^{\ast }=\mathsf{GB}+\forall x(x\in \omega \rightarrow x\in 
\mathsf{C}_{\mathsf{Most}}).$}\medskip

\noindent \textbf{(b)} $\mathsf{GB}^{\prime }=\mathsf{GB}+\forall X\ \forall
k\in \omega \ \exists T\ \mathrm{Tr}(T,X,k)$, where $\mathrm{Tr}(T,X,k)$
expresses \textquotedblleft $T$ is a $\mathrm{Depth}_{k}$-truth class over $%
( $\textrm{V}$,\in ,X)$\textquotedblright .\medskip

\noindent The following is a special case of \cite[Theorem 3.15]%
{Ali-Mostowski-Bridge}.\medskip

\noindent \textbf{7.7.~Theorem. }$\mathsf{GB}^{\ast }\vdash \forall \varphi
\left( \mathrm{Prov}_{\mathsf{ZF}}(\varphi )\rightarrow \mathrm{T}_{\mathrm{%
Most}}(\varphi )\right) $.\medskip

\noindent \textbf{7.8.~Definition. }Suppose $(\mathcal{M},T\mathfrak{)}%
\models \mathsf{CT}^{-}[\mathsf{ZF}]$. \ $\medskip $

\noindent \textbf{(a)} For each $\varphi (x,v)\in \mathrm{Form}^{\mathcal{M}%
} $, and $p\in M,$ let $\varphi ^{T}(x,\dot{p}):=\left\{ m\in M:(\mathcal{M}%
,T,p)\models \left[ \varphi (\dot{m},\dot{p})\in \mathsf{T}\right] \right\}
.\medskip $

\noindent \textbf{(b)} $\mathfrak{X}_{\mathrm{Def}_{T}\mathrm{(}\mathcal{M}%
\mathrm{)}}$ is the collection of subsets of $M$ that are of the form $%
\varphi ^{T}(x,p)$ for some unary $\varphi (x,v)\in \mathrm{Form}^{\mathcal{M%
}}$ and some parameter $p\in M.\medskip $

\noindent \textbf{7.9.~Theorem.}~\textit{Suppose }$(\mathcal{M},T\mathfrak{)}%
\models \mathsf{CT}^{-}[\mathsf{ZF}]$. \textit{The following are equivalent}:%
$\medskip $

\begin{enumerate}
\item[$(a)$] $(\mathcal{M},T\mathfrak{)}\models \mathsf{Int}$-$\mathsf{Repl}%
.\medskip $

\item[$(b)$] $(\mathcal{M},\mathfrak{X}_{\mathrm{Def}_{T}\mathrm{(}\mathcal{M%
}\mathrm{)}}\mathfrak{)}\models \mathsf{GB}.$\medskip
\end{enumerate}

\noindent \textbf{Proof.}~We only sketch the proof. The direction $%
(b)\Rightarrow (a)$ is straightforward. The other direction follows from the
following two observations:\medskip

\begin{enumerate}
\item[$(1)$] If $(\mathcal{M},T\mathfrak{)}\models \mathsf{CT}^{-}[\mathsf{ZF%
}]$, then $(\mathcal{M},\mathfrak{X}_{\mathrm{Def}_{T}\mathrm{(}\mathcal{M}%
\mathrm{)}}\mathfrak{)}$ satisfies all of the axioms of $\mathsf{GB}$ with
the possible exception of the replacement axiom.\medskip

\item[$(2)$] For each $\varphi (v,x,y)\in \mathrm{Form}^{\mathcal{M}}$, and $%
X=$ $\varphi ^{T}(\dot{p},x,y)$ for some parameter $p\in M$, then $(\mathcal{%
M},X)\models \mathsf{ZF}(\mathsf{X})$ iff $(\mathcal{M},T\mathfrak{)}\models 
\mathsf{T}(\mathsf{Repl}_{\varphi }).$
\end{enumerate}

\hfill $\square $\medskip

\noindent \textbf{7.10.}~\textbf{Definition.}~Within $\mathsf{KP}$, let $%
\left\langle \mathrm{True}_{k}:k\in \omega ^{\mathcal{M}}\right\rangle $ be
the sequence of $\mathcal{L}_{\mathrm{set}}$-formulae that is defined by the
following recursive clauses:

\begin{equation*}
\mathrm{True}_{1}(\varphi ):=\exists y\exists z\,\left[ \left( \left(
\varphi =\ulcorner \dot{y}=\dot{z}\urcorner \right) \wedge (y=z)\right) \vee
\left( \left( \varphi =\ulcorner \dot{y}\in \dot{z}\urcorner \right) \wedge
(y\in z)\right) \right] .
\end{equation*}

\noindent For $k\geq 2,$ 
\begin{equation*}
\mathrm{True}_{k}(\varphi ):=\left[ 
\begin{array}{c}
\left[ \left( \mathrm{depth}(\varphi )=1\right) \wedge \mathrm{True}%
_{1}(\varphi )\right] \vee \smallskip \\ 
\bigvee\limits_{0<r<k}\left[ \left( \mathrm{depth}(\varphi )=r\right) \wedge
\left( \mathrm{Neg}_{r}(\varphi )\vee \mathrm{Disj}_{r}(\varphi )\vee 
\mathrm{Exist}_{r}(\varphi )\right) \right] ,%
\end{array}%
\right]
\end{equation*}%
where:

\begin{center}
$\mathrm{Neg}_{r}(\varphi ):=\left[ \exists \psi \left( \left( \varphi
=\lnot \psi \right) \wedge \lnot \mathrm{True}_{r}(\psi )\right) \right] ,$%
\medskip

$\mathrm{Disj}_{r}(\varphi ):=\left[ \exists \psi _{1}\,\exists \psi
_{2}\left( \left( \varphi =\psi _{1}\vee \psi _{2})\right) \wedge \left( 
\mathrm{True}_{r}(\psi _{1})\vee \mathrm{True}_{r}(\psi _{2})\right) \right) %
\right] $, and\medskip

$\mathrm{Exist}_{r}(\varphi ):=\left[ \exists \psi \ \left( \left( \varphi
=\exists x\,\psi (x)\right) \wedge \exists v\,\mathrm{True}_{r}(\psi (\dot{v}%
)\right) \right] $

$.$\medskip
\end{center}

\noindent \textbf{7.11.}~\textbf{Lemma.}~\textit{Suppose }$k,n\in \mathbb{N}$
\textit{and} $\varphi $ \textit{is an }$n$-\textit{ary} $\mathcal{L}_{%
\mathrm{set}}$ \textit{formula of depth} $k$. \textit{Then} 
\begin{equation*}
\mathsf{KP}\vdash \left[ \forall x_{0}\cdot \cdot \cdot \forall x_{n-1}\,%
\left[ \varphi (x_{0},\cdot \cdot \cdot ,x_{n-1})\leftrightarrow \mathrm{True%
}_{k}\left( \varphi \left( \dot{x}_{0},\cdot \cdot \cdot ,\dot{x}%
_{n-1}\right) \right) \right] \right] .
\end{equation*}%
\noindent \textbf{Proof.}~This can be established by a straightforward
induction in the metatheory.

\hfill $\square $\medskip

\noindent ~\textbf{7.12.~Lemma.}~$\mathsf{KP}\vdash \left[ \forall k,n\in
\omega \ \forall \varphi \in \mathrm{Form}_{n}\cap \mathrm{Depth}_{k}\
\forall \alpha \in \mathrm{Asn}(\varphi )\ \mathrm{Prov}_{\mathrm{KP}}\left[
\left( \varphi \ast \alpha \leftrightarrow \mathrm{True}_{k}\left( \varphi
\ast \alpha \right) \right) \right] \right] $.\medskip

\noindent \textbf{Proof.}~This is based on the observation that Lemma 7.11
states a simple syntactic fact that is readily provable in $\mathsf{KP}$.

\hfill $\square $\medskip

\noindent \textbf{7.13.}~\textbf{Lemma.}~$\mathsf{CT}^{-}[\mathsf{KP}%
]+\Delta _{0}^{\mathrm{fin}}$-$\mathsf{Ind(T)}\vdash \left[ \forall \sigma
\in \mathrm{Depth}_{k}\cap \mathrm{Sent}^{+}\,\left( \mathsf{T}\left( 
\mathrm{True}_{k}(\dot{\sigma})\right) \leftrightarrow \mathsf{T}\left(
\sigma \right) \right) \right] .$ \medskip

\noindent \textbf{Proof.}~This follows from putting Lemma 7.12 together with 
$\mathsf{GRef}_{\mathsf{T}}$; recall that $\mathsf{GRef}_{\mathsf{T}}$ was
shown in Theorem 6.10 to be provable in $\mathsf{CT}^{-}[\mathsf{KP}]+\Delta
_{0}^{\mathrm{fin}}$-$\mathsf{Ind(T)}$.

\hfill $\square $\medskip

\noindent \textbf{7.14.~Corollary.}~\textit{If }$\left( \mathcal{M},T\right)
\models \mathsf{CT}^{-}[\mathsf{KP}]+\Delta _{0}^{\mathrm{fin}}$-$\mathsf{%
Ind(T)}$, \textit{then for each }$k\in \omega ^{\mathcal{M}},$ $T\left( 
\mathrm{True}_{k}(x)\right) $ \textit{is a }$\mathrm{Depth}_{k}$-\textit{%
truth class over} $\mathcal{M}.$\footnote{%
This is the set-theoretic analogue of \cite[Lemma 3.7]{Bartek+Mateusz on CT0}%
. It can also be proved by taking advantage of the veracity of $\mathsf{DC}$
in $\mathsf{CT}^{-}[\mathsf{KP}]+\Delta _{0}^{\mathrm{fin}}$-$\mathsf{Ind(T).%
}$}\medskip

\noindent \textbf{Proof.}~This follows directly from Lemma 7.13.

\hfill $\square $\medskip

\noindent \textbf{7.15.~Theorem.}~$\mathsf{GB}^{\ast }$, $\mathsf{GB}%
^{\prime },$ \textit{and} $\mathsf{CT}_{\ast }[\mathsf{ZF}]$ \textit{are
pairwise mutually} $\mathrm{V}$-\textit{interpretable.} \medskip

\noindent \textbf{Proof.}~This is proved using Corollary 7.14, with the same
strategy as \cite[Theorem 3.15]{Ali-Mostowski-Bridge}.\footnote{%
As noted by \L e\l yk, by using Lemma 7.13 instead of Corollary 7.14 one can
show that the interpretations witnessing the mutual $\mathrm{V}$%
-interpretability of $\mathsf{GB}^{\prime }$ and $\mathsf{CT}_{\ast }[%
\mathsf{ZF}]$ witness that $\mathsf{CT}_{\ast }[\mathsf{ZF}]$ is a
\textquotedblleft $\mathrm{V}$-retract\textquotedblright\ of $\mathsf{GB}%
^{\prime }.$ The arithmetical analogue of this fact follows from the
provability of \textquotedblleft $\mathsf{T}$ is closed under
proofs\textquotedblright\ in $\mathsf{CT}_{0}[\mathsf{PA}]$, established by 
\L e\l yk \cite{Mateusz-prolongable}.} \hfill $\square $\medskip 

\noindent \textbf{7.16.~Corollary.} $\mathsf{GB}^{\ast }$, $\mathsf{GB}%
^{\prime },$ \textit{and} $\mathsf{CT}_{\ast }[\mathsf{ZF}]$ \textit{have
the same }$\mathcal{L}_{set}$-\textit{consequences.}\medskip

\noindent \textbf{7.17.~Remark.}~The results in this section formulated for
the $\mathrm{Depth}_{n}$ hierarchy also hold for the $\Sigma _{n}$%
-hierarchy.\bigskip \bigskip

\begin{center}
\textbf{8.~SET THEORETICAL CONSEQUENCES OF} $\mathsf{CT}_{\ast }[\mathsf{ZFC}%
]$\bigskip
\end{center}

\noindent The main results of Sections 6 and 7 suggest that the
set-theoretical counterpart of the theory $\mathsf{CT}_{0}[\mathsf{PA}]$ is $%
\mathsf{CT}_{\ast }[\mathsf{ZF}]$, and not the much stronger theory $\mathsf{%
CT}_{0}[\mathsf{ZF}].$ The main result of this section (Theorem 8.8) adds
further evidence to this claim, since it shows that similar to $\mathsf{CT}%
_{0}[\mathsf{PA}]$, the purely set-theoretical consequences of $\mathsf{CT}%
_{\ast }[\mathsf{ZFC}]$ can be characterized in terms of iterations of an
appropriate reflection scheme. More specifically, Kotlarski \cite{Kotlarski
bounded} used a mix of model-theoretic and proof-theoretic methods to
characterize the purely arithmetical consequences of $\mathsf{CT}_{0}[%
\mathsf{PA}]$ in terms of $\omega $-iterations of the uniform reflection
scheme. Subsequent proofs of Kotlarski's theorem were given by Beklemishev
and Pakhomov \cite{Beklemishev-Pakhomov} (with proof-theoretic machinery),
and by \L e\l yk \cite{Mateusz-prolongable} (using model-theoretic tools).
Our proof of Theorem 8.8 was directly inspired by \L e\l yk's proof of
Kotlarski's theorem, which we advise the reader to consult. \medskip

\noindent Recall from Definition 7.10 and Lemma 7.11 that for each $n\in 
\mathbb{N},$ there is an $\mathcal{L}_{\mathrm{set}}$-formula $\mathrm{True}%
_{n}(x)$ such that, provably in $\mathsf{KP}$, $\mathrm{True}_{n}(x)$ serves
as a truth predicate for $\mathrm{Depth}_{n}$-sentences in $\mathrm{Sent}^{+}
$ (with no restraint on the constant symbols that appear in them).\footnote{%
We can also work with the truth predicate $\mathrm{True}_{\Sigma _{n}}$ in
this section (see Theorem 3.2.9 and Remark 3.2.10). We opted for $\mathrm{%
True}_{n}$ since the formula defining it in $\mathsf{KP}$ is simpler.}
\medskip 

\noindent \textbf{8.1.~Definition.~}Let $\mathsf{U}$ be a recursive theory
extending $\mathsf{KP}$, let $\mathrm{U}(x)$ be the formula expressing $x\in 
\mathsf{U}$, and 
\begin{equation*}
\mathrm{Prov}_{\mathrm{U}+\mathrm{True}_{n}}(\psi )
\end{equation*}%
be the $\mathcal{L}_{\mathrm{set}}$-formula whose only free variable is $%
\psi $ that expresses:\medskip

\begin{center}
$\psi \in \mathrm{Sent}^{+}$ and $\psi $ is derivable in first order logic
from $\left\{ x:\mathrm{U}(x)\vee \mathrm{True}_{n}(x)\right\} $. \medskip
\end{center}

\noindent We write: 
\begin{equation*}
\mathrm{Con}_{\mathrm{U}+\mathrm{True}_{n}}(\psi )
\end{equation*}%
as shorthand for $\lnot \mathrm{Prov}_{\mathrm{U}+\mathrm{True}_{n}}(\lnot
\psi ).$\medskip

\begin{enumerate}
\item[\textbf{(a)}] $\mathsf{REF}(\mathsf{U}):=\mathsf{U}+\{\mathrm{REF}%
_{n,\varphi }:\varphi \in \mathrm{Form},\ n\in \mathbb{N}\}$, where: \medskip
\end{enumerate}

\begin{center}
$\mathrm{REF}_{n,\varphi }:=\forall x\,\left[ \mathrm{Prov}_{\mathrm{U}+%
\mathrm{True}_{n}}\left( \varphi (\dot{x}))\rightarrow \varphi (x)\right) %
\right] .$\footnote{%
It is well-known that in arithmetical setting, this description of $\mathsf{%
REF(U)}$ is equivalent to the seemingly weaker version that makes no
reference to the internal partial truth predicates $\mathrm{True}_{n}$,
i.e., the collection of sentences of the form $\forall x\,\left[ \mathrm{Prov%
}_{\mathrm{U}}\left( \varphi (\dot{x}))\rightarrow \varphi (x)\right) \right]
,$ as $\varphi (x)$ ranges over arithmetical formulae.}\medskip
\end{center}

\begin{enumerate}
\item[\textbf{(b)}] For each $n\in \mathbb{N},$ $\mathsf{REF}^{n+1}(\mathsf{%
U):=REF}(\mathsf{REF}^{n}(\mathsf{U)).}$\medskip

\item[\textbf{(c)}] $\mathsf{REF}^{\omega }(\mathsf{U}):=\bigcup\limits_{n%
\in \mathbb{N}}\mathsf{REF}^{n}(\mathsf{U}).$\medskip

\item[\textbf{(d)}] $\mathsf{CON}(\mathsf{U}):=\mathsf{U}+\{\mathrm{CON}%
_{n,\varphi }:\varphi (x)\in \mathrm{Form}_{1},\ n\in \mathbb{N}\}$, where: 
\begin{equation*}
\mathrm{CON}_{n,\varphi }:=\forall x\,\left[ \varphi (x)\rightarrow \mathrm{%
Con}_{\mathrm{U}+\mathrm{True}_{n}}(\varphi (\dot{x}))\right] .
\end{equation*}

\item[\textbf{(e)}] For each $n\in \mathbb{N},$ $\mathsf{CON}^{n+1}(\mathsf{U%
}):=\mathsf{CON}(\mathsf{CON}^{n}(\mathsf{U)).}$\medskip

\item[\textbf{(f)}] $\mathsf{CON}^{\omega }(\mathsf{U}):=\bigcup\limits_{n%
\in \mathbb{N}}\mathsf{CON}^{n}(\mathsf{U}).$\medskip
\end{enumerate}

\noindent \textbf{8.2.~Remark.}~It is easy to see that $\mathsf{REF}(\mathsf{%
U})$ and $\mathsf{CON}(\mathsf{U)}$ are deductively equivalent, which in
turn makes it clear that the same goes for $\mathsf{REF}^{\omega }(\mathsf{U}%
)$ and $\mathsf{CON}^{\omega }(\mathsf{U)}$. Since $\mathsf{U}$ is assumed
to be a recursive theory, $\mathsf{REF}^{\omega }(\mathsf{U})$ and $\mathsf{%
CON}^{\omega }(\mathsf{U)}$ are also recursive theories. It is also known
that the deductive closures of $\mathsf{REF}^{\omega }(\mathsf{U})$ and $%
\mathsf{CON}^{\omega }(\mathsf{U)}$ remain the same by replacing $\mathrm{%
True}_{n}$ by $\mathrm{True}_{\Sigma _{n}}$ in their definitions.\medskip

\noindent \textbf{8.3.~Lemma.}~\textit{Let} $\mathsf{U}$ \textit{be a
recursive theory extending }$\mathsf{KP}$\textit{, and} $\mathrm{U}(x)$ 
\textit{be the formula expressing} $x\in \mathsf{U}$\textit{.} \textit{Then
we have}: \medskip

\begin{center}
$\mathsf{CT}^{-}[\mathsf{KP}]+\Delta _{0}^{\mathrm{fin}}$-$\mathsf{Ind(T)}%
+\forall x\left( \mathrm{U}(x)\rightarrow \mathsf{T}(x)\right) \vdash \sigma 
$,\medskip

where $\sigma :=\left[ \forall \varphi \,\forall k\in \omega \,\forall
\alpha \in \mathrm{Asn}(\varphi )\,\left[ \mathrm{Prov}_{\mathrm{U+True}%
_{k}}\left( \varphi \ast \alpha )\rightarrow \mathsf{T}(\varphi \ast \alpha
\right) \right] \right] \mathsf{.}$\footnote{%
Recall that $\varphi \ast \alpha $ is the sentence obtained by replacing
each free variable v in $\varphi $ with the constant $\dot{m}$, where $%
\alpha (v)=m.$}\medskip 
\end{center}

\noindent \textit{In particular,} $\mathsf{CT}^{-}[\mathsf{KP}]+\Delta _{0}^{%
\mathrm{fin}}$-$\mathsf{Ind(T)}+\forall x\,\left( \mathrm{U}(x)\rightarrow 
\mathsf{T}(x)\right) \vdash \mathsf{REF}(\mathsf{U}).$\medskip

\noindent \textbf{Proof.} Recall from Theorem 6.10 that $\mathsf{GRef}_{%
\mathsf{T}}$ is provable in $\mathsf{CT}^{-}[\mathsf{KP}]+\Delta _{0}^{%
\mathrm{fin}}$-$\mathsf{Ind(T)}$. The first statement follows from $\mathsf{%
GRef}_{\mathsf{T}}$ and Lemma 7.13, which assures us that within $\mathsf{CT}%
^{-}[\mathsf{KP}]+\Delta _{0}^{\mathrm{fin}}$-$\mathsf{Ind(T)}$, for all $%
k\in \omega $, $\mathsf{T}\left( \mathrm{True}_{k}(\dot{\sigma})\right)
\leftrightarrow \mathsf{T}\left( \sigma \right) $ for any sentence $\sigma $
in\ \textrm{Sent}$^{+}$ of depth $k.$ The second claim immediately follows
from the first. \hfill $\square $\medskip

\noindent \textbf{8.4.~Corollary}. \textit{Let} $\mathsf{U}$ \textit{be a
recursive theory extending }$\mathsf{KP}$\textit{, and} $\mathrm{U}(x)$ 
\textit{be the formula expressing} $x\in \mathsf{U.}$ \textit{Then we have}%
:\medskip

\begin{center}
$\mathsf{CT}^{-}[\mathsf{KP}]+\Delta _{0}^{\mathrm{fin}}$-$\mathsf{Ind(T)}%
+\forall x\left( \mathrm{U}(x)\rightarrow \mathsf{T}(x)\right) \vdash 
\mathsf{REF}^{\omega }(\mathsf{U}).$\medskip
\end{center}

\noindent \textit{In particular}, $\mathsf{CT}_{\ast }[\mathsf{ZF}]\vdash 
\mathsf{REF}^{\omega }(\mathsf{ZF}).$ \medskip

\noindent \textbf{Proof}. This follows from Lemma 8.3.

\hfill $\square $\medskip

\noindent \textbf{8.5.~Corollary}.\textbf{~}\textit{If }$\mathsf{ZF}^{+}$ 
\textit{is a finite extension of} $\mathsf{ZF}$ (\textit{in the same language%
}),\textit{\ then} $\mathsf{CT}_{\ast }[\mathsf{ZF}^{+}]\vdash \mathsf{REF}%
^{\omega }(\mathsf{ZF}^{+}).$ \textit{In particular}, 
\begin{equation*}
\mathsf{CT}_{\ast }[\mathsf{ZFC}]\vdash \mathsf{REF}^{\omega }(\mathsf{ZFC}).
\end{equation*}

\noindent \textbf{Proof}.\textbf{~}This is an immediate consequence of
Corollary 8.4 and the deduction theorem\textsf{.}

\hfill $\square $\medskip

\noindent The following definition will be used in the proof of Theorem
8.8.\medskip

\noindent \textbf{8.6.~Definition}.\textbf{~}Let $\mathsf{Z}$ be a new unary
predicate, and $\mathsf{GR}(x\mathsf{,T,Z})$ (GR for `Global Reflection') be
the formula in the language $\mathcal{L}_{\mathrm{set}}(\mathsf{T},\mathsf{Z}%
)$ that expresses:\medskip

\begin{center}
$x\in \omega $, $\mathsf{Z}\subseteq \mathrm{Sent}^{+}$, and for all $%
\varphi \in \mathrm{Sent}^{+}$ of depth at most $x$,

if $\varphi $ is derivable from the assumptions in $\mathsf{T}\cup \mathsf{Z}
$, then $\varphi \in \mathsf{T}$. \medskip
\end{center}

\noindent More formally,

\begin{equation*}
\mathsf{GR}(x\mathsf{,T,Z}):=\forall \varphi \in \mathrm{Sent}^{+}\,\left[
\left( \mathrm{Prov}_{\mathrm{T}+\mathrm{Z}}\left( \varphi \right) )\wedge
\varphi \in \mathrm{Depth}_{x}\right) \rightarrow \mathsf{T}(\varphi )\right]
.
\end{equation*}%
\medskip

\noindent \textbf{8.7.~Remark}.\textbf{~}Definition 8.1 makes it clear that
if $\mathcal{M}\models \mathsf{REF}(\mathsf{U})$, then:\medskip

\begin{center}
$\mathbb{\ }\forall n\in \mathbb{N}\ \left( \mathcal{M},\mathrm{True}%
_{n}\right) \models \mathsf{GR}\left( n\mathsf{,T,T}+\mathrm{U}\right) .$%
\footnote{%
The definition of $\mathsf{REF(U)}$ given in Definition 8.1, makes direct
reference to the predicates $\mathrm{True}_{n},$ which makes this fact easy
to see. In the arithmetical context (as in \cite{Beklemishev-Pakhomov} and 
\cite{Mateusz-prolongable}), where the formulation of $\mathsf{REF(U)}$\
does not explicitly mention $\mathrm{True}_{n},$ this property also holds,
but requires a proof.}\medskip
\end{center}

\noindent In particular, if $\mathcal{M}\models \mathsf{REF}^{\omega }(%
\mathsf{ZFC})$, then $\forall n\in \mathbb{N}\ \left( \mathcal{M},\mathrm{%
True}_{n}\right) \models \mathsf{GR}\left( n\mathsf{,T,T}+\mathrm{REF}^{n}%
\mathrm{(ZFC})\right) .$\medskip

\noindent \textbf{8.8.~Theorem.}~\textit{The }$\mathcal{L}_{\mathrm{set}}$-%
\textit{consequences of }$\mathsf{CT}_{\ast }[\mathsf{ZFC}]$ \textit{is
axiomatized by} $\mathsf{REF}^{\omega }(\mathsf{ZFC})$.\medskip

\noindent \textbf{Proof.~}In light of Corollary 8.5, it suffices to show
that every countable model $\mathcal{K}$ of $\mathsf{REF}^{\omega }(\mathsf{%
ZFC})$ has an elementary extension $\mathcal{M}$ such that for some $%
T\subseteq M$, $(\mathcal{M},T)$ satisfies $\mathsf{CT}_{\ast }[\mathsf{ZFC}%
].$ We will show that there is an infinite sequence of structures the
form:\medskip

\begin{center}
$\left\langle \left( \mathcal{M}_{n},T_{n},k_{n}\right) :n\in \mathbb{N}%
\right\rangle $, \medskip
\end{center}

\noindent such that the following conditions hold for all $n\in \mathbb{N}$.
In condition $\mathbf{R}_{2}(n)$ below, we use the notation $\mathsf{CT}%
^{-}\left( \mathsf{F}\right) $ introduced in Definition 3.2.4, for the
formula expressing \textquotedblleft $\mathsf{T}$ satisfies compositional
clauses for formulae in $\mathsf{F}$\textquotedblright . \medskip

$\mathbf{R}_{1}(n):$ $\left( \mathcal{M}_{n},T_{n}\right) \models \mathsf{ZFC%
}(\mathsf{T}).$\medskip

$\mathbf{R}_{2}(n):$ $\left( \mathcal{M}_{n},T_{n}\right) \models \mathsf{CT}%
^{-}\left( \mathrm{Depth}_{k_{n}}\right) .$\medskip

$\mathbf{R}_{3}(n):$ $T_{n}\cap \mathrm{Depth}_{k_{n}}^{\mathcal{M}%
_{n}}=M_{n}\cap \mathrm{Depth}_{k_{n}}^{\mathcal{M}_{n+1}}\cap T_{n+1}$.%
\footnote{%
Note that by condition $\mathbf{R}_{6}(n)$, the $\mathcal{L}_{\mathrm{set}}$%
-formulae (without constants) of depth $k_{n}$ in $\mathcal{M}_{n}$ coincide
with $\mathcal{L}_{\mathrm{set}}$-formulae of depth $k_{n}$ in $\mathcal{M}%
_{n+1}.$}\medskip

$\mathbf{R}_{4}(n):$ $\mathcal{K}\prec \mathcal{M}_{n}\prec \mathcal{M}%
_{n+1}.$\medskip

$\mathbf{R}_{5}(n):$ $k_{n}$ is a nonstandard finite ordinal of $\mathcal{M}%
_{n},$ and $k_{n+1}$ dominates all finite ordinals of $\mathcal{M}_{n}.$%
\medskip

$\mathbf{R}_{6}(n):$ $\left( \omega ,\in \right) ^{\mathcal{M}%
_{n}}\subsetneq _{\mathrm{end}}\left( \omega ,\in \right) ^{\mathcal{M}%
_{n+1}}.$\footnote{%
This condition states that $\left( \omega ,\in \right) ^{\mathcal{M}_{n+1}}$
is a proper end extension of $\left( \omega ,\in \right) ^{\mathcal{M}_{n}}$%
; thus the `new' finite ordinals of $\mathcal{M}_{n+1}$ dominate finite
ordinals of $\mathcal{M}_{n}$.}\medskip

$\mathbf{R}_{7}(n):\left( \mathcal{M}_{n},T_{n}\right) \models \left[
\forall x\,\left( \left( x\in \mathrm{ZFC}\wedge x\in \mathrm{Depth}%
_{k_{n}}\right) \rightarrow \mathsf{T}(x)\right) \right] .$\footnote{%
Recall that we construe $\mathsf{ZFC}$ as the result of enriching a certain
finite set of axioms $\left\{ \varphi _{1},\cdot \cdot \cdot ,\varphi
_{n}\right\} $ with the schemes of separation and collection. Thus, $x\in 
\mathrm{ZFC}$ in the unary formula that expresses \textquotedblleft $x$ is
either an instance of the separation scheme, or an instance of the
collections scheme", or $x=\varphi _{1},$ or $x=\varphi _{2},\cdot \cdot
\cdot ,$ or $x=\varphi _{n}$\textquotedblright .}\medskip

\noindent Note that the proof of the theorem will be complete if there is
such a sequence $\left\langle \left( \mathcal{M}_{n},T_{n},k_{n}\right)
:n\in \mathbb{N}\right\rangle $, since we can then readily show that:\medskip

\begin{center}
$(\mathcal{M}_{\infty },T_{\infty })\models \mathsf{CT}_{\ast }[\mathsf{ZFC}%
] $ and $\mathcal{M\prec M}_{\infty },$\medskip
\end{center}

\noindent where: \medskip

\begin{center}
$\mathcal{M}_{\infty }:=\bigcup\limits_{n\in \mathbb{N}}\mathcal{M}_{n}$,
and $T_{\infty }:=\bigcup\limits_{n\in \mathbb{N}}\widehat{T}_{n}$, where $%
\widehat{T}_{n}:=T_{n}\cap \mathrm{Depth}_{k_{n}}^{\mathcal{M}_{n}}.$
\medskip
\end{center}

\begin{itemize}
\item Our construction of $\left( \mathcal{M}_{0},T_{0},k_{0}\right) $ takes
place in the real world. However, for the later stages of the construction,
we will build $\left( \mathcal{M}_{n+1},T_{n+1},k_{n+1}\right) $ \textit{%
inside} $\left( \mathcal{M}_{n},T_{n},k_{n}\right) $. Since the internal
construction of $\left( \mathcal{M}_{n+1},T_{n+1},k_{n+1}\right) $ within $%
\left( \mathcal{M}_{n},T_{n},k_{n}\right) $ follows the same steps as the
ones that are carried out in the real world for building $\left( \mathcal{M}%
_{0},T_{0},k_{0}\right) ,$ we provide full details of the construction of $%
\left( \mathcal{M}_{0},T_{0},k_{0}\right) $ so that the reader gains an
intuition of the later stages of construction.\medskip\ 

\item One might think that all that is required of $\left( \mathcal{M}%
_{0},T_{0},k_{0}\right) $ is that $\mathcal{K\prec M}_{0}$, and for some
nonstandard $k_{0}\in \omega ^{\mathcal{M}_{0}},$ there is a $\mathrm{Depth}%
_{k_{0}}^{\mathcal{M}_{0}}$-truth class $T_{0}$ over $\mathcal{M}_{0}$
satisfying condition $(\ast )$ below:\medskip
\end{itemize}

\noindent $(\ast )$ $\ \ \left( \mathcal{M}_{0},T_{0},k_{0}\right) \models
\forall x\,\left[ \left( x\in \mathrm{ZFC}\wedge x\in \mathrm{Depth}%
_{k_{0}}\right) \rightarrow \mathsf{T}(x)\right] .$\medskip

\noindent However, $(\ast )$ turns out to be insufficient for our purposes.
We will instead arrange the much stronger condition $(\ast \ast )$ below,
which will enable us to use $\left( \mathcal{M}_{0},T_{0},k_{0}\right) $ as
the `launching pad' of our construction of length $\omega $, whose `engine'
is Lemma $(\nabla ).$\medskip 

\noindent $(\ast \ast )$ $\ \ \left( \mathcal{M}_{0},T_{0},k_{0}\right)
\models \mathsf{GR}(k_{0}\mathsf{,T,T}+\mathrm{REF}^{k_{0}}(\mathrm{ZFC})).$%
\footnote{%
The binary formula $x\in \mathsf{REF}^{y}(\mathrm{ZFC})$ is defined by an
internal recursion within $\mathsf{ZF}$.}\medskip

\noindent For this purpose, consider the following collection of sentences
in the language obtained by enriching $\mathcal{L}_{\mathrm{set}}$ with a
new predicate $\mathsf{T}$, as well as constants for each element of $K$,
and a fresh constant $c$. \medskip

\noindent $\Gamma _{1}:=\mathrm{ED}(\mathcal{K}).$\medskip

\noindent $\Gamma _{2}:=\mathsf{ZF}(\mathsf{T})$. \medskip

\noindent $\Gamma _{3}:=$ $\left\{ c\in \omega \right\} \cup \{n\in c:n\in 
\mathbb{N}\}$. \medskip

\noindent $\Gamma _{4}:=\mathsf{CT}^{-}\left( \mathrm{Depth}_{c}\right) .$%
\medskip

\noindent $\Gamma _{5}:=\mathsf{GR}(c\mathsf{,T,T}+\mathrm{REF}^{c}\mathrm{%
(ZFC)}).$\medskip

\noindent Let $\Gamma $ be the union of $\Gamma _{1}$ through $\Gamma _{5}.$
We will show that $\Gamma $ is consistent by verifying the consistency of
every finite subset of $\Gamma $. Note that since $\mathcal{K}$ is a model
of $\mathsf{REF}^{\omega }(\mathsf{ZFC})$, for each $n\in \mathbb{N}$, the
structure $\left( \mathcal{K},\mathrm{True}_{n}^{\mathcal{K}}\right) $
satisfies the following properties:\medskip 

\noindent (a) $\mathsf{ZF}(\mathsf{T}).$\medskip

\noindent (b) $\mathsf{CT}^{-}\left( \mathrm{Depth}_{n}\right) .$\medskip

\noindent (c) $\mathsf{GR}(n\mathsf{,T,T}+\mathrm{REF}^{n}(\mathrm{ZFC}))$%
.\medskip

\noindent Property (a) holds since $\mathrm{True}_{n}^{\mathcal{K}}$ is
definable in $\mathcal{K}$, property (b) holds by Corollary 7.13, and
property (c) holds by Remark 8.7. This makes it clear that given an
arbitrary finite subset $\Gamma ^{\prime }$ of $\Gamma $, if $n^{\prime }$
is chosen as the largest natural number mentioned in $\Gamma ^{\prime },$
then we have:%
\begin{equation*}
\left( \mathcal{K},\mathrm{True}_{n^{\prime }}^{\mathcal{K}}\right) \models
\Gamma ^{\prime }.
\end{equation*}%
This concludes our verification of the consistency of $\Gamma $. Therefore
there is some countable model $(\mathcal{M}_{0},T_{0},k_{0})$ of $\Gamma $
(where $T_{0}$ is the interpretation of $\mathsf{T}$, and $k_{0}$ is the
interpretation of $c$). In particular, we have:%
\begin{equation*}
\left( \mathcal{M}_{0},T_{0},k_{0}\right) \models \left[ \mathsf{CT}%
^{-}\left( \mathrm{Depth}_{k_{0}}\right) \wedge \mathsf{GR}(k_{0}\mathsf{,T,T%
}+\mathrm{REF}^{k_{0}}(\mathrm{ZFC}))\right] .
\end{equation*}

\noindent We now prove Lemma $\mathbf{(}\nabla )$ below, which provides the
engine for the recursive construction of the desired $\left( \mathcal{M}%
_{n},T_{n},k_{n}\right) $ for $n\geq 1$.\footnote{%
Lemma\textbf{\ }$\mathbf{(}\nabla )$ is the analogue of Lemma 5.6 of \cite%
{Mateusz-prolongable}; the main difference in their proofs is the use of
forcing here to arrange a global well-ordering so that the set-theoretical
counterpart of the arithmetical completeness theorem can be carried out.}%
\medskip

\noindent \textbf{Lemma }$\mathbf{(}\nabla )$\textbf{.~}\textit{Suppose} $%
\mathcal{M}$ \textit{is a countable model of} $\mathsf{ZFC}$, $T\subseteq M,$
$k$\textit{\ is a nonstandard element of} $\omega ^{\mathcal{M}}$, \textit{%
and the following three conditions hold}:\medskip

\begin{enumerate}
\item[$(1)$] $(\mathcal{M},T)\models \mathsf{ZF(T)}$. \medskip

\item[$(2)$] $(\mathcal{M},T)\models \mathsf{CT}^{-}\left( \mathrm{Depth}%
_{k}\right) $. \medskip

\item[$(3)$] $(\mathcal{M},T)\models \mathsf{GR}(k\mathsf{,T,T}+\mathrm{REF}(%
\mathrm{U}))).$ \medskip
\end{enumerate}

\noindent \textit{Then} \textit{there is some countable model} $\mathcal{M}%
^{\ast }\succ \mathcal{M}$, \textit{such that }$\left( \omega ,\in \right) ^{%
\mathcal{M}}\subsetneq _{\mathrm{end}}\left( \omega ,\in \right) ^{\mathcal{M%
}^{\ast }}$, \textit{and for some new finite ordinal in} $k^{\ast }\in
\omega ^{\mathcal{M}^{\ast }}$, \textit{and some} $T^{\ast }\subseteq
M^{\ast }$, \textit{the following four conditions are satisfied}:\medskip 

\begin{enumerate}
\item[$(i)$] $(\mathcal{M}^{\ast },T^{\ast })\models \mathsf{ZF(T)}$.\medskip

\item[$(ii)$] $(\mathcal{M}^{\ast },T^{\ast })\models \mathsf{CT}^{-}\left( 
\mathrm{Depth}_{k^{\ast }}\right) .$\medskip

\item[$(iii)$] $(\mathcal{M}^{\ast },T^{\ast })\models \mathsf{GR}(k^{\ast }%
\mathsf{,T,T}+\mathrm{U}).$\footnote{%
Note that the superscript here is $p-1$, whereas the superscript in (3) is $%
p $. Since $p$ is a nonstandard finite ordinal, this assures us that for
each $n\in \mathbb{N},$ $n$ application of the lemma results in $p-n$, which
remains a nonstandard finite ordinal.}\medskip

\item[$(iv)$] $T\cap \mathrm{Depth}_{k}^{\mathcal{M}}=T^{\ast }\cap \mathrm{%
Depth}_{k}^{\mathcal{M}^{\ast }}\cap M.$\footnote{%
Note that since $\left( \omega ,\in \right) ^{\mathcal{M}}\subsetneq _{%
\mathrm{end}}\left( \omega ,\in \right) ^{\mathcal{M}^{\ast }}$ is among the
conclusions of the Lemma, the $\mathcal{L}_{\mathrm{set}}$-formulae (without
constants) of depth $k$ in $\mathcal{M}$ coincide with $\mathcal{L}_{\mathrm{%
set}}$-formulae of depth $k$ in $\mathcal{M}^{\ast }.$}\medskip
\end{enumerate}

\noindent \textbf{Proof.~}Let $\left( \mathcal{M},T,k\right) $ be as in the
assumptions of the Lemma. By condition $(3)$ of the lemma, the `theory'\ $%
\mathsf{T}+\mathrm{REF}(\mathrm{U})$ is consistent from the point of view of 
$(\mathcal{M},T)$. However, $\mathsf{T}+\mathrm{REF}(\mathrm{U})$ is a
proper class of sentences of $\mathcal{M}$ (since all the constants naming
the elements of the universe of $\mathcal{M}$ occur in $\mathsf{T}$), we
cannot use the completeness theorem of first order logic for such a large
consistent `theory'\ in $\mathsf{ZF(T)}$ alone. The obstacle we are facing
is due to the fact that in order to carry out the Henkin proof of a
class-sized theory, we need to have access to a global well-ordering within $%
\left( \mathcal{M},T\right) .$ Such a well-ordering is available if we
further assume that \textrm{U} includes the axiom $\mathsf{V}=\mathsf{L}$,
or more generally $\exists p\left( \mathsf{V}=\mathsf{HOD}(p)\right) $, but
it is well-known that $\mathsf{ZFC}$ alone does not guarantee the existence
of a global well-ordering. \medskip 

\noindent There is a `magical' way to circumvent the above obstacle: we can
use forcing to expand $\left( \mathcal{M},T\right) $ to: 
\begin{equation*}
\left( \mathcal{M},T,<_{M}\right) \models \mathsf{ZF}(\mathsf{T,<})+\mathsf{%
GW}(<),
\end{equation*}%
where $\mathsf{GW}(<)$ is the sentence asserting that $<$ is a set-like%
\footnote{%
In other words, every proper initial segment of $<$ is a set (and not a
proper class).} linear order of the universe, and every nonempty set has a $%
\mathsf{<}$-least element.\footnote{%
The existence of such an ordering $<_{M}$ was proved by Felgner \cite%
{Felgner} for countable models $\mathcal{M}$ of $\mathsf{ZFC}$, but the
proof works equally well for models of $\mathsf{ZFC}+\mathsf{Sep}(\mathcal{L}%
)+\mathsf{Coll}(\mathcal{L})$ for any finite language $\mathcal{L}$
extending $\mathcal{L}_{\mathrm{set}}.$ Also notice that the existence of
such an expansion $(\mathcal{M},T,<_{M})$ is equivalent to arranging an
expansion $(\mathcal{M},T,f)$ satisfying $\mathsf{ZF}$ in the extended
language such that $f$ is a bijection between $M$ and $\mathrm{Ord}^{%
\mathcal{M}}$.} This allows $\left( \mathcal{M},T,<_{M}\right) $ to define a
model $\mathcal{K}$ of $\mathsf{T}+\mathrm{REF}(\mathrm{U})$ with the
additional bonus that the entire elementary diagram of $\mathcal{K}$
(incorporating also nonstandard sentences of $\mathcal{M}$) is definable in $%
\left( \mathcal{M},T,<_{M}\right) .$ In other words, $\left( \mathcal{M}%
,T,<_{M}\right) $ \textit{strongly interprets} a model $\mathcal{K}$ of $%
\mathsf{T}+\mathrm{REF}(\mathrm{U}),$ i.e., there is a definable class in $%
\left( \mathcal{M},T,<_{M}\right) $ that serves as $\mathrm{ED}(\mathcal{K}).
$ Thanks to condition $(2)$ and the nonstandardness of $k$, from the real
world point of view $\mathsf{T}$ includes the elementary diagram of $%
\mathcal{M}$, thus viewed from the real world, $\mathcal{M\prec K}.$\medskip 

\noindent Next, we carry out an \textit{internal variant} of the compactness
argument used earlier for the construction of $(\mathcal{M}_{0},T_{0},k_{0})$
in order to construct -- within $\left( \mathcal{M},T,<_{M},k\right) $ --
the desired structure $\left( \mathcal{M}^{\ast },T^{\ast },k^{\ast }\right)
.$\footnote{%
By arranging $\mathcal{M}^{\ast }$ to have countable cofinality as viewed
from $\left( \mathcal{M},T,<_{M}\right) $, it is possible to interpret $%
\left( \mathcal{M}^{\ast },T^{\ast },<_{M^{\ast }},k^{\ast }\right) $ in $%
\left( \mathcal{M},T,<_{M},k\right) ,$ where $<_{M^{\ast }}$ is a generic
global well--order, but we do not need this extra flourish for our
construction.} For this purpose we consider the analogue $\Gamma ^{\ast }$
of the set $\Gamma $ of sentences used in the construction of $(\mathcal{M}%
_{0},T_{0},k_{0}).$ $\Gamma ^{\ast }$ is defined \textit{within} $\left( 
\mathcal{M},T,<_{M}\right) $ as the union of the following collection of
sentences $\Gamma _{1}^{\ast }$ through $\Gamma _{5}^{\ast }$. Thus $\Gamma
^{\ast }$ is formulated in the language obtained by enriching $\mathcal{L}_{%
\mathrm{set}}$ with a new predicate $\mathsf{T}$, constants for each element
of $K$, and a fresh constant $c$. \medskip 

\noindent $\Gamma _{1}^{\ast }:=\mathrm{ED}(\mathcal{K}).$\medskip

\noindent $\Gamma _{2}^{\ast }:=\mathsf{ZF}(\mathsf{T})$. \medskip

\noindent $\Gamma _{3}^{\ast }:=$ $\left\{ c\in \omega \right\} \cup \{n\in
c:n\in \mathbb{N}\}$. \medskip

\noindent $\Gamma _{4}^{\ast }:=\mathsf{CT}^{-}\left( \mathrm{Depth}%
_{c}\right) .$\medskip

\noindent $\Gamma _{5}^{\ast }:=\mathsf{GR}(c\mathsf{,T,T}+\mathrm{REF}^{c}(%
\mathrm{ZFC})).$\medskip

\noindent $\left( \mathcal{M},T,<_{M}\right) $ has access to the full
elementary diagram of $\mathcal{K}$, and it can verify that every finite
subset of $\Gamma ^{\ast }$ can be interpreted in $\mathcal{K}$. Thus $%
\left( \mathcal{M},T,<_{M}\right) $ views $\Gamma ^{\ast }$ as a consistent
theory. More specifically, $\mathcal{K}$ is a model of all the axioms of $%
\mathsf{ZF}$ (including nonstandard ones) in the eyes of $(\mathcal{M}%
,T,<_{M})$, and therefore for each $p\in \omega ^{\mathcal{K}}$, the
`formula' $\mathrm{True}_{p}\in \mathrm{Form}^{\mathcal{K}}$ (as in
Definition 7.10) will be assessed by $(\mathcal{M},T,<_{M})$ to give rise to
a $\mathrm{Depth}_{p}$-truth class $T_{p}$ on $\mathcal{K}$. Moreover, since 
$\mathrm{ED}(\mathcal{K})$ is definable in $\left( \mathcal{M}%
,T,<_{M}\right) ,$ and as seen by $\left( \mathcal{M},T,<_{M}\right) $, $%
T_{p}$ is definable in $\mathcal{K}$, we can conclude that:\medskip 

\begin{center}
$\left( \mathcal{M},T,<_{M}\right) \models $ $[$\textquotedblleft $(\mathcal{%
K},T_{k})\models \mathsf{ZF(T)}$\textquotedblright $],$ \medskip
\end{center}

\noindent or more formally, $\left( \mathcal{M},T,<_{M}\right) \models \left[
\forall \psi \,\left( \psi \in \mathrm{ZF(T)}\ \psi \in \mathrm{ED}(\mathcal{%
K},T_{k})\right) \right] .$ Moreover, since \medskip

\begin{center}
$\left( \mathcal{M},T,<_{M}\right) \models $ $[$\textquotedblleft $\mathcal{K%
}\models \mathsf{T}+\mathrm{REF}(\mathrm{U})$\textquotedblright $]$, \medskip
\end{center}

\noindent as in Remark 8.7, $\left( \mathcal{M},T,<_{M}\right) \models
\lbrack \forall p\in \omega $ \textquotedblleft $\left( \mathcal{K},\mathrm{%
True}_{p}^{\mathcal{K}}\right) \models \mathrm{GR}\left( p\mathsf{,T,T}+%
\mathrm{U}\right) $\textquotedblright ]$.$\medskip

\noindent We can therefore conclude that $\left( \mathcal{M},T,<_{M}\right) $
sees $\Gamma ^{\ast }$ as a consistent proper class of sentences. The global
well-ordering $<_{M}$, now comes handy for a second time, since it allows $%
\left( \mathcal{M},T,<_{M}\right) $ to strongly interpret a structure $%
\left( \mathcal{M}^{\ast },T^{\ast },k^{\ast }\right) $ satisfying $\Gamma
^{\ast }$ (where $k^{\ast }$ is the interpretation of the constant symbol $c$%
)$.$ Note that by general considerations\footnote{%
The argument here is similar to the argument used in showing that
conservative extensions of models of $\mathsf{PA}$ are end extensions.}, the
definability of $\mathcal{M}^{\ast }$ in $(\mathcal{M},T,<_{M})$ implies
that all the new finite ordinals of $\mathcal{M}^{\ast }$ exceed all the
finite ordinals of $\mathcal{M}$, so $\left( \omega ,\in \right) ^{\mathcal{M%
}}\subsetneq _{\mathrm{end}}\left( \omega ,\in \right) ^{\mathcal{M}^{\ast
}}.$ In summary, the internal version of the completeness theorem was used
twice in $\left( \mathcal{M},T,<_{M}\right) $, first to build $\mathcal{K}$,
and then to build $(\mathcal{M}^{\ast },T^{\ast },k^{\ast })$. This
concludes the proof of Lemma $\mathbf{(}\nabla ).$

\begin{itemize}
\item Since $\mathrm{REF}^{p}(\mathrm{ZFC}))=\mathrm{REF}(\mathrm{REF}^{p-1}(%
\mathrm{ZFC}))$, by using $\mathrm{U}:=\mathrm{REF}^{p-1}(\mathrm{ZFC}))$ in
the Lemma $\mathbf{(}\nabla )$, we obtain Lemma $\mathbf{(}\heartsuit )$ as
an immediate consequence.\medskip
\end{itemize}

\noindent \textbf{Lemma }$\mathbf{(}\heartsuit )$\textbf{.~}\textit{Suppose} 
$\mathcal{M}$ \textit{is a countable }$\omega $-\textit{nonstandard model of}
$\mathsf{ZFC}$, $k$ \textit{and} $p$ \textit{are nonstandard elements of} $%
\omega ^{\mathcal{M}}$, \textit{and the following conditions hold:}\medskip

\begin{enumerate}
\item[$(1)$] $(\mathcal{M},T)\models \mathsf{ZF(T)}$. \medskip

\item[$(2)$] $(\mathcal{M},T)\models \mathsf{CT}^{-}\left( \mathrm{Depth}%
_{k}\right) $. \medskip

\item[$(3)$] $(\mathcal{M},T)\models \mathsf{GR}(k\mathsf{,T,T}+\mathrm{REF}%
^{p}(\mathrm{ZFC}))).$ \medskip
\end{enumerate}

\noindent \textit{Then} \textit{there is some countable model} $\mathcal{M}%
^{\ast }\succ \mathcal{M}$, \textit{such that }$\left( \omega ,\in \right) ^{%
\mathcal{M}}\subsetneq _{\mathrm{end}}\left( \omega ,\in \right) ^{\mathcal{M%
}^{\ast }}$, \textit{and for some new finite ordinal in} $k^{\ast }\in
\omega ^{\mathcal{M}^{\ast }}$, \textit{and some} $T^{\ast }\subseteq
M^{\ast }$, \textit{the following conditions are satisfied}:\medskip 

\begin{enumerate}
\item[$(i)$] $(\mathcal{M}^{\ast },T^{\ast })\models \mathsf{ZF(T)}$.\medskip

\item[$(ii)$] $(\mathcal{M}^{\ast },T^{\ast })\models \mathsf{CT}^{-}\left( 
\mathrm{Depth}_{k^{\ast }}\right) .$\medskip

\item[$(iii)$] $(\mathcal{M}^{\ast },T^{\ast })\models \mathsf{GR}(k^{\ast }%
\mathsf{,T,T}+\mathrm{REF}^{p-1}(\mathrm{ZFC})).$\footnote{%
Note that the superscript here is $p-1$, whereas the superscript in (3) is $p
$. Since $p$ is a nonstandard finite ordinal, this assures us that for each $%
n\in \mathbb{N},$ $n$ applications of the lemma results in $p-n$, which
remains a nonstandard finite ordinal.}\medskip 

\item[$(iv)$] $T\cap \mathrm{Depth}_{k}^{\mathcal{M}}=T^{\ast }\cap \mathrm{%
Depth}_{k}^{\mathcal{M}^{\ast }}\cap M.$\footnote{%
Since $\left( \omega ,\in \right) ^{\mathcal{M}}\subsetneq _{\mathrm{end}%
}\left( \omega ,\in \right) ^{\mathcal{M}^{\ast }}$ is among the conclusions
of the lemma, the $\mathcal{L}_{\mathrm{set}}$-formulae (without constants)
of depth $k$ in $\mathcal{M}$ coincide with $\mathcal{L}_{\mathrm{set}}$%
-formulae of depth $k$ in $\mathcal{M}^{\ast }.$}\medskip
\end{enumerate}

\noindent Starting with $\left( \mathcal{M}_{0},T_{0},k_{0}\right) $, which
earlier in the proof\footnote{%
The special case of (d) for $n=0$ was tagged as $(\ast \ast )$.} was
arranged to satisfy (d) below for $n=0$, repeated applications of Lemma $%
\mathbf{(}\heartsuit )$ with $p=k_{0}$, yields the sequence $\left\langle
\left( \mathcal{M}_{n},T_{n},k_{n}\right) :n\in \mathbb{N}\right\rangle $
such that the following holds for each $n\in \mathbb{N}:$\medskip

\begin{enumerate}
\item[$(a)$] $(\mathcal{M}_{n},T_{n})\models \mathsf{ZF(T).}$\medskip

\item[$(b)$] $(\mathcal{M}_{n},T_{n})\models \mathsf{CT}^{-}\left( \mathrm{%
Depth}_{k_{n}}\right) .$\medskip

\item[$(c)$] $\left( \omega ,\in \right) ^{\mathcal{M}_{n}}\subsetneq _{%
\mathrm{end}}\left( \omega ,\in \right) ^{\mathcal{M}_{n+1}},$ $k_{n+1}\in
\omega ^{\mathcal{M}_{n+1}}$, and $k_{n+1}\notin M_{n}.$\medskip 

\item[$(d)$] $(\mathcal{M}^{\ast },T^{\ast })\models \mathsf{GR}(k_{n}%
\mathsf{,T,T}+\mathrm{REF}^{k_{0}-n}(\mathrm{ZFC})).$\medskip 

\item[$(e)$] $T\cap \mathrm{Depth}_{k_{n}}^{\mathcal{M}_{n}}=T_{n+1}\cap 
\mathrm{Depth}_{k_{n}}^{\mathcal{M}_{n+1}}\cap M_{n}.$\medskip 
\end{enumerate}

\noindent This makes it clear that conditions $\mathbf{R}_{1}(n)$ through $%
\mathbf{R}_{8}(n)$ (specified early in the proof) are met for all $n\in 
\mathbb{N}.$ As explained earlier, this is sufficient to establish Theorem
8.8.

\hfill $\square $\medskip

\noindent \textbf{8.9.~Remark.~}Using the methodology of arithmetizing the
model-theoretic proof of conservativity of $\mathsf{CT}^{-}[\mathsf{PA}]$
over $\mathsf{PA}$ used in \cite{Trio on feasible red.}, or the one in \cite%
{Ali-curious}, one should be able to show that Theorem 8.8 can be verified
in $\mathsf{WKL}_{0}$, and therefore in Primitive Recursive
Arithmetic.\medskip

\noindent \textbf{8.10.~Remark.~}The proof strategy of Theorem 8.8 can be
used to show (without any appeal to forcing) that if $\mathsf{U}$ is a
recursive extension of $\mathsf{KP}+\mathsf{V}=\mathsf{L}$, then the purely
set-theoretical consequences of $\mathsf{CT}^{-}[\mathsf{U}]+\Delta _{0}^{%
\mathrm{fin}}$-$\mathsf{Ind(T)}+\forall x\left( \mathrm{U}(x)\rightarrow 
\mathsf{T}(x)\right) $ coincides with $\mathsf{REF}^{\omega }(\mathsf{U})$.

\bigskip \bigskip

\begin{center}
\textbf{9.~BETWEEN} $\mathsf{CT}_{\ast }[\mathsf{ZF}]$ \textbf{AND} $\mathsf{%
CT}_{0}[\mathsf{ZF}]$\bigskip
\end{center}

\noindent Recall from Section 5 that $\mathsf{FRef}$ is provable in $\mathsf{%
CT}_{0}[\mathsf{ZF}]\mathsf{,}$ where:\medskip

\begin{center}
$\mathsf{CT}_{0}[\mathsf{ZF}]:=\mathsf{CT}^{-}[\mathsf{ZF}]+\mathsf{\Delta }%
_{0}$-$\mathsf{Sep(T)+\Delta _{0}}$-$\mathsf{Coll(T).}$\medskip
\end{center}

\noindent In this section we examine the strength of $\mathsf{CT}_{\ast }[%
\mathsf{ZF}]+\Delta _{0}$-$\mathsf{Sep}(\mathsf{T})$, and $\mathsf{CT}_{\ast
}[\mathsf{ZF}]+\mathsf{FRef}$. As we will see, these theories lie strictly
between $\mathsf{CT}_{\ast }[\mathsf{ZF}]$ and $\mathsf{CT}_{0}[\mathsf{ZF}]$%
, both in deductive power, and consistency strength. \medskip

\noindent \textbf{9.1.~Theorem.~}\textit{Over} $\mathsf{CT}^{-}[\mathsf{KP}]$%
, \textit{the following are deductively equivalent}:\medskip

\begin{enumerate}
\item[$(a)$] $\forall x\ \exists y\ \left( y=x\cap \mathsf{T}\right) .$%
\footnote{%
This condition is often read as \textquotedblleft $\mathsf{T}$ is piecewise
coded\textquotedblright . In the context of set theory, this condition can
be thought of expressing that $\forall x\ \left( \mathrm{Th}\mathsf{(}%
\mathrm{V}\mathsf{,\in ,}a\mathsf{)}_{a\in x}\in \mathrm{V}\right) .$ Here $%
\mathsf{(}\mathrm{V}\mathsf{,\in ,}a\mathsf{)}_{a\in x}$ is the result of
expanding $\mathsf{(}\mathrm{V}\mathsf{,\in )}$ by the elements of $x.$
Thus, $\mathsf{(}\mathrm{V}\mathsf{,\in ,}a\mathsf{)}_{a\in x}$ is an $%
\mathcal{L}$-structure, where $\mathcal{L}$ is the result of adding
constants for the elements of $x$ to $\mathcal{L}_{\mathrm{set}}.$}\medskip

\item[$(b)$] $\Delta _{0}$-$\mathsf{Sep}(\mathsf{T}).$\medskip
\end{enumerate}

\noindent \textbf{Proof.~}$(a)\Rightarrow (b)$ is established by a
straightforward induction of the depth of $\Delta _{0}(\mathsf{T})$%
-formulae. The closure of the universe under the operation of relative
complementation takes care of the negation step; the closure of the universe
under unions takes care of the disjunction step; and the closure of the
universe under Cartesian products and projections takes care of the
existential step. $(b)\Rightarrow (a)$ is trivial.

\hfill $\square $\medskip

\noindent \textbf{9.2.~Remark.~}Clearly the consistency of $\mathsf{ZF}$ is
provable in $\mathsf{CT}_{\ast }[\mathsf{ZF}]$. Also, note that every $%
\omega $-model of $\mathsf{CT}^{-}[\mathsf{ZF}]$ is a model of $\mathsf{CT}%
^{\ast }[\mathsf{ZF}]+\mathsf{Ind}(\mathsf{T}).$\medskip

\noindent \textbf{9.3.~Proposition.~}\textit{Assuming the consistency of }$%
\mathsf{CT}_{\ast }[\mathsf{ZF}],$ $\mathsf{CT}_{\ast }[\mathsf{ZF}]$ 
\textit{cannot prove that} $\mathsf{ZF}$ \textit{has an} $\omega $-\textit{%
model.}\medskip 

\noindent \textbf{Proof.~}This follows from putting the second assertion in
Remark 9.2 together with G\"{o}del's second incompleteness theorem.

\hfill $\square $\medskip

\noindent \textbf{9.4.~Remark.~}As we shall explain, assuming that $\mathsf{%
ZF}$\ has an $\omega $-model, there is an $\omega $-model $\mathcal{M}$ of $%
\mathsf{ZF}$ that satisfies \textquotedblleft $\mathsf{ZF}$\ has no $\omega $%
-model\textquotedblright . By considering $(\mathcal{M},T)$, where $T$ is
the elementary diagram of $\mathcal{M}$, we obtain a stronger form of
Proposition 9.3, since $(\mathcal{M},T)$\ \textit{is an} $\omega $-\textit{%
model of }$\mathsf{CT}_{\ast }[\mathsf{ZF}]$ \textit{in which there is no} $%
\omega $-\textit{model of} $\mathsf{ZF}$. The existence of the desired $%
\omega $-model $\mathcal{M}$ follows from the abstract form of G\"{o}del's
second incompleteness theorem that states that if $\Gamma $ is a consistent
theory extending Robinson's $\mathsf{Q}$ that supports a unary predicate $%
\theta (x)$ satisfying the three Hilbert-Bernays-L\"{o}b provability
conditions below, then $\Gamma \nvdash $ $\theta (\ulcorner 0=1\urcorner )$.%
\footnote{%
See, e.g., \cite[Ch.~18]{Boolos et al.}, for the presentation of such a
general form of G\"{o}del's second incompleteness theorem.}\medskip 

\begin{enumerate}
\item[HBL-1] $\Gamma \vdash \varphi \Longrightarrow \Gamma \vdash \mathsf{%
\theta }(\ulcorner \varphi \urcorner )$.\medskip

\item[HBL-2] $\Gamma \vdash \mathsf{\theta }(\ulcorner \varphi \rightarrow
\psi \urcorner )\rightarrow (\mathsf{\theta }(\ulcorner \varphi \urcorner
)\rightarrow \mathsf{\theta }(\ulcorner \psi \urcorner ))$.\medskip

\item[HBL-3] $\Gamma \vdash \mathsf{\theta }(\ulcorner \varphi \urcorner
)\rightarrow \mathsf{\theta }(\ulcorner \mathsf{\theta (}\ulcorner \varphi
\urcorner )\urcorner )$.\medskip
\end{enumerate}

\noindent More explicitly, let $\Omega _{\mathrm{ZF}}(x)$ be the $\mathcal{L}%
_{\mathrm{set}}$-formula that expresses \textquotedblleft $x$ is an $\omega $%
-model of $\mathsf{ZF}$\textquotedblright ; $\Gamma :=\mathsf{ZF}$, and let $%
\theta (v)$ be the $\mathcal{L}_{\mathrm{set}}$-formula expressing
\textquotedblleft $v$ is true in all $\omega $-models of $\mathsf{ZF}$%
\textquotedblright , i.e., the formula:\medskip

\begin{center}
\textquotedblleft $v$ is the (code of) an $\mathcal{L}_{\mathrm{set}}$%
-sentence and $\forall x\,(\Omega _{\mathrm{ZF}}(x)\rightarrow v\in \mathrm{%
Th}(x)).$ \medskip
\end{center}

\noindent Then conditions HBL-1 through HBL-3 are straightforward to verify,
thanks to provability in $\mathsf{ZF}$ of the following statement:\medskip

\begin{center}
For all $\mathcal{L}_{\mathrm{set}}$-structures $\mathcal{M}$, if $\Omega _{%
\mathrm{ZF}}(\mathcal{M})$, then: \medskip

for all $\mathcal{N}\in \mathcal{M}$, if $\mathcal{M}\models \Omega _{%
\mathrm{ZF}}(\mathcal{N})$, then $\Omega _{\mathrm{ZF}}(\mathcal{M})$.%
\footnote{%
This formula expresses \textquotedblleft an $\omega $-model of $\mathsf{ZF}$
of an $\omega $-model of $\mathsf{ZF}$, is an $\omega $-model of $\mathsf{ZF}
$\textquotedblright .}\medskip
\end{center}

\noindent \textbf{9.5.~Remark.~}We saw in Theorem 5.7 that $\mathsf{FRef}$
is provable in $\mathsf{CT}_{0}[\mathsf{ZF}]$. On the other hand, Theorem
9.1 makes it clear that $\mathsf{FRef}$ implies $\Delta _{0}$-$\mathsf{Sep}(%
\mathsf{T})$ in the presence of $\mathsf{CT}^{-}[\mathsf{ZF}].$ Thus we
have:\medskip

\begin{center}
$\mathsf{CT}_{0}[\mathsf{ZF}]\ \vdash \ \mathsf{CT}^{-}[\mathsf{ZF}]+\mathsf{%
FRef}\ \vdash \ \mathsf{CT}_{\ast }[\mathsf{ZF}]+\Delta _{0}$-$\mathsf{Sep}(%
\mathsf{T})\ \vdash \ \mathsf{CT}_{\ast }[\mathsf{ZF}].$\medskip
\end{center}

\noindent The above will be refined in Theorem 9.10.\medskip

\noindent \textbf{9.6.~Examples}. \medskip

\begin{enumerate}
\item[\textbf{(a)}] \textit{Separating} $\mathsf{CT}_{0}[\mathsf{ZF}]$ 
\textit{from} $\mathsf{CT}_{\ast }[\mathsf{ZF}]+\mathsf{FRef}$: Let $\kappa $
be a strongly inaccessible cardinal. It is well-known that there is a closed
unbounded subset $C$ of $\kappa $ such that:%
\begin{equation*}
C=\left\{ \alpha <\kappa :\left( \mathrm{V\!}_{\alpha },\in \right) \prec
\left( \mathrm{V\!}_{\kappa },\in \right) \right\} .
\end{equation*}%
Enumerate $C$ in increasing order as $\left\langle \alpha _{\delta }:\delta
\in \kappa \right\rangle $. Let $T$ be the elementary diagram of $\left( 
\mathrm{V\!}_{\alpha _{\omega }},\in \right) .$ Then $\left( \mathrm{V\!}%
_{\alpha _{\omega }},\in ,T\right) $ satisfies $\mathsf{CT}[\mathsf{ZF}]+%
\mathsf{FRef+Sep(T),}$ but it is not a model of $\mathsf{CT}_{0}[\mathsf{ZF}%
] $ since the map $n\mapsto \alpha _{n}$, which maps $\omega $ cofinally
into $\alpha _{\omega }$, is $\Delta _{0}(\mathsf{T})$-definable in $\left( 
\mathrm{V\!}_{\alpha _{\omega }},\in ,T\right) .$\medskip

\item[\textbf{(b)}] \textit{Separating} $\mathsf{CT}_{\ast }[\mathsf{ZF}]+%
\mathsf{FRef}$ \textit{from} $\mathsf{CT}_{\ast }[\mathsf{ZF}]+\Delta _{0}$-$%
\mathsf{Sep}(\mathsf{T})$: Let $\alpha $ be the first ordinal such that $%
\left( \mathrm{V\!}_{\alpha },\in \right) $ is a model of $\mathsf{ZF}$, and
let $T$ be the elementary diagram of $\left( \mathrm{V\!}_{\alpha },\in
\right) .$ Then $\left( \mathrm{V\!}_{\alpha },\in ,T\right) $ satisfies $%
\mathsf{CT}_{\ast }[\mathsf{ZF}]+\mathsf{Sep(T)},$ but it does not satisfy $%
\mathsf{FRef}$ by the choice of $\alpha .$\medskip

\item[\textbf{(c)}] \textit{Separating} $\mathsf{CT}_{\ast }[\mathsf{ZF}%
]+\Delta _{0}$-$\mathsf{Sep}(\mathsf{T})$ from $\mathsf{CT}_{\ast }[\mathsf{%
ZF}]$: Let $\alpha $ be the first ordinal with the property that $\left( 
\mathrm{L}_{\alpha },\in \right) $ is a model of $\mathsf{ZF}$ (where $%
\mathrm{L}_{\alpha }$ is the $\alpha $-th approximation to the constructible
universe). Thus $\left( \mathrm{L}_{\alpha },\in \right) $ is the venerable
Shepherdson-Cohen minimal model of $\mathsf{ZF}$. It is well-known that this
model is pointwise definable. Let $T$ be the elementary diagram of $\left( 
\mathrm{L}_{\alpha },\in \right) .$ Then $\left( \mathrm{L}_{\alpha },\in
,T\right) $ is a model of $\mathsf{CT}_{\ast }[\mathsf{ZF}]$ in which $%
\mathrm{Th}\mathsf{(}\mathrm{V}\mathsf{)}\in \mathrm{V}$ fails. This follows
from pointwise definability of $\left( \mathrm{L\!}_{\alpha },\in \right) $
together with Undefinability of Truth Theorem.\medskip
\end{enumerate}

\noindent \textbf{9.7.~Theorem.~}$\mathsf{CT}_{\ast }[\mathsf{ZF}]+\Delta
_{0}$-$\mathsf{Sep}(\mathsf{T})$ \textit{proves that} $\mathsf{ZF}$ \textit{%
has a well-founded model }(\textit{and thus, by Mostowski collapse, }$%
\mathsf{ZF}$\textit{\ has a transitive model})\textit{. }\medskip

\noindent \textbf{Proof.~}We reason in an arbitrary model $(\mathcal{M},%
\mathsf{T})$ of $\mathsf{CT}_{\ast }[\mathsf{ZF}]+\Delta _{0}$-$\mathsf{Sep}(%
\mathsf{T})$. By $\Delta _{0}$-$\mathsf{Sep}(\mathsf{T)}$, there is a
(countable) element $t$ of $\mathcal{M}$ such that: 
\begin{equation*}
(\mathcal{M},\mathsf{T})\models \left[ t=\mathrm{Th}\mathsf{(}\mathrm{V,\in }%
\mathsf{)}\right] .
\end{equation*}%
In light of the assumption that $\mathsf{GRef}_{\mathsf{ZF}}$ holds in $(%
\mathcal{M},\mathsf{T})$, $\mathcal{M}$ satisfies:\medskip

\begin{center}
$t$ is a complete consistent extension of $\mathsf{ZF}$. \medskip
\end{center}

\noindent Now, within $\mathcal{M}$, we can use the Omitting Types Theorem
to build a countable \textit{Paris} model $\mathcal{M}_{0}$ of $t$, i.e., a
model $\mathcal{M}_{0}$ of $t$ such that every ordinal of $\mathcal{M}_{0}$
is pointwise definable from the point of view of $\mathcal{M}$; see \cite%
{Ali-Paris} . We will show that $\mathcal{M}_{0}$ is well-founded from the
point of view of $\mathcal{M}$ using a proof by contradiction. If $\mathcal{M%
}_{0}$ is ill-founded from the point of view of $\mathcal{M}$, then there is
a function $f$ in $\mathcal{M}$ such that:%
\begin{equation*}
\mathcal{M}\models \left[ \left( f:\omega \rightarrow M_{0}\right) \wedge
\forall k\in \omega \ \left( f(k+1)\in ^{\mathcal{M}_{0}}f(k)\right) \right]
.
\end{equation*}%
Here we don't need dependent choice in $\mathcal{M}$ to get hold of $f$,
since $\mathcal{M}_{0}$ is countable in $\mathcal{M}$ and is therefore
well-orderable in $\mathcal{M}$. Let $g(k)=\rho ^{\mathcal{M}_{0}}(f(k))$, \
where $\rho $ is the usual ordinal-valued rank function on sets (as in part
(d) of Definition 2.2). Then:%
\begin{equation*}
\mathcal{M}\models \left[ g:\omega \rightarrow \mathrm{Ord}^{\mathcal{M}%
_{0}}\wedge \forall k\in \omega \ \left( g(k+1)\in ^{\mathcal{M}%
_{0}}g(k)\right) \right] .
\end{equation*}%
Arguing within $\mathcal{M}$, since $\mathcal{M}_{0}$ is a Paris model, the
existence of the function $g$ above makes it clear there is a sequence of
formulae $\left\langle \varphi _{k}(x):k\in \omega \right\rangle $ such
that, for all $k\in \omega $, $t$ includes sentences of the following form: 
\begin{equation*}
\exists !x\,\varphi _{k}(x)\wedge \exists !y\,\varphi _{k+1}(y)\wedge \left(
y\in x\right) .
\end{equation*}%
Let $\alpha _{0}\in \mathrm{Ord}^{\mathcal{M}}$ such that $(\mathcal{M}%
,T)\models \mathsf{T}(\varphi _{0}(\dot{\alpha}_{0}$$)).$ Since $(\mathcal{M}%
,T)\models \Delta _{0}$-$\mathsf{Sep}(\mathsf{T})$, there is some set $s$ in 
$\mathcal{M}$ such that $(\mathcal{M},T)$ satisfies:%
\begin{equation*}
(\mathcal{M},T)\models \left[ s=\left\{ \alpha \in \alpha _{0}:\exists k\in
\omega \,\mathsf{T}(\varphi _{k}(\dot{\alpha}))\right\} \right] .
\end{equation*}%
\noindent It is evident that $\mathcal{M}$ views $s$ as having no $\in $%
-minimal element, which contradicts the axiom of foundation in $\mathcal{M}$.

\hfill $\square $\medskip

\noindent \textbf{9.8.~Remark.~}Consider the sequence of theories $\mathsf{U}%
_{n}$ defined as follows: \medskip

\begin{center}
$\mathsf{U}_{1}:=\mathsf{CT}_{\ast }[\mathsf{ZF}]+\left[ \mathrm{Th}(\mathrm{%
V,\in })\in \mathrm{V}\right] $,\ \ $\mathsf{U}_{2}:=\mathsf{U}_{1}+\left[ 
\mathrm{Th}(\mathrm{V},\in ,\mathrm{Th}(\mathrm{V}))\in \mathrm{V}\right] 
\mathrm{,}$ etc. \medskip
\end{center}

\noindent Thus \textsf{U}$_{1}$ includes the axiom stating that the set of
true sentences (with no constants) exists as a set; and $\mathsf{U}_{2}$
includes the axiom stating that the set of true sentences with at most one
constant naming $\mathrm{Th}(\mathrm{V,\in })$ exists as a set. Then for all 
$n\in \mathbb{N}$ we have:\medskip

\begin{center}
$\mathsf{CT}_{\ast }[\mathsf{ZF}]+\Delta _{0}$-$\mathsf{Sep}(\mathsf{T}%
)\vdash \mathsf{U}_{n}.$ \medskip
\end{center}

\noindent Moreover, using the proof technique of Theorem 9.7, we can show
that the existence of a transitive model of each $\mathsf{U}_{n}$ is
provable in $\mathsf{CT}_{\ast }[\mathsf{ZF}]+\Delta _{0}$-$\mathsf{Sep}(%
\mathsf{T}).\mathsf{\ }$This shows that there is a natural hierarchy of
theories between $\mathsf{CT}_{\ast }[\mathsf{ZF}]$ and $\mathsf{CT}_{\ast }[%
\mathsf{ZF}]+\Delta _{0}$-$\mathsf{Sep}(\mathsf{T}).$\medskip 

\noindent \textbf{9.9.~Theorem.~}\textit{The existence of an }$\omega $%
\textit{-model of} $\mathsf{ZF}$ \textit{is provable in} $\mathsf{CT}_{\ast
}[\mathsf{ZF}]+\left[ \mathrm{Th}(\mathrm{V,\in })\in \mathrm{V}\right] 
\mathrm{.}$ \medskip

\noindent \textbf{Proof.~}This can be established with the proof strategy of
Theorem 9.7, together with the fact that by Theorems 6.10 and 6.12, $\Delta
_{0}^{\mathrm{fin}}$-$\mathsf{Sep}(\mathsf{T})$ is provable in $\mathsf{CT}%
_{\ast }[\mathsf{ZF}].$ The role of the foundation axiom in the proof of
Theorem 9.7 is replaced here with the Pigeonhole Principle, in the basic
form that asserts that given any $k\in \omega $, there is no injection from
the elements of $k+1$ to the elements of $k$.

\hfill $\square $\medskip

\noindent \textbf{9.10.~Theorem.~}\textit{For recursively axiomatized
theories} $\mathsf{U}$ and $\mathsf{V}$ \textit{including }$\mathsf{KP}$,%
\textbf{\ }\textit{let} $\mathsf{U}\blacktriangleleft \mathsf{V}$ \textit{%
stand for the conjunction of} $\mathsf{V}\vdash \mathsf{U}$ \textit{and }$%
\mathsf{V}\vdash \mathrm{Con}(\mathsf{U}).$ \textit{Then we have}:\medskip

\begin{center}
$\mathsf{CT}_{\ast }[\mathsf{ZF}]\blacktriangleleft \mathsf{CT}_{\ast }[%
\mathsf{ZF}]+\Delta _{0}$-$\mathsf{Sep}(\mathsf{T})\blacktriangleleft 
\mathsf{CT}_{\ast }[\mathsf{ZF}]+\mathsf{FRef\blacktriangleleft CT}_{0}[%
\mathsf{ZF}].$\medskip
\end{center}

\noindent \textit{More explicitly}:\textbf{\ }\medskip

\begin{enumerate}
\item[$(a)$] $\mathsf{CT}_{0}[\mathsf{ZF}]$ \textit{proves that} $\mathsf{CT}%
_{\ast }[\mathsf{ZF}]+\mathsf{FRef}$ \textit{has a model of the form} $(%
\mathrm{V\!}_{\alpha },\in ,T).$ \textit{In particular,\smallskip } 
\begin{equation*}
\mathsf{CT}_{0}[\mathsf{ZF}]\vdash \mathrm{Con}\left( \mathsf{CT}_{\ast }[%
\mathsf{ZF}]+\mathsf{FRef}+\mathsf{Sep(T)}\right) .
\end{equation*}

\item[$(b)$] $\mathsf{CT}_{\ast }[\mathsf{ZF}]+\mathsf{FRef}$ \textit{proves
that} $\mathsf{CT}_{\ast }[\mathsf{ZF}]+\Delta _{0}$-$\mathsf{Sep}(\mathsf{T}%
)$ \textit{has a model of the form} $(\mathrm{V\!}_{\alpha },\in ,T).$ 
\textit{In particular,\smallskip } 
\begin{equation*}
\mathsf{CT}_{\ast }[\mathsf{ZF}]+\mathsf{FRef}\vdash \mathrm{Con}\left( 
\mathsf{CT}_{\ast }[\mathsf{ZF}])+\mathsf{Sep(T)}\right) .
\end{equation*}

\item[$(c)$] $\mathsf{CT}_{\ast }[\mathsf{ZF}]+\Delta _{0}$-$\mathsf{Sep}(%
\mathsf{T})$ \textit{proves that} $\mathsf{CT}_{\ast }[\mathsf{ZF}]$ \textit{%
has a model of the form }$(m,\in T)$ \textit{for some }$m.$ \textit{In
particular,}\smallskip
\end{enumerate}

\begin{center}
$\mathsf{CT}[\mathsf{ZF}]+\Delta _{0}$-$\mathsf{Sep}(\mathsf{T})\vdash 
\mathrm{Con}\left( \mathsf{CT}_{\ast }[\mathsf{ZF}])+\mathsf{Ind}(\mathsf{T}%
)\right) .$\smallskip
\end{center}

\begin{enumerate}
\item[$(d)$] \textit{The existence of an} $\omega $-\textit{model of }$%
\mathsf{ZF}$ \textit{is not provable in} $\mathsf{CT}_{\ast }[\mathsf{ZF}]$, 
\textit{but\smallskip\ }%
\begin{equation*}
\mathsf{CT}_{\ast }[\mathsf{ZF}]\vdash \mathrm{Con}\left( \mathsf{CT}^{-}[%
\mathsf{ZF}]+\mathsf{Sep}(\mathsf{T})\right) .
\end{equation*}
\end{enumerate}

\noindent \textbf{Proof.~}$(a)$ follows from the provability of $\mathsf{FRef%
}^{2}$ in $\mathsf{CT}_{0}[\mathsf{ZF}]$, established in Theorem 5.9. To
show $(b)$, we reason in $\mathsf{CT}_{\ast }[\mathsf{ZF}]+\mathsf{FRef}.$
By $\mathsf{FRef}$ there is an ordinal $\alpha $ such that $\left( \mathrm{%
V\!}_{\alpha },\in \right) \models \mathsf{ZF.}$ Given such an $\alpha $,
let $T$ be the elementary diagram for $\left( \mathrm{V\!}_{\alpha },\in
\right) .$ By Remark 9.2, $\left( \mathrm{V\!}_{\alpha },\in ,T\right)
\models \mathsf{CT}_{\ast }[\mathsf{ZF}]$. So the proof of (b) is complete
once we observe that (provably in $\mathsf{ZF}$) if $X\subseteq \mathrm{V\!}%
_{\alpha }$, then we have:%
\begin{equation*}
\forall a\in \mathrm{V\!}_{\alpha }\ X\cap a\in \mathrm{V\!}_{\alpha }.
\end{equation*}%
\noindent To verify $(c)$, we reason in $\mathsf{CT}_{\ast }[\mathsf{ZF}%
]+\Delta _{0}$-$\mathsf{Sep}(\mathsf{T}).$ By Theorem 9.7, there is some $%
m\in M$ such that $(m,\in )$ is a model of $\mathsf{ZF}$. Let $T$ be the
elementary diagram for $\left( m,\in \right) .$ In light of Remark 9.2, $%
\left( m,\in ,T\right) \models \mathsf{CT}_{\ast }[\mathsf{ZF}]+\mathsf{Ind}(%
\mathsf{T})$.\medskip 

\noindent The first assertion of $(d)$ follows from Proposition 9.3. Since
the consistency of $\mathsf{ZF}$ is provable in \textit{\ }$\mathsf{CT}%
_{\ast }[\mathsf{ZF}]$, the second assertion of (d) follows from the fact
that Theorem 4.3 is provable in $\mathsf{ZF}$.

\hfill $\square $\medskip

\noindent \textbf{9.11.~Remark.~}Each of the theories $\mathsf{CT}^{-}[%
\mathsf{ZF}]+\mathsf{Sep(T)}$, and $\mathsf{CT}^{-}[\mathsf{ZF}]~+$ $\mathsf{%
Int}$\textsf{-}$\mathsf{Repl}$\ is conservative over $\mathsf{ZF}$, but
their union is not, as it implies $\mathsf{CT}_{\ast }[\mathsf{ZF}]+\Delta
_{0}$-$\mathsf{Sep}(\mathsf{T})$. Note that by part (b) of Theorem 9.10, $%
\mathsf{CT}^{-}[\mathsf{ZF}]+\mathsf{FRef}$ proves the consistency of $%
\mathsf{CT}^{-}[\mathsf{ZF}]+\mathsf{Sep(T)}+\mathsf{Int}$\textsf{-}$\mathsf{%
Repl.}$\bigskip \bigskip

\begin{center}
\textbf{10.~CONSERVATIVITY OF }$\mathsf{CT}^{-}[\mathsf{ZF}]+\mathsf{Coll(T)}
$\textbf{\ OVER }$\mathsf{ZF}$\bigskip
\end{center}

\noindent In this section we establish the conservativity of $\mathsf{CT}%
^{-}[\mathsf{ZF}]+\mathsf{Coll(T)}$\textbf{\ }over\textbf{\ }$\mathsf{ZF}$
in Theorem 10.2. This result complements the fact that $\mathsf{CT}^{-}[%
\mathsf{ZF}]+\mathsf{Sep(T)}$ is conservative over $\mathsf{ZF}$ (Theorem
4.3). The proof of Theorem 10.2 is based on combining a key result due to
Keisler on elementary end extensions of models of $\mathsf{ZF}$ with the
model-theoretic method introduced in \cite{Ali+albert-short} for the
construction of full satisfaction classes. The proof strategy was inspired
by Wcis\l o's proof \cite{Wcislo-Collection} of the conservativity of $%
\mathsf{CT}^{-}[\mathsf{PA}]+\mathsf{Coll(T)}$\textbf{\ }over\textbf{\ }$%
\mathsf{PA}$. We begin by reviewing Keisler's theorem.\footnote{%
The model theory of the collection scheme (in the guise of the regularity
scheme) is studied in \cite{Enayat+Mohsenipour}.} \medskip

\noindent \textbf{10.1.~Theorem} (Keisler).\textbf{~}\textit{Let }$\mathcal{N%
}$ \textit{be a countable model of} $\mathsf{ZF}(\mathcal{L})$ \textit{for
some countable }$\mathcal{L}\supseteq \mathcal{L}_{\mathrm{set}}$. \textit{%
Then the following hold}:\medskip

\begin{enumerate}
\item[$(a)$] $\mathcal{N}$ \textit{has a countable elementary end extension}%
. \medskip

\item[$(b)$] $\mathcal{N}$ \textit{has an} $\aleph _{1}$-\textit{like
elementary end extension}.\footnote{$\mathcal{M}$ is said to be $\aleph _{1}$%
-like if the universe $M$ of $\mathcal{M}$ has cardinality $\aleph _{1}$,
but for each $m\in M,$ $\left\{ x\in M:\mathcal{M}\models x\in m\right\} $
is countable.}\medskip
\end{enumerate}

\noindent \textbf{Proof outline}.\textbf{~}The Keisler-Morley Theorem as
proved in \cite{Keisler-Morley} has a number of versions; the one that is
usually stated concerns elementary end extensions of countable models of $%
\mathsf{ZF}$, e.g., as in the exposition in the Chang-Keisler canonical
reference in Model Theory \cite[Theorem 2.2.18]{Chang-Keisler}. The omitting
types proof in the aforementioned reference (which takes advantage of the
provability of the collection scheme in $\mathsf{ZF}$ in the guise of the
regularity scheme) shows the more general result that if $\mathcal{L}$\ is a
countable language extending $\mathcal{L}_{\mathrm{set}}$, then every
countable model $\mathcal{N}$ of $\mathsf{ZF}(\mathcal{L})$ has an
elementary end extension. As a consequence, by using this theorem $\aleph
_{1}$-times (while taking unions at limit ordinals) $\mathcal{N}$ has an $%
\aleph _{1}$-like elementary end extension.

\hfill $\square $\medskip

\noindent \textbf{10.2.~Theorem.~}$\mathsf{CT}_{\ast }[\mathsf{ZF}]+\mathsf{%
Int}$-$\mathsf{Repl+Coll}(\mathsf{T})$ \textit{is conservative over} $%
\mathsf{ZF}$.\footnote{%
The meta-theory for this result is $\mathsf{ZFC}$. It can be reduced to $%
\mathsf{Z}_{3}$ (third order arithmetic) plus enough choice to guarantee
that $\aleph _{1}$ is a regular cardinal. This result can be further
strengthened by adding further `good behavior' axioms to $\mathsf{CT}^{-}[%
\mathsf{ZF}]+\mathsf{Int}$-$\mathsf{Repl+Coll}(\mathsf{T})$, such as $%
\mathsf{EC}$ (existential correctness, see the proof of $(a)\Rightarrow (b)$
of Theorem 6.10), and the agreement of $\mathsf{T}$ with definable partial
truth predicates.}\medskip

\noindent \textbf{Proof}.\textbf{~}Since $\aleph _{1}$ is a regular cardinal
(in $\mathsf{ZFC}$), it is a theorem of $\mathsf{ZFC}$ that every $\aleph
_{1}$-like $\mathcal{L}$-structure satisfies $\mathsf{Coll}(\mathcal{L}).$
Also recall that if $S$ is a full extensional satisfaction class over a
model $\mathcal{M}$, then the associated truth predicate $T_{S}$ (as in
Proposition 3.2.6) is a full truth class on $\mathcal{M}$. In light of the
completeness theorem of first order logic, this makes it evident that the
proof of Theorem 10.2 is complete once we show that every countable model $%
\mathcal{M}_{0}\models \mathsf{ZF}$ has an elementary extension $\mathcal{M}%
^{\ast }$ that has an expansion $\left( \mathcal{M}^{\ast },S^{\ast }\right) 
$ satisfying the following three properties:\medskip

\noindent (1) $S^{\ast }$ is a full extensional satisfaction class over $%
\mathcal{M}^{\ast }$.\medskip

\noindent (2) $\mathsf{Int}$-$\mathsf{Repl}$ is deemed true by $S^{\ast }$%
.\medskip

\noindent (3) $\mathcal{M}^{\ast }$ is $\aleph _{1}$-like.\medskip

\noindent To construct the desired $\mathcal{M}^{\ast }$ satisfying (1)
through (3) above we argue as follows:\medskip

\noindent STEP 1\textbf{.} By Proposition 3.2.6 and Corollary 4.7, there is
a countable elementary extension $\mathcal{M}_{1}$ of $\mathcal{M}_{0}$ such
that $\mathcal{M}_{1}$ has an expansion $\left( \mathcal{M}_{1},S_{1}\right) 
$ such that $S_{1}$ is a full extensional satisfaction class over $\mathcal{M%
}_{1}$ and $\left( \mathcal{M}_{1},S_{1}\right) $ satisfies $\mathsf{Int}$-$%
\mathsf{Repl}$.\medskip

\noindent STEP 2. Let $F_{1}=\mathrm{Form}^{\mathcal{M}_{1}}$, and for each $%
\varphi \in F_{1}$, let $X_{\varphi }=\{\alpha \in \mathrm{Asn}^{\mathcal{M}%
_{1}}:\left\langle \varphi ,\alpha \right\rangle \in S\}.$ Note that: 
\begin{equation*}
S_{1}=\left\{ \left\langle \varphi ,\alpha \right\rangle :\varphi \in
F_{1},\ \alpha \in X_{\varphi }\right\} .
\end{equation*}%
Since internal replacement holds in $\left( \mathcal{M}_{1},S_{1}\right) $,
by Theorem 7.9 $(\mathcal{M}_{1},X_{\varphi })_{\varphi \in F_{1}}\models 
\mathsf{ZF}(\mathcal{L})$, where $\mathcal{L}$ is the result of augmenting $%
\mathcal{L}_{\mathrm{set}}$ with predicate symbols $\mathrm{X}_{\varphi }$
for each $\varphi \in F_{1}\mathfrak{.}$ \medskip 

\noindent STEP 3. The countability of both $\mathcal{M}_{1}$ and $\mathfrak{X%
}$, together with the fact that $(\mathcal{M}_{1},X_{\varphi })_{\varphi \in
F_{1}}\models \mathsf{ZF}(\mathcal{L})$ allows us to invoke Theorem 10.1 to
get hold of $(\mathcal{M}^{\ast },X_{\varphi }^{\ast })_{\varphi \in F_{1}}$
such that:\medskip

\noindent $(i)$ $\ \ $ $\mathcal{M}^{\ast }$ is $\aleph _{1}$-like and $(%
\mathcal{M}_{1},X_{\varphi })_{\varphi \in F_{1}}\prec _{\mathrm{end}}(%
\mathcal{M}^{\ast },X_{\varphi }^{\ast })_{\varphi \in F_{1}}$.\medskip

\noindent The fact that $\mathcal{M}^{\ast }$ is an end extension of $%
\mathcal{M}_{1}$ assures us that $\omega ^{\mathcal{M}_{1}}$ is not enlarged
in the passage from $\mathcal{M}_{1}$ to $\mathcal{M}^{\ast }$, i.e.,
\medskip 

\noindent $(ii)$ $\ \ \left\{ x\in M_{1}:\mathcal{M}_{1}\models \left[ x\in
\omega \right] \right\} =\left\{ x\in M^{\ast }:\mathcal{M}^{\ast }\models %
\left[ x\in \omega \right] \right\} ,$\medskip

\noindent which in turn assures us that:\medskip

\noindent $(iii)$ $\ \ \mathcal{M}_{1}$ and $\mathcal{M}^{\ast }$ have the
same $\mathcal{L}_{\mathrm{set}}$-formulae, i.e., 
\begin{equation*}
\left\{ x\in M_{1}:\mathcal{M}_{1}\models \left[ x\in \mathrm{Form}\right]
\right\} =\left\{ x\in M^{\ast }:\mathcal{M}^{\ast }\models \left[ x\in 
\mathrm{Form}\right] \right\} .
\end{equation*}%
However, $\mathcal{M}^{\ast }$ has far more assignments than $\mathcal{M}_{1}
$.\medskip 

\noindent STEP 4. The fact that $S_{1}$ is an extensional satisfaction class
on $\mathcal{M}_{1}$ that satisfies $\mathsf{Int}$-$\mathsf{Repl}$ makes it
clear that $(\mathcal{M}_{1},X_{\varphi })_{\varphi \in F_{1}}$ satisfies
the \textit{universal generalizations} of the statements (1) through (7)
below, in which $\varphi ,$ $\psi ,$ $\psi _{1},$ $\psi _{2}$ vary over
(codes of) $\mathcal{L}_{\mathrm{set}}$-formulae, and $\alpha $ varies over
assignments.\footnote{%
See Definition 3.2.1 for the abbreviations used in (1) through (7).}\medskip

\begin{enumerate}
\item[(1)] $\left[ \varphi =\ulcorner x=y\urcorner \right] \longrightarrow %
\left[ \mathsf{X}_{\varphi }(\alpha )\leftrightarrow \left[ \mathrm{Asn}%
(\alpha ,\varphi )\wedge \alpha (x)=\alpha (y)\right] \right] $.\medskip

\item[(2)] $\left[ \varphi =\ulcorner x\in y\urcorner \right]
\longrightarrow \left[ \mathsf{X}_{\varphi }(\alpha )\leftrightarrow \left[ 
\mathrm{Asn}(\alpha ,\varphi )\wedge \alpha (x)\in \alpha (y)\right] \right]
.$\medskip

\item[(3)] $\left[ \varphi =\lnot \psi \right] \longrightarrow \left[ 
\mathsf{X}_{\varphi }(\alpha )\leftrightarrow \left[ \mathrm{Asn}(\alpha
,\varphi )\wedge \lnot \mathsf{X}_{\psi }(\alpha )\right] \right] .$\medskip

\item[(4)] $\left[ \varphi =\left( \psi _{1}\vee \psi _{2}\right) \right]
\longrightarrow $ $\left[ \mathsf{X}_{\varphi }(\alpha )\leftrightarrow %
\left[ \mathrm{Asn}(\alpha ,\varphi )\wedge \left( \mathsf{X}_{\psi
_{1}}(\alpha \upharpoonright \mathrm{FV}(\psi _{1}))\vee \mathsf{X}_{\psi
_{2}}(\alpha \upharpoonright \mathrm{FV}(\psi _{2}))\right) \right] \right]
. $\medskip

\item[(5)] $\left[ \varphi =\exists v\,\psi \right] \longrightarrow \left[ 
\mathsf{X}_{\varphi }(\alpha )\leftrightarrow \left[ \mathrm{Asn}(\alpha
,\varphi )\wedge \exists \beta \supseteq \alpha \,\mathsf{X}_{\psi }(\beta
)\wedge \mathrm{Asn}(\beta ,\psi )\right] \right] .$\medskip

\item[(6)] $\left[ \psi =\mathrm{Repl}_{\varphi (v,x,y)}\right] \rightarrow %
\left[ \forall \alpha \left( \mathrm{Asn}(\alpha ,\psi )\rightarrow \mathsf{X%
}_{\varphi }(\alpha )\right) \right] .$\footnote{$\mathsf{Repl}_{\varphi
(v,x,y)}$ was defined in part (d) of Definition 2.1. Note that (6) ensures
that $\forall v\,\mathsf{Repl}_{\varphi (v,x,y)}$ is deemed true by the
satisfaction predicate $S^{\ast }$ described in Step 5 of the proof.}\medskip

\item[(7)] $\left[ (\varphi _{0},\alpha _{0})\thicksim (\varphi _{1},\alpha
_{1})\right] \longrightarrow \left[ \mathsf{X}_{\varphi _{0}}(\alpha
_{0})\leftrightarrow \mathsf{X}_{\varphi _{1}}(\alpha _{1})\right] .$\medskip
\end{enumerate}

\noindent STEP 5. Since $(\mathcal{M}_{1},X_{\varphi })_{\varphi \in
F_{1}}\prec (\mathcal{M}^{\ast },X_{\varphi }^{\ast })_{\varphi \in F_{1}}$
by $(i)$ of Step 3, we can conclude that $(\mathcal{M}^{\ast },X_{\varphi
}^{\ast })_{\varphi \in F_{1}}$ also satisfies conditions (1) through (7).
This fact makes it clear that by `gluing' the family $\left\{ X_{\varphi
}^{\ast }:\varphi \in F_{1}\right\} $ together as:%
\begin{equation*}
S^{\ast }=\left\{ \left\langle \varphi ,\alpha \right\rangle :\varphi \in
F_{1},\ \alpha \in X_{\varphi }^{\ast }\right\} ,
\end{equation*}
we obtain an extensional full satisfaction class $S^{\ast }$ on $\mathcal{M}%
^{\ast }$ that validates internal replacement. Note that the `fullness' of $%
S^{\ast }$ is assured by $(iii)$ of Step 3. In light of Proposition 3.2.6,
if $T^{\ast }:=\mathcal{T}(S^{\ast })$ is the truth class corresponding to $%
S^{\ast }$, then we have $(\mathcal{M}^{\ast },T^{\ast })\models \mathsf{CT}%
^{-}[\mathsf{ZF}]+\mathsf{Int}$-$\mathsf{Repl.}$

\hfill $\square $\bigskip \bigskip

\begin{center}
\textbf{11.~QUESTIONS}\bigskip
\end{center}

\noindent The following two questions are motivated by the model-theoretic
proof of Theorem 8.8, which characterizes purely set-theoretical
consequences of $\mathsf{CT}_{\ast }[\mathsf{ZFC}].$ Note that the proof of
Theorem 8.8 breaks down when applied to $\mathsf{CT}_{\ast }[\mathsf{ZF}]$.
See also Remark 8.9 and 8.10.\medskip

\noindent \textbf{11.1.~Question.~}\textit{Does }$\mathsf{CON}^{\omega }%
\mathsf{(ZF)}$ \textit{axiomatize} \textit{the set of purely set-theoretical
consequences of }$\mathsf{CT}_{\ast }[\mathsf{ZF}]?$\medskip

\noindent \textbf{11.2.~Question.~}\textit{Can Theorem 8.8 also be
demonstrated by proof-theoretic methods?}\medskip

\noindent The following question is motivated by the provability of
\textquotedblleft $\mathsf{ZF}$ has a well-founded model\textquotedblright\
in $\mathsf{CT}_{\ast }[\mathsf{ZF}]+\Delta _{0}$-$\mathsf{Sep(T)}$,
established in Theorem 9.7, and the provability of \textquotedblleft $%
\mathsf{ZF}$ has an $\omega $-model\textquotedblright\ in $\mathsf{CT}_{\ast
}[\mathsf{ZF}]+\left[ \mathrm{Th}(\mathrm{V,\in })\in \mathrm{V}\right] $,
established in Theorem 9.9. \medskip

\noindent \textbf{11.3.~Question.~}\textit{Can the theory }$\mathsf{CT}%
_{\ast }[\mathsf{ZF}]+\left[ \mathrm{Th}(\mathrm{V,\in })\in \mathrm{V}%
\right] $ \textit{prove the} \textit{existence of a well-founded model of }$%
\mathsf{ZF}$\textit{?}\medskip

\noindent The conservativity of $\mathsf{CT}^{-}[\mathsf{ZF}]+\mathsf{Coll(T)%
}$ over $\mathsf{ZF}$ (Theorem 10.2) was established by a model-theoretic
argument involving uncountable models. The highly nonfinitary nature of the
proof motivates the following questions. See also Theorem 4.9 and Remark
4.10.\medskip

\noindent \textbf{11.4.~Question.~}\textit{Is} \textit{the conservativity of}
$\mathsf{CT}^{-}[\mathsf{ZF}]+\mathsf{Coll(T)}$ \textit{over} $\mathsf{ZF}$ 
\textit{provable in} $\mathsf{PA?}$\medskip

\noindent \textbf{11.5.~Question.~}\textit{Is} $\mathsf{CT}^{-}[\mathsf{ZF}]+%
\mathsf{Coll(T)}$ \textit{interpretable in }$\mathsf{ZF?}$\medskip

\noindent \textbf{11.6.~Question.~}\textit{Does }$\mathsf{CT}^{-}[\mathsf{ZF}%
]+\mathsf{Coll(T)}$ \textit{exhibit superpolynomial speed-up over} $\mathsf{%
ZF}$?

\noindent \textsc{Ali Enayat, Department of Philosophy, Linguistics, and
Theory of Science, University of Gothenburg}, \textsc{\ Sweden; }\texttt{%
email: ali.enayat@gu.se}

\end{document}